\documentclass{article}

\usepackage {amsmath, amssymb, amsthm}
\usepackage {mathrsfs}
\usepackage [top=1.5in, right=1in, bottom=1.5in, left=1in] {geometry}
\usepackage {graphicx}
\usepackage {tikz}
\usetikzlibrary {snakes}
\usepackage {hyperref}
\hypersetup {
    colorlinks = true,
    linktoc = all,
    linkcolor = blue,
    citecolor = red
}

\usepackage [nottoc, notlot, notlof] {tocbibind}

\usepackage {libertine}
\usepackage [T1]{fontenc}

\newcommand {\db} [1] {\left\{\!\!\left\{ #1 \right\}\!\!\right\}}

\newcommand {\freeprod} {\mathop{\hbox{\scalebox{1.7}{$\ast$}}}}

\theoremstyle {plain}

\newtheorem {lem} {Lemma}
\newtheorem {thm} {Theorem}
\newtheorem {cor} {Corollary}
\newtheorem {prop} {Proposition}

\theoremstyle {definition}
\newtheorem* {ex} {Example}

\newtheorem {defn} {Definition}

\newtheorem* {rmk} {Remark}

\title{Non-Commutative Integrability \\of the \\Grassmann Pentagram Map}
\author{Nicholas Ovenhouse \thanks{Department of Mathematics, Michigan State University}}
\date{}

\begin {document}

\maketitle

\begin {abstract}
    The pentagram map is a discrete integrable system first introduced by Schwartz in 1992.
    It was proved to be intregable by Schwartz, Ovsienko, and Tabachnikov in 2010.
    Gekhtman, Shapiro, and Vainshtein studied Poisson geometry associated to certain networks
    embedded in a disc or annulus, and its relation to cluster algebras. Later, Gekhtman et al.
    and Tabachnikov reinterpreted the pentagram map in terms of these networks, and used the
    associated Poisson structures to give a new proof of integrability. In 2011, Mari Beffa and Felipe 
    introduced a generalization of the pentagram map to certain Grassmannians, and proved
    it had a Lax representation. We reinterpret this Grassmann pentagram map in terms of noncommutative algebra,
    in particular the double brackets of Van den bergh,
    and generalize the approach of Gekhtman et al. to establish a noncommutative version
    of integrability.
\end {abstract}

\tableofcontents

\section {Introduction}

\bigskip

\subsection {Background}

\bigskip

The pentagram map was first introduced by Schwartz \cite{schwartz92} as a transformation of
the moduli space of labeled polygons in $\Bbb{P}^2$. After labeling the vertices of an $n$-gon
with the numbers $1,\dots,n$, the diagonals are drawn which connects each $i$ to $i+2$ 
(see \textbf{Figure \ref{fig:hex_map}}).
The intersection points of these diagonals are taken to be the vertices of a new $n$-gon.
Later, Schwartz, along with Ovsienko and Tabachnikov, generalized the map to the space of
``twisted'' $n$-gons in $\Bbb{P}^2$, and proved that this generalized pentagram map is
completely integrable in the Liouville sense \cite{ost10}. In 2011, Max Glick interpreted the pentagram
map in terms of ``$Y$-mutations'' of a certain cluster algebra \cite{glick11}. In the literature,
there are other names for $Y$-mutations or $Y$-dynamics. Fomin and Zelevinsky have called
them ``coefficient dynamics'' \cite{fz1} \cite{fz4}. Gekhtman, Shapiro, and Vainshtein
refer to them as ``$\tau$-coordinate mutations'' \cite{gsv_book}. Fock and Goncharov refer
to them as ``$\mathcal{X}$-variable mutations'' \cite{fock_goncharov}.

\bigskip

In 2016, Gekhtman, Shapiro, Vainshtein, and Tabachnikov, building on the work of both Glick and Postnikov \cite{postnikov_06},
interpreted a certain set of coordinates on the space of twisted polygons as weights on some
directed graph, and the pentagram map in these coordinates as a certain sequence of ``Postnikov moves''
applied to this graph \cite{gsvt_16}. They used their previous work on the Poisson geometry of the space of edge weights
of such graphs \cite{gsv_09} \cite{gsv_10} to give a new proof of the integrability of the pentagram map.

\bigskip

Recently, Mari Beffa and Felipe considered a generalization of the pentagram map, where the twisted polygons were taken to
be in the Grassmannian instead of projective space \cite{bf15}, and they demonstrated a Lax representation for
this version of the pentagram map, establishing integrability for this generalized version. 
The purpose of the present paper is to interpret this Grassmannian pentagram map
as a transformation of a set of matrix-valued variables (and more generally as a formal noncommutative rational transformation), 
and use ideas from non-commutative Poisson geometry
-- namely the ``double brackets'' of Van den Bergh \cite{bergh_08} and ``$H_0$-Poisson structures'' of Crawley-Boevey\cite{crawley-boevey_11} -- 
to formulate a non-commutative version of
integrability, generalizing the approach of Gekhtman, Shapiro, Vainshtein, and Tabachnikov, using weighted directed graphs.
This non-commutative algebraic structure is a generalization of both the original pentagram map and this more recent
Grassmann version. Furthermore, the integrability of both examples follows from this non-commutative structure,
by projecting to representation spaces of a certain non-commutative algebra.

\bigskip

\subsection {Structure of the Paper}

\bigskip

The structure of the paper is as follows. The first three sections are a condensed review of the basics of the
pentagram map. We mainly follow the notations and conventions of
Gekhtman, Shapiro, Vainshtein, and Tabachnikov \cite{gsv_09} \cite{gsv_10} \cite{gsvt_16}. The later sections
are an introduction to the Grassmannian pentagram map, and the development of the new theory. We now give a more
detailed description of the sections.

\bigskip

In section 1, we review the definition of ``twisted'' polygons and the
pentagram map. We also construct a system of coordinates
on the moduli space of twisted polygons (following \cite{gsvt_16}), 
and express the pentagram map in these coordinates, where it is seen to be a relatively simple rational expression.

\bigskip

In section 2, we define the weighted directed graphs that we will be using, and review the associated Poisson structure
on the space of edge weights introduced in \cite{gsv_10}. We also review the notion of ``boundary measurements'' for
such graphs, and the ``Postnikov moves'' and ``gauge transformations'' which leave the boundary measurements invariant \cite{postnikov_06}.

\bigskip

In section 3, we review the observation from \cite{gsvt_16} that after identifying the coordinates from section 1 with
certain weights on some graph, the pentagram map can be interpreted as a sequence of Postnikov moves and gauge transformations.
This observation allows us to find invariants of the pentagram map, defined in terms of the boundary measurements.
These invariants also turn out to be in involution with respect to the Poisson structure, giving the integrability
of the pentagram map.

\bigskip

In section 4, we introduce the notion of twisted polygons in the Grassmannian, and define the pentagram map
for these twisted Grassman polygons, following the presentation of Mari Beffa and Felipe \cite{bf15}. 
We then identify the moduli space with a quotient of the space of matrix-valued weights on certain graphs,
and use this identification to coordinatize the moduli space.
Next, we write the pentagram map in terms of
the matrix weights, in an attempt to mimic, formally, the approach of \cite{gsvt_16}.
We show that in these (non-commutative) variables, the pentagram map is a non-commutative 
version of the formula from section 1.

\bigskip

In section 5, we review the definitions and properties of ``double brackets'' on non-commutative algebras
(due to Van den Bergh \cite{bergh_08}), and ``$H_0$-Poisson structures'' (due to Crawley-Boevey \cite{crawley-boevey_11}).
This is the formalism we will use in place of a usual Poisson structure
in this non-commutative setting.

\bigskip

In section 6, we consider the same directed graph from section 2, but we consider non-commutative weights.
We also consider the analogues of the boundary measurements, Postnikov moves,
and gauge transformations in this non-commutative setting. To make rigorous sense of the Postnikov moves,
we introduce the free skew field in a set of formal non-commuting indeterminates.
We then define a double bracket on the free skew field generated by the edges. This algebra acts
as a sort of ``space of non-commutative edge weights''. The relationship with Goldman's bracket
on character varieties \cite{goldman_86} is discussed, giving the double bracket (and its induced brackets)
a geometric interpretation in terms of intersection numbers on a certain surface. Finally, we find a set of
elements/variables in the free skew field which are non-commutative analogues to the coordinates from section 1.
The same sequence of Postnikov moves as in the classical case transforms these variables precisely by
the formula from section 4 for the Grassman pentagram map.

\bigskip

In section 7, analogous to the classical case, we find noncommutative ``invariants'' of the map
which are defined in terms of the boundary measurements. They are in fact invariant only modulo commutators. 
We also show that these invariants are, 
in some sense, involutive with respect to the non-commutative Poisson structure, giving a non-commutative
version of integrability.

\bigskip

In section 8, we discuss how our noncommutative Poisson structures induce usual ``commutative'' Poisson structures
on the moduli space of representations of a certain free algebra. We review a discussion from section 4
which identifies the moduli space of polygons with a subset of this representation space. Thus, the moduli
space of polygons inherits a Poisson structure. The invariants obtained by Mari Beffa and Felipe are the
traces of our non-commutative invariants, and they form a Liouville integrable system under the induced Poisson bracket.

\bigskip

Finally, in section 9, we pose some lingering questions and directions for further research.

\bigskip

\section {The ``Classical'' Pentagram Map}

\bigskip

\subsection {Twisted Polygons and the Pentagram Map}

\bigskip

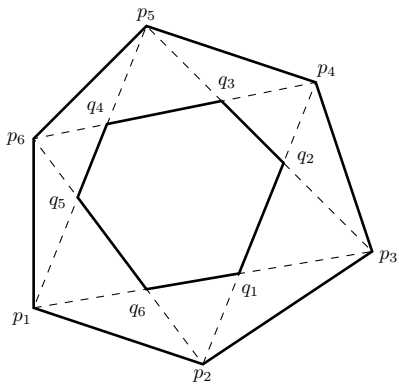
\begin {figure}[h!]
    \centering
    \caption{The pentagram map for a hexagon}
    \label{fig:hex_map}
    \begin {tikzpicture}[scale=0.75, every node/.style={scale=0.75}]
        \draw [line width = 1pt] (1,1) -- (4,0) -- (7,2) -- (6,5) -- (3,6) -- (1,4) -- cycle;
        \draw [dashed] (1,1) -- (7,2) -- (3,6) -- cycle;
        \draw [dashed] (4,0) -- (6,5) -- (1,4) -- cycle;
        \draw [line width = 1pt] (4.64,1.61) -- (5.43,3.57) -- (4.33,4.67) -- (2.3,4.26) -- (1.78,2.96) -- (3,1.33) -- cycle;

        \draw (0.8,0.8) node {$p_1$};
        \draw (4, -0.2) node {$p_2$};
        \draw (7.3, 1.9) node {$p_3$};
        \draw (6.2, 5.2) node {$p_4$};
        \draw (3,6.2) node {$p_5$};
        \draw (0.7,4) node {$p_6$};

        \draw (4.84,1.31) node {$q_1$};
        \draw (5.83,3.67) node {$q_2$};
        \draw (4.43,5) node {$q_3$};
        \draw (2.1,4.56) node {$q_4$};
        \draw (1.43,2.86) node {$q_5$};
        \draw (2.85,0.95) node {$q_6$};
    \end {tikzpicture}
\end {figure}

We refer to the pentagram map as the ``classical'' pentagram map to distinguish it from the Grassmann pentagram map (the subject of the later sections of the paper).
First we give the basic idea, and then extend to so-called ``twisted'' polygons. The map was originally introduced by Schwartz \cite{schwartz92}, but here we mainly
follow the notations and conventions from \cite{gsvt_16}.
Consider an $n$-gon in $\Bbb{P}^2 = \Bbb{RP}^2$, with vertices labeled
$p_1$ through $p_n$. Draw the diagonals of the $n$-gon which connect $p_i$ to $p_{i+2}$ (with indices read cyclically).
Label the intersection of the lines $\overline{p_i p_{i+2}}$ and $\overline{p_{i+1} p_{i+3}}$ as $q_i$. 
Then the pentagram map, which we denote by $T$, sends the first labeled polygon $P$
to the labeled polygon $Q = T(P)$ whose vertices are $q_i$. An example for $n=6$ is shown below in \textbf{Figure \ref{fig:hex_map}}.

More generally, we define ``twisted polygons''. 

\begin {defn}
    A \emph{twisted $n$-gon} in $\Bbb{P}^2$ is a bi-infinte sequence $(p_i)_{i \in \Bbb{Z}}$ of points in $\Bbb{P}^2$
    such that $p_{i+n} = M p_i$ for all $i$, where $M$ is some projective transformation in $\mathrm{PGL}_3(\Bbb{R})$,
    referred to as the \emph{monodromy} of the twisted $n$-gon. 
\end {defn}

\begin {rmk}
    The usual notion of an $n$-gon can be thought of as the case
    where $M = \mathrm{Id}$, and the sequence is periodic. In this case we call the polygon a \emph{closed $n$-gon}. 
\end {rmk}
    
The pentagram map, $T$, can be defined analogously for twisted $n$-gons.
We always assume twisted $n$-gons are generic in the sense
that no three consecutive points $p_i,p_{i+1},p_{i+2}$ are collinear. 

\medskip

\begin {defn}
    Two twisted $n$-gons $(p_i)$ and $(q_i)$ said to be \emph{projectively equivalent} if there is some $G \in \mathrm{PGL}_3(\Bbb{R})$ so that
    $q_i = G p_i$ for every $i$. We will denote by $\mathcal{P}_n$ the moduli space of twisted $n$-gons under projective equivalence.
\end {defn}

\begin {rmk}
    The pentagram map commutes with projective transformations in the sense that for $G \in \mathrm{PGL}_3(\Bbb{R})$, the polygons
    $T(P)$ and $T(G \cdot P)$ are projectively equivalent. Thus it induces a well-defined map $T \colon \mathcal{P}_n \to \mathcal{P}_n$.
\end {rmk}

\medskip

There is the following important result about this map:

\bigskip

\begin {thm} \cite{ost10}
    The pentagram map $T \colon \mathcal{P}_n \to \mathcal{P}_n$ is completely integrable, in the Liouville sense. That is,
    there exists a $T$-invariant Poisson structure on $\mathcal{P}_n$ and sufficiently many independent functions $f_i \in C^\infty(\mathcal{P}_n)$
    for which $\{f_i,f_j\} = 0$.
\end {thm}

\bigskip

\subsection {Corrugated Polygons and Higher Pentagram Maps}

\bigskip

In \cite{gsvt_16}, the authors define \emph{generalized higher pentagram maps}, which we define and discuss now. Instead of working
in $\Bbb{P}^2$, we generalize to $\Bbb{P}^{k-1}$, with the usual pentagram map being the specialization to $k=3$. In the same way as
before, we define twisted polygons in $\Bbb{P}^{k-1}$.

\begin {defn}
    A \emph{twisten $n$-gon} in $\Bbb{P}^k$ is a bi-infinite sequences of points $(p_i)_{i \in \Bbb{Z}}$ in $\Bbb{P}^{k-1}$ with monodromy
    $M \in \mathrm{PGL}_{k}$ so that $p_{i+n} = M p_i$ for all $i$. Define $\mathcal{P}_{k,n}$ to be the set of all projective equivalence
    classes of twisted $n$-gons in $\Bbb{P}^{k-1}$ with the genericity condition that any consecutive $k$ points do not lie in a proper
    projective subspace. Keeping with the previous notation, we define $\mathcal{P}_n := \mathcal{P}_{3,n}$.
\end {defn}

\medskip

\begin {defn}
    A \emph{corrugated polygon} is a twisted $n$-gon in $\Bbb{P}^{k-1}$ with the additional property that
    $p_i$, $p_{i+1}$, $p_{i+k-1}$, $p_{i+k}$ span a projective plane for all $i$.
    In other words, for any lift $v_i$ of $p_i$ to $\Bbb{R}^k$, the four vectors $v_i$, $v_{i+1}$, $v_{i+k-1}$, $v_{i+k}$ span a $3$-dimensional
    subspace. In particular, when $k=3$, any twisted polygon is automatically corrugated.
\end {defn}

\medskip

\begin {defn}
    Define $\mathcal{P}^0_{k,n}$ to be the subset of corrugated polygons with the property that for each $i$, any $3$ of the $4$
    points $p_i$, $p_{i+1}$, $p_{i+k-1}$, $p_{i+k}$ are not collinear. In other words, for any lift $v_i$, and any $i$, any $3$ of the $4$
    vectors $v_i$, $v_{i+1}$, $v_{i+k-1}$, $v_{i+k}$ are linearly independent. 
\end {defn}

\medskip

Define $\mathscr{L}_i$ to be the line containing $p_i$ and $p_{i+k-1}$.
Then since $p_i$, $p_{i+1}$, $p_{i+k-1}$, $p_k$ span a projective plane in $\Bbb{P}^{k-1}$, the lines $\mathscr{L}_i$
and $\mathscr{L}_{i+1}$ must intersect. We call the intersection point $q_i$. We can then define a
pentagram map $T \colon \mathcal{P}^0_{k,n} \to \mathcal{P}_{k,n}$ which sends $(p_i)$ to $(q_i)$. Note the codomain is not
$\mathcal{P}^0_{k,n}$, since the image my be degenerate.

\bigskip

\subsection {Coordinates on the Moduli Space}

\bigskip

Next we will construct a system of coordinates on $\mathcal{P}_n$, following the presentation in \cite{gsvt_16}. 
The following construction will also give a system of coordinates on $\mathcal{P}^0_{k,n}$ for the generalized higher
pentagram maps, but we present only the case $k=3$, for the sake of simplicity.
Given a twisted $n$-gon
$(p_i)$, we lift it to a bi-infinite sequence of vectors $(v_i)$ in $\Bbb{R}^3$. That is, $v_i$ projects to $p_i$ under the canonical
map $\Bbb{R}^3 \setminus \{0\} \to \Bbb{P}^2$.
The genericity assumption guarantees that for each $i$,
$\{v_i,v_{i+1},v_{i+2}\}$ is a basis for $\Bbb{R}^3$. Thus we have for each $i$ a linear dependence relation:
\begin {equation} \label{eq:lin_dep}
    v_{i+3} = a_i v_i + b_i v_{i+1} + c_i v_{i+2}
\end {equation}
We call the lift ``twisted'' if $v_{i+n} = A v_i$, where $A \in \mathrm{GL}_3$ is any matrix representing the monodromy.
If the lift is twisted, then the sequences $a_i$, $b_i$, and $c_i$ are $n$-periodic,
and none of the $a_i$, $b_i$, or $c_i$ are zero. The following result is
proved in more generality in \cite{gsvt_16}, but we prove it here for the sake of presentation.

\begin {prop} \label{prop:normal_lift} \cite{gsvt_16}
    The lift $(v_i)$ can be chosen so that $c_i = 1$ for all $i$, and the remaining coefficients are $n$-periodic.
\end {prop}
\begin {proof}
    We first start by choosing a twisted lift $(v_i)$, which is always possible.
    Any other lift $(\hat{v}_i)$ differs by rescaling. That is, there are non-zero constants $\lambda_i$
    so that $\hat{v}_i = \lambda_i v_i$. Then \textbf{Equation \ref{eq:lin_dep}} becomes:
    \[ \hat{v}_{i+3} = \left(a_i \frac{\lambda_{i+3}}{\lambda_i} \right) \hat{v}_i + 
                       \left(b_i \frac{\lambda_{i+3}}{\lambda_{i+1}} \right) \hat{v}_{i+1} + 
                       \left(c_i \frac{\lambda_{i+3}}{\lambda_{i+2}} \right) \hat{v}_{i+2} \]
    We want to choose $\lambda_i$ so that $c_i \frac{\lambda_{i+3}}{\lambda_{i+2}} = 1$ for all $i$. Re-arranging this equation, and
    re-indexing gives the recurrence $\lambda_{i+1} = \frac{\lambda_i}{c_{i-2}}$. We may therefore set $\lambda_0 = 1$,
    and then define the rest by this recurrence.

    \medskip

    Solving this recurrence gives $\lambda_{i+n} = \frac{\lambda_i}{c_1 \cdots c_n}$. The factors of $c_1 \cdots c_n$ therefore cancel in
    the expressions of the new coefficients, and since $a_i$, $b_i$, $c_i$ were periodic, so are the new coefficients.
\end {proof}

With the previous result in mind, we change notation and let $x_i := b_i$ and $y_i := a_i$, so that \textbf{Equation \ref{eq:lin_dep}} becomes
\begin {equation} \label{eq:lin_dep_xy}
    v_{i+3} = y_i v_i + x_i v_{i+1} + v_{i+2}
\end {equation}

\medskip

\begin {lem} \cite{gsvt_16}
    $\dim \mathcal{P}^0_{k,n} = 2n$, with the $x_i$, $y_i$ being a system of coordinates.
\end {lem}

\bigskip

Next we will write the pentagram map in these coordinates, and see that it is a rational map. 
Since we consider the pentagram map acting on labeled polygons, we will abuse notation slightly
and write $T(p_i) = q_i$ for the individual vertices of the polygon. Taking the abuse a step further, if we have lifts $(v_i)$ of $(p_i)$ and
$(w_i)$ of $(q_i) = (T(p_i))$, we will also write $T(v_i) = w_i$. To see how the pentagram map transforms the $x_i$ and $y_i$ coordinates, we will
look at how the pentagram map acts on a lift. First let's introduce some notation. As before, let $\mathscr{L}_i$ be the line between $p_i$ and $p_{i+2}$.
Recall that the pentagram map is given by $q_i = T(p_i) = \mathscr{L}_i \cap \mathscr{L}_{i+1}$.
Let $\mathscr{P}_i$ be the plane in $\Bbb{R}^3$ which projects onto $\mathscr{L}_i$. Then any vector in the
intersection $\mathscr{P}_i \cap \mathscr{P}_{i+1}$ can be taken to be the lift of the image $w_i = T(v_i)$. 
In particular, re-arrangeing \textbf{Equation \ref{eq:lin_dep_xy}} gives a candidate which can be written in two ways:
\begin {equation} \label{eq:pent_lifts}
    v_{i+3} - x_i v_{i+1} = y_i v_i + v_{i+2} \in \mathscr{P}_i \cap \mathscr{P}_{i+1}
\end {equation}

\begin {prop} \label{prop:pentagram_in_coords} \cite{gsvt_16}
    The pentagram map $T \colon \mathcal{P}_n \to \mathcal{P}_n$ is given in the $x_i,y_i$ coordinates by
    \begin {align*}
        x_i &\mapsto x_i \, \frac{x_{i+2} + y_{i+3}}{x_i + y_{i+1}} \\
        y_i &\mapsto y_{i+1} \, \frac{x_{i+2} + y_{i+3}}{x_i + y_{i+1}}
    \end {align*}
\end {prop}
\begin {proof}
    The lifts of the image polygon's vertices will also satisfy the linear dependence relation as in \textbf{Equation \ref{eq:lin_dep}}:
    \[ T(v_{i+3}) = Y_i T(v_i) + X_i T(v_{i+1}) + Z_i T(v_{i+2}) \]
    In the above equation, substitute for $T(v_i)$ and $T(v_{i+2})$ the expressions from \textbf{Equation \ref{eq:pent_lifts}}
    belonging to $\mathscr{P}_{i+1}$, and substitute for $T(v_{i+1})$ and $T(v_{i+3})$ the corresponding expressions belonging to $\mathscr{P}_i$.
    Doing so, one obtains the equation
    \[ v_{i+5} + y_{i+3} v_{i+3} = \left( X_i y_{i+1} - Y_i x_i \right) v_{i+1} + \left( X_i+Y_i-Z_ix_{i+2}\right)v_{i+3} + Z_i v_{i+5} \]
    Comparing coefficients on both sides, we can solve to get that
    \begin {align*}
        Z_i &= 1 \\
        X_i &= x_i \frac{x_{i+2} + y_{i+3}}{x_i + y_{i+1}} \\
        Y_i &= y_{i+1} \frac{x_{i+2} + y_{i+3}}{x_i + y_{i+1}}
    \end {align*}
\end {proof}

\bigskip

\section {Networks and Poisson Structures}

\bigskip

In this section, we review the relevant definitions and constructions needed to formulate
the pentagram map in terms of edge-weighted directed graphs. First, we review the types of weighted
directed graphs we will be working with, and an important class of transformations of such graphs.
Then we review the Poisson structures introduced in \cite{gsv_09}\cite{gsv_10} on the
space of weights of such graphs.

\bigskip

\subsection {Weighted Directed Fat Graphs}

\bigskip

\begin {defn}
    A \emph{fat graph} is a graph, together with a prescribed cyclic ordering of the half-edges incident to each vertex.
\end {defn}

\medskip

In particular, any graph embedded onto a two-dimensional oriented surface is naturally a fat graph, with the orientation
induced from that of the surface. Later on, we will only be considering such graphs which are embedded on a surface, so all graphs
will assumed to be fat graphs, even if not explicitly stated. 

\bigskip

\begin {defn}
    A \emph{quiver} (or \emph{directed graph}) is a tuple $Q = (Q_0,Q_1,s,t)$ where $Q_0$ is the set of vertices of the underlying graph,
    $Q_1$ the set of edges, and $s,t \colon Q_1 \to Q_0$ are the ``source'' and ``target'' maps, indicating the direction of the arrows.
\end {defn}

\medskip

\begin {defn}
    We will use the term \emph{network} to mean a weighted directed fat graph. A weighting means an 
    assignment $Q_1 \to \Bbb{R}$ of a real number to each edge.
\end {defn}

\bigskip

Postnikov considered what he called \emph{perfect planar networks} \cite{postnikov_06}, which are networks embedded in a disk such that
\begin {enumerate}
    \item All vertices on the boundary are univalent
    \item All internal vertices are trivalent and are neither sources nor sinks
    \item All edge weights are positive real numbers
\end {enumerate}

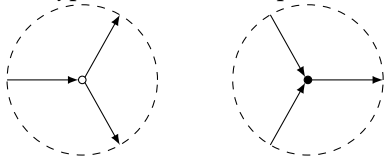
\begin {figure}[h!]
\centering
\caption {Types of vertices in a perfect network}
\label {fig:white_black_vertices}
\begin {tikzpicture}
    \draw [dashed]     (0,0) circle (1cm);
    \draw (0,0) circle (0.05cm);
    \draw [-latex]     (-1,0) -- (-0.05,0);
    \draw [-latex]     (0.035,0.035) -- (0.5,0.866);
    \draw [-latex]     (0.035,-0.035) -- (0.5,-0.866);

    \draw [dashed]     (3,0) circle (1cm);
    \draw [fill=black] (3,0) circle (0.05cm);
    \draw [-latex]     (3.05,0) -- (4,0);
    \draw [-latex]     (2.5,0.866) -- (2.965,0.035);
    \draw [-latex]     (2.5,-0.866) -- (2.965,-0.035);
\end {tikzpicture}
\end {figure}

Note that the second condition, that internal vertices are trivalent, and are not sources or sinks, implies that they are one of two types:
either a vertex has one incoming and two outgoing edges, or it has one outgoing and two incoming edges. When we draw networks, we will
picture the former as white vertices, and the latter as black. This is pictured in \textbf{Figure \ref{fig:white_black_vertices}}.

\bigskip

\begin {defn} \cite{postnikov_06}
    For an acyclic perfect network $Q$, let $a_1,\dots,a_k$ denote the sources on the boundary, and let $b_1,\dots,b_\ell$ be the sinks on the boundary.
    The \emph{boundary measurement} $b_{ij}$ is defined to be the sum of the weights of all paths from $a_i$ to $b_j$, where the weight
    of a path is the product of the weights of its edges. 
    We arrange the boundary measurements into a $k \times \ell$ matrix $\mathscr{B}(Q) = (b_{ij})$, called the \emph{boundary measurement matrix}
    of $Q$.
\end {defn}

\medskip

Later on, we will be interested only in a particular acyclic network, so we will not discuss the
more general definition for when $Q$ has directed cycles. 

\bigskip

\begin {figure}[h!]
\centering
\caption {Local Postnikov moves}
\label {fig:postnikov_moves}
\begin {tikzpicture}[scale=0.8, every node/.style={scale=0.8}]
    % Label (I)
    \draw (-4,0) node {( I )};

    % Before
    \draw [dashed]     (0,0) circle (2cm);
    \draw (0.707,0.707) circle (0.05cm);
    \draw[fill=black] (0.707,-0.707) circle (0.05cm);
    \draw[fill=black] (-0.707,0.707) circle (0.05cm);
    \draw (-0.707,-0.707) circle (0.05cm);

    \draw [-latex] (-1.414,1.414) -- (-0.75,0.75);
    \draw [-latex] (0.75,0.75) -- (1.414,1.414);
    \draw [-latex] (0.75,-0.75) -- (1.414,-1.414);
    \draw [-latex] (-1.414,-1.414) -- (-0.74,-0.74);

    \draw [-latex] (-0.65,0.707) -- (0.65,0.707);
    \draw [-latex] (-0.65,-0.707) -- (0.65,-0.707);
    \draw [-latex] (-0.707,-0.65) -- (-0.707,0.65);
    \draw [-latex] (0.707,0.65) -- (0.707,-0.65);

    \draw (-1,0) node {$a$};
    \draw (0,-1) node {$b$};
    \draw (1,0) node {$c$};
    \draw (0,1) node {$d$};

    \draw (-0.9,1.3) node {$e_1$};
    \draw (-0.9,-1.3) node {$e_2$};
    \draw (0.9,-1.3) node {$e_3$};
    \draw (0.9,1.3) node {$e_4$};

    % Arrow

    \draw [-latex] (2.5,0) -- (3.5,0);

    % After
    \draw [dashed]    (6+0,0)           circle (2cm);
    \draw[fill=black] (6+0.707,0.707)   circle (0.05cm);
    \draw[fill=white] (6+0.707,-0.707)  circle (0.05cm);
    \draw[fill=white] (6+-0.707,0.707)  circle (0.05cm);
    \draw[fill=black] (6+-0.707,-0.707) circle (0.05cm);

    \draw [-latex] (6+-1.414,1.414)  -- (6+-0.75,0.75);
    \draw [-latex] (6+0.75,0.75)     -- (6+1.414,1.414);
    \draw [-latex] (6+0.75,-0.75)    -- (6+1.414,-1.414);
    \draw [-latex] (6+-1.414,-1.414) -- (6+-0.74,-0.74);

    \draw [-latex] (6+-0.65,0.707)  -- (6+0.65,0.707);
    \draw [-latex] (6+-0.65,-0.707) -- (6+0.65,-0.707);
    \draw [-latex] (6+-0.707,0.65) -- (6+-0.707,-0.65);
    \draw [-latex] (6+0.707,-0.65)   -- (6+0.707,0.65);

    \draw (6-0.1+-1,0) node {$d \, \frac{c}{\Delta}$};
    \draw (6+0,-1) node {$\Delta$};
    \draw (6+0.1+1,0)  node {$d \, \frac{a}{\Delta}$};
    \draw (6+0,1)  node {$d \, \frac{b}{\Delta}$};

    \draw (6+-0.9,1.3) node {$e_1$};
    \draw (6+-0.9,-1.3) node {$e_2$};
    \draw (6+0.9,-1.3) node {$e_3$};
    \draw (6+0.9,1.3) node {$e_4$};

    \draw (10,0) node {$(\Delta := b + adc)$};

    % Type (II)

    \draw (-4,-5) node {( II )};

    % Before
    \draw [dashed]     (0,-5)      circle (2cm);
    \draw [fill=white] (-0.667,-5) circle (0.05cm);
    \draw [fill=white] (0.667,-5)  circle (0.05cm);

    \draw [-latex] (-2,-5) -- (-0.667-0.05,-5);
    \draw [-latex] (-0.667+0.05,-5) -- (0.667-0.05,-5);
    \draw [-latex] (0.667+0.05,-5) -- (2,-5);
    \draw [-latex] (-0.667,-5+0.05) -- (-0.667,-5+1.886);
    \draw [-latex] (0.667,-5-0.05) -- (0.667,-5-1.886);

    \draw (-1.33,-5.2) node {$a$};
    \draw (0,-5.2)     node {$b$};
    \draw (1.33,-5.2)  node {$c$};
    \draw (-0.85,-4)   node {$x$};
    \draw (0.85,-6)    node {$y$};

    % Arrow
    \draw [-latex] (2.5,-5) -- (3.5,-5);

    % After
    \draw [dashed]     (6+0,-5)      circle (2cm);
    \draw [fill=white] (6+-0.667,-5) circle (0.05cm);
    \draw [fill=white] (6+0.667,-5)  circle (0.05cm);

    \draw [-latex] (6+-2,-5)          -- (6+-0.667-0.05,-5);
    \draw [-latex] (6+-0.667+0.05,-5) -- (6+0.667-0.05,-5);
    \draw [-latex] (6+0.667+0.05,-5)  -- (6+2,-5);
    \draw [-latex] (6+0.667,-5+0.05) -- (6+-0.667,-5+1.886);
    \draw [-latex] (6-0.667,-5-0.05)  -- (6+0.667,-5-1.886);

    \draw (6+-1.33,-5.2) node {$a$};
    \draw (6+0,-5.2)     node {$b$};
    \draw (6+1.33,-5.2)  node {$c$};
    \draw (6+0.3,-5+1)   node {$\frac{x}{b}$};
    \draw (6-0.3,-5-1)   node {$yb$};

    % Type (III)

    \draw (-4,-10) node {( III )};

    % Before
    \draw [dashed]     (0,-10)      circle (2cm);
    \draw [fill=black] (-0.667,-10) circle (0.05cm);
    \draw [fill=black] (0.667,-10)  circle (0.05cm);

    \draw [-latex] (-2,-10) -- (-0.667-0.05,-10);
    \draw [-latex] (-0.667+0.05,-10) -- (0.667-0.05,-10);
    \draw [-latex] (0.667+0.05,-10) -- (2,-10);
    \draw [-latex] (-0.667,-10+1.886) -- (-0.667,-10+0.05);
    \draw [-latex] (0.667,-10-1.886) -- (0.667,-10-0.05);

    \draw (-1.33,-10.2) node {$a$};
    \draw (0,-10.2)     node {$b$};
    \draw (1.33,-10.2)  node {$c$};
    \draw (-0.85,-9)    node {$x$};
    \draw (0.85,-11)    node {$y$};

    % Arrow
    \draw [-latex] (2.5,-10) -- (3.5,-10);

    % After
    \draw [dashed]     (6+0,-10)      circle (2cm);
    \draw [fill=black] (6+-0.667,-10) circle (0.05cm);
    \draw [fill=black] (6+0.667,-10)  circle (0.05cm);

    \draw [-latex] (6+-2,-10)          -- (6+-0.667-0.05,-10);
    \draw [-latex] (6+-0.667+0.05,-10) -- (6+0.667-0.05,-10);
    \draw [-latex] (6+0.667+0.05,-10)  -- (6+2,-10);
    \draw [-latex] (6+-0.667,-10+1.886) -- (6+0.667,-10+0.05);
    \draw [-latex] (6+0.667,-10-1.886) -- (6-0.667,-10-0.05);

    \draw (6+-1.33,-10.2) node {$a$};
    \draw (6+0,-10.2)     node {$b$};
    \draw (6+1.33,-10.2)  node {$c$};
    \draw (6+0.3,-10+1)   node {$xb$};
    \draw (6-0.3,-10-1)   node {$\frac{y}{b}$};

\end {tikzpicture}
\end {figure}
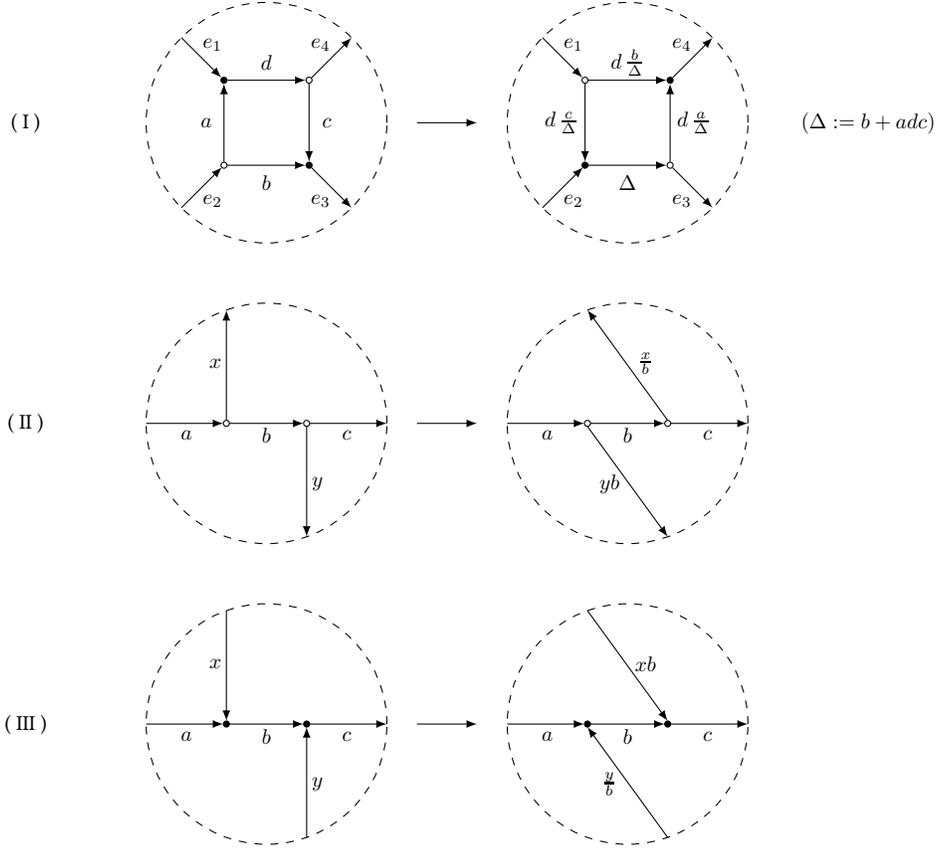

Postnikov described several local transformations of networks which leave the boundary measurements invariant.
The moves are pictured in \textbf{Figure \ref{fig:postnikov_moves}}, and are labeled with the corresponding edge weights.

\bigskip

We refer to the type I move as the ``square move'', 
the type II move as ``white-swapping'', and the type III move as ``black-swapping''.
It is an easy exercise to check that these three moves do not change the boundary measurements from a source to a sink in these pictures.

\bigskip

There is another type of local move called a \emph{gauge transformation}, which changes the edge weights but not the underlying graph. 
Let $\lambda$ be any non-zero number.
Then at any vertex we may multiply all incoming weights by $\lambda$ and all outgoing weights by $\lambda^{-1}$. This obviously
does not change the boundary measurements. The group generated by all gauge transformations is called the \emph{gauge group}.

\bigskip

From now on, we will only consider networks where the sources and sinks are not interlaced, so that all sources can be pictured
on the left side, and all sinks on the right. The reason for restricting to these networks is that they have the following nice
property.
One may consider ``concatenating'' two such networks in a disk, by glueing segments of their boundaries in a way that each sink
on one is identified with a source of the other. This is pictured in \textbf{Figure \ref{fig:concat_local}}. 
The resulting edge after identification is weighted by the product of the two edges.
If the sources and sinks are not interlaced, as mentioned above, then the boundary measurement matrix of the resulting network
after concatenation is the product of the two boundary measurement matrices.

\bigskip

We may also consider weighted directed networks embedded on an annulus \cite{gsv_10}. Later, when we talk about the pentagram map, it will be
these annular networks that we will consider. For simplicity, we will again only consider acyclic networks.
Similar to the disk case, consider \emph{perfect} networks, meaning any internal vertices are trivalent (and not sources or sinks), and that all boundary vertices are univalent.
We also assume all vertices on the inner boundary component are sources, and all vertices on the outer boundary component are sinks.
Then in a similar way as before, we define a boundary measurement matrix, whose entries are the sums of path weights from a given source to a given sink.
Since the Postnikov moves are local, and do not take into account the global topology of the surface on which the network is embedded, 
they may also be applied in this case, and of course they still do not change the boundary measurements.

\bigskip

As in the disk case, we may consider concatenation of two annular networks, by identifying the outer boundary circle of one with the inner
boundary circle of the other, in a way that each sink is identified with a source. If we assume, as mentioned above, that all sources are
on the inner boundary circle and all sinks on the outer, then the resulting boundary measurement matrix
is the product of the two boundary measurement matrices.

\bigskip

If there are the same number of sources and sinks, then it will sometimes be convenient to glue the inner and outer boundary circles together so that 
pairs of sources and sinks are identified, to obtain a network on the torus. 

\bigskip

Lastly, we define the \emph{modified edge weights}, which are elements of the  Laurent polynomial ring $\Bbb{R}[\lambda^\pm]$
in the indeterminate $\lambda$, defined as follows. 
Choose an oriented curve $\rho$, called the \emph{cut}, which connects the inner and outer boundary components of the annulus.
Suppose an oriented edge $\alpha \in Q_1$ has weight $x_\alpha \in \Bbb{R}$. It may
happen that $\alpha$ intersects the cut. If $i$ is an intersection point, define $\varepsilon_i = 1$
if the tangent vectors to $\alpha$ and $\rho$ at $i$ form an oriented basis of the tangent space, and define $\varepsilon_i = -1$
if they have the opposite orientation. Then the modified edge weight of $\alpha$ is defined as
\[  x_\alpha \lambda^{\varepsilon_i} \]
If we use the modified edge weights, then the boundary measurements become Laurent polynomials in $\lambda$.
From now on, we always assume the boundary measurements are in terms of modified edge weights.
We will use the notation $\mathscr{B}_Q(\lambda)$ to emphasize that the entries are functions of $\lambda$.

\medskip

\begin {rmk}
    As mentioned above, if there are the same number of sources and sinks, we may glue the two boundary circles together
    to obtain a network on the torus. The two boundary circles become a single loop on the torus, which we call the \emph{rim}.
    If we choose a new rim, and cut the torus along this new rim, we get a new (possibly different) network on an annulus.
    The boundary measurement matrices differ by a re-factorization, which we will now explain. If we draw the new rim on the old annulus,
    we can realize the annulus as the concatenation of two annular networks $Q'$ and $Q''$, one on each side of the new rim.
    Let $A$ and $B$ be the boundary measurement matrices of $Q'$ and $Q''$ respectively.
    Then the boundary measurement matrix is $AB$. Glueing into a torus and then cutting along the new rim amounts to concatenating
    the pictures of $Q'$ and $Q''$ in the opposite order. The new boundary measurement matrix is thus $BA$.
    It was observed by Izosimov \cite{izosimov_18} that many different generalizations of the pentagram map may be described
    in terms of these types of matrix re-factorizations in Poisson Lie groups.
\end {rmk}

\bigskip

\subsection {Poisson Structures on the Space of Edge Weights}

\bigskip

Given a perfect network in an annulus, as in the last section, we may forget the particular choice of edge weights and consider the underlying quiver $Q = (Q_0,Q_1)$.
The \emph{space of edge weights} of the quiver $Q$, denoted $\mathscr{E}_Q := (\mathbb{R}^*)^{Q_1}$, is the space of all possible choices of non-zero weights.
Gekhtman, Shapiro, and Vainshtein defined a family of Poisson structures on $\mathscr{E}_Q$ \cite{gsv_09}, which we will review now.

\bigskip

The Poisson structures are defined locally at each vertex, and then the local Poisson brackets are shown to be ``compatible'' with concatenation
(this will be made more precise later), and so they can be combined to give a global Poisson bracket. This is outlined in \cite{gsv_09} and \cite{gsv_10},
but we present the construction here in detail, since we will mimic it very closely later on when we define a non-commutative analogue. Recall that
for networks on an annulus, we require internal vertices are trivalent, and that they are neither sources nor sinks, and so they are either
white or black, as mentioned before:

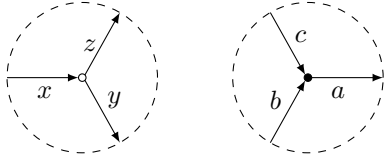
\begin {center}
\begin {tikzpicture}
    \draw [dashed]     (0,0) circle (1cm);
    \draw (0,0) circle (0.05cm);
    \draw [-latex]     (-1,0) -- (-0.05,0);
    \draw [-latex]     (0.035,0.035) -- (0.5,0.866);
    \draw [-latex]     (0.035,-0.035) -- (0.5,-0.866);

    \draw (-0.5,-0.2) node {$x$};
    \draw (0.1,0.43) node {$z$};
    \draw (0.43,-0.3) node {$y$};

    \draw [dashed]     (3,0) circle (1cm);
    \draw [fill=black] (3,0) circle (0.05cm);
    \draw [-latex]     (3.05,0) -- (4,0);
    \draw [-latex]     (2.5,0.866) -- (2.965,0.035);
    \draw [-latex]     (2.5,-0.866) -- (2.965,-0.035);

    \draw (3.4,-0.2) node {$a$};
    \draw (2.57,-0.3) node {$b$};
    \draw (2.9,0.55) node {$c$};
\end {tikzpicture}
\end {center}

Let $\mathscr{E}_\circ$ be the space $(\Bbb{R}^*)^3$ with coordinates $x,y,z$. We will think of the variables $x,y,z$
as representing the edge weights at a white vertex picture above. Similarly, let $\mathscr{E}_\bullet = (\Bbb{R}^*)^3$ with coordinates $a,b,c$
corresponding to a black vertex. Of course they are diffeomorphic, but we will define different Poisson brackets on them.
Choose any $w_1,w_2,w_3 \in \Bbb{R}$, and define a \emph{log-canonical} Poisson bracket $\{-,-\}_\circ$ on $\mathscr{E}_\circ$ by

\[ \{x,y\}_\circ = w_1 xy, \quad \{x,z\}_\circ = w_2  xz, \quad \{y,z\}_\circ = w_3 yz \]

Similarly, for scalars $k_1,k_2,k_3 \in \Bbb{R}$, define $\{-,-\}_\bullet$ on $\mathscr{E}_\bullet$ by

\[ \{a,b\}_\bullet = k_1 ab, \quad \{a,c\}_\bullet = k_2 ac, \quad \{b,c\}_\bullet = k_3 bc \]

Next we will consider the operation of glueing/concatenating these local pictures together, as described before. 
This means identifying part of the boundary of one picture with part of another,
in such a way that a sink is identified with a source.
We then erase the common glued boundary, and remove the common identified vertex, identifying the two incident edges
(which have a consistent orientation by construction). The weight of the new edge is associated with the product of
the weights of the edges being glued. \textbf{Figure \ref{fig:concat_local}} illustrates glueing a black and white local picture:

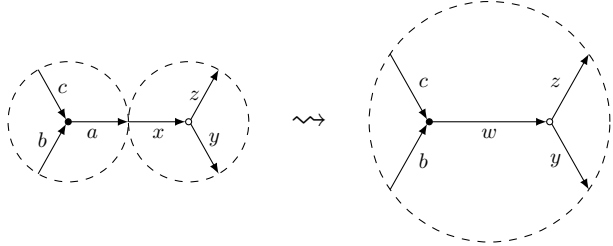
\begin {figure}[h!]
\centering
\caption {Concatenation of Local Pictures}
\label {fig:concat_local}
\begin {tikzpicture}[scale=0.8, every node/.style={scale=0.8}]
    % Left Circle
    \draw [dashed]     (0,0) circle (1cm);
    \draw [fill=black] (0,0) circle (0.05cm);
    \draw [-latex]     (0.05,0) -- (1,0);
    \draw [-latex]     (-0.5,0.866) -- (2.965-3,0.035);
    \draw [-latex]     (-0.5,-0.866) -- (2.965-3,-0.035);

    \draw (0.4,-0.2) node {$a$};
    \draw (2.57-3,-0.3) node {$b$};
    \draw (2.9-3,0.55) node {$c$};

    % Right Circle
    \draw [dashed]     (2,0) circle (1cm);
    \draw (2,0) circle (0.05cm);
    \draw [-latex]     (1,0) -- (2-0.05,0);
    \draw [-latex]     (2.035,0.035) -- (2.5,0.866);
    \draw [-latex]     (2.035,-0.035) -- (2.5,-0.866);

    \draw (1.5,-0.2) node {$x$};
    \draw (2.1,0.43) node {$z$};
    \draw (2.43,-0.3) node {$y$};

    \draw (4,0) node {\LARGE $\rightsquigarrow$};

    % Combined Circle
    \draw[dashed]     (7,0) circle (2cm);
    \draw[fill=black] (6,0) circle (0.05cm);
    \draw[fill=white] (8,0) circle (0.05cm);

    \draw[-latex] (6.05,0) -- (7.95,0);
    \draw[-latex] (8.035,0.035) -- (8.6514,1.128);
    \draw[-latex] (8.035,-0.035) -- (8.6514,-1.128);
    \draw[-latex] (6-0.6514,1.128) -- (6-0.035,0.035);
    \draw[-latex] (6-0.6514,-1.128) -- (6-0.035,-0.035);

    \draw (7,-0.2) node {$w$};
    \draw (5.9,-0.65) node {$b$};
    \draw (5.9,0.65) node {$c$};
    \draw (8.1,0.65) node {$z$};
    \draw (8.1,-0.65) node {$y$};
\end {tikzpicture}
\end {figure}

The definitions of $\{-,-\}_\circ$ and $\{-,-\}_\bullet$ are ``compatible'' with this glueing procedure in a way that we now make precise.
Recall that $\mathscr{E}_Q = (\mathbb{R}^*)^{Q_1}$ is the space of edge weights. We can associate to it the algebra $E_Q := \Bbb{R}(Q_1)$,
which is the field of rational functions with indeterminates corresponding to the arrows of the quiver.
We also define $\mathscr{H}_Q := (\Bbb{R}^*)^{Q_1 \coprod Q_1}$
to be the \emph{space of half-edge weights}. This is an assignment of a non-zero number to each half-edge. Each half edge is associated to
the vertex which it is incident to. For an arrow $\alpha \in Q_1$, we call the corresponding half edges $\alpha_s$ and $\alpha_t$, for the
source and target ends of the arrow. The algebra of functions $H_Q$ corresponding to $\mathscr{H}_Q$ is then the field of rational functions
in twice as many variables, corresponding to $\alpha_s$ and $\alpha_t$ for $\alpha \in Q_1$.

\medskip

Recall that $Q_0$ denotes the set of vertices of the quiver $Q$. We can partition this into
$Q_0 = V_\circ \cup V_\bullet \cup V_\partial = V_i \cup V_\partial$, where $V_\circ$ is the set of white vertices,
$V_\bullet$ the black vertices, and $V_\partial$ the boundary vertices, and $V_i = V_\circ \cup V_\bullet$ is the
set of internal vertices.
Note that since all interal vertices are trivalent, and all boundary vertices are univalent,
the set of half-edges is in bijection with $V_\partial \coprod V_i \coprod V_i \coprod V_i$.
Recall we defined the spaces $\mathscr{E}_\circ \cong \mathscr{E}_\bullet \cong (\Bbb{R}^*)^3$.
We now additionally define $\mathscr{E}_\partial \cong \Bbb{R}^*$ with trivial Poisson bracket, corresponding
to boundary vertices. Then the space of half-edges can be thought of as 
\[ \mathscr{H}_Q \cong \prod_{\ast \in V_\partial}   \mathscr{E}_\partial \, \times \, 
                       \prod_{\circ \in V_\circ}     \mathscr{E}_\circ    \, \times \, 
                       \prod_{\bullet \in V_\bullet} \mathscr{E}_\bullet \]
Under this identification, we realize the algebra $H_Q$ as the tensor product
\[ H_Q \cong \bigotimes_{\ast \in V_\partial}   E_\partial \, \otimes \, 
             \bigotimes_{\circ \in V_\circ}     E_\circ    \, \otimes \, 
             \bigotimes_{\bullet \in V_\bullet} E_\bullet \]
From the Poisson brackets defined above on the algebras $E_\circ$, $E_\bullet$, and $E_\partial$,
we get a natural induced Poisson bracket on $H_Q$, where $\{\alpha,\beta\} = 0$ if $\alpha$ and
$\beta$ are half-edges which are not incident to a common vertex, and the boundary half-edges
are casimirs.

\medskip

We define a ``glueing'' map $g \colon \mathscr{H}_Q \to \mathscr{E}_Q$, represented by \textbf{Figure \ref{fig:concat_local}}, 
given by $(\alpha_s,\alpha_t)_{\alpha \in Q_1} \mapsto (\alpha_s \alpha_t)_{\alpha \in Q_1}$. The earlier claim that
the local brackets on $\mathscr{E}_\circ$ and $\mathscr{E}_\bullet$ are compatible with glueing is made precise by the following.

\bigskip

\begin {prop}
    There is a unique Poisson bracket on $\mathscr{E}_Q$ such that the glueing map $g \colon \mathscr{H}_Q \to \mathscr{E}_Q$
    is a Poisson morphism.
\end {prop}
\begin {proof}
    The statement that $g$ is a Poisson map is equivalent to the pull-back map $g^* \colon E_Q \to H_Q$ being a 
    homomorphism of Poisson algebras. Note that for $\alpha \in Q_1$, the map $g^*$ is given by $g^*(\alpha) = \alpha_s \alpha_t$.
    To be a Poisson algebra homomorphism would mean that for any $\alpha,\beta \in Q_1$,
    \[ \{g^*(\alpha),g^*(\beta)\}_{H_Q} = g^* \left( \{ \alpha, \beta \}_{E_Q} \right) \]
    We claim that $\{\alpha,\beta\}_{E_Q}$ is determined by this property. If $\alpha$ and $\beta$ do not share a common vertex,
    then obviously $\{\alpha,\beta\} = 0$. There are still many cases to consider, depending on which end (source or target)
    of each $\alpha$ and $\beta$ meet at the common vertex, and where in the cyclic ordering those half-edges are at that vertex.
    For example, let's consider the picture in \textbf{Figure \ref{fig:concat_local}}, and try to define $\{w,y\}_{E_Q}$. We only
    see one end of $y$ in the figure, so let's say $g^*(y) = y_s y_t$, where $y_s$ is the end we see at the white vertex. Then the 
    condition that $g^*$ be Poisson means we must have
    \begin {align*} 
        g^*\left( \{w,y\}_{E_Q} \right) &= \{g^*(w),g^*(y)\}_{H_Q} \\
            &= \{ax, y_s y_t\}_{H_Q} \\
            &= a y_t \{x, y_s\}_{H_Q} \\
            &= w_1 \, ax y_s y_t \\
            &= w_1 g^*(w) g^*(y) \\
            &= g^* \left( w_1 \, wy \right)
    \end {align*}
    Since $g^*$ is injective, this uniquely defines $\{w,y\}_{E_Q}$. The calculations for all other cases are similarly simple.
    We see that for $\alpha,\beta \in Q_1$ with a common vertex, the bracket $\{\alpha,\beta\}_{E_Q}$ is given by the same expression
    as the bracket of the corresponding half-edges in $H_Q$.
\end {proof}

\bigskip

We will also want to consider the \emph{doubled quiver}.

\medskip

\begin {defn}
    The \emph{doubled quiver}, $\overline{Q}$, is defined to have the same vertex set as $Q$, whereas the edge set of $\overline{Q}$
    is the disjoint union of two copies of $Q_1$. For each arrow $\alpha \in Q_1$, there are two arrows $\alpha$ and $\alpha^*$ in $\overline{Q}_1$.
    The arrow $\overline{\alpha}$ is called the ``opposite'' arrow of $\alpha$, with $s(\alpha^*) = t(\alpha)$ and $t(\alpha^*) = s(\alpha)$. 
\end {defn}

\medskip

We associate $\alpha^*$ with the function $\alpha^{-1}$ on $\mathscr{E}_Q$. Then a path in $\overline{Q}$ may be represented as a
Laurent monomial, which is the product of the edge weights along the path (allowing inverses). Since Poisson brackets extend uniquely
to localizations, we may extend the Poisson bracket defined above to all Laurent polynomials (and indeed any rational functions) on $\mathscr{E}_Q$.

\bigskip

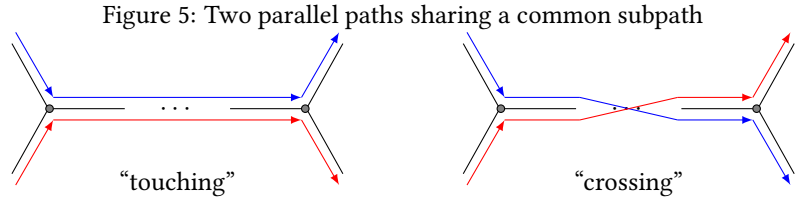
\begin {figure}[h!]
\centering
\caption {Two parallel paths sharing a common subpath}
\label {fig:paths_meet}

\begin {tikzpicture}

    % -----------------
    % touching
    % -----------------

    % Left Circle
    \draw [fill=gray] (0,0) circle (0.05cm);
    \draw (0.05,0) -- (1,0);
    \draw (-0.5,0.866) -- (2.965-3,0.035);
    \draw (-0.5,-0.866) -- (2.965-3,-0.035);

    % Right Circle
    \draw [fill=gray] (3.4,0) circle (0.05cm);
    \draw (2.4,0) -- (3.4-0.05,0);
    \draw (3.435,0.035) -- (3.9,0.866);
    \draw (3.435,-0.035) -- (3.9,-0.866);

    \draw (1.7,0) node {$\cdots$};

    \draw (1.7,-1) node {``touching''};

    \draw [blue, -latex] (0.05,0.15) -- (3.35,0.15);
    \draw [blue, -latex] (-0.45,1.016) -- (0.05,0.15);
    \draw [blue, -latex] (3.35,0.15) -- (3.85,1.016);

    \draw [red, -latex] (0.05,-0.15) -- (3.35,-0.15);
    \draw [red, -latex] (-0.45,-1.016) -- (0.05,-0.15);
    \draw [red, -latex] (3.35,-0.15) -- (3.85,-1.016);

    % -----------------
    % crossing
    % -----------------

    % Left Circle
    \draw [fill=gray] (6 + 0,0) circle (0.05cm);
    \draw (6 + 0.05,0)      -- (6 + 1,0);
    \draw (6 + -0.5,0.866)  -- (6 + 2.965-3,0.035);
    \draw (6 + -0.5,-0.866) -- (6 + 2.965-3,-0.035);

    % Right Circle
    \draw [fill=gray] (6 + 3.4,0) circle (0.05cm);
    \draw (6 + 2.4,0)        -- (6 + 3.4-0.05,0);
    \draw (6 + 3.435,0.035)  -- (6 + 3.9,0.866);
    \draw (6 + 3.435,-0.035) -- (6 + 3.9,-0.866);

    \draw (6 + 1.7,0) node {$\cdots$};

    \draw (6 + 1.7,-1) node {``crossing''};

    \draw [blue]         (6 + 0.05,0.15)   -- (6 + 1.05,0.15);
    \draw [blue]         (6 + 1.05,0.15)   -- (6 + 2.35,-0.15);
    \draw [blue, -latex] (6 + 2.35,-0.15)  -- (6 + 3.35,-0.15);
    \draw [blue, -latex] (6 + -0.45,1.016) -- (6 + 0.05,0.15);
    \draw [blue, -latex] (6 + 3.35,-0.15)  -- (6 + 3.85,-1.016);

    \draw [red]          (6 + 0.05,-0.15)   -- (6 + 1.05,-0.15);
    \draw [red]          (6 + 1.05,-0.15)   -- (6 + 2.35,0.15);
    \draw [red, -latex]  (6 + 2.35,0.15)    -- (6 + 3.35,0.15);
    \draw [red, -latex]  (6 + -0.45,-1.016) -- (6 + 0.05,-0.15);
    \draw [red, -latex]  (6 + 3.35,0.15)    -- (6 + 3.85,1.016);

\end {tikzpicture}
\end {figure}

We will now be interested in describing the Poisson bracket of two paths, thought of as rational functions on $\mathscr{E}_Q$.
To do so, it will be convenient to introduce the constants
\begin {align*}
    A_\bullet &= k_1 - k_2 - k_3 \\
    A_\circ   &= w_1 - w_2 - w_3 
\end {align*}
Consider two paths $f$ and $g$. Whenever the paths meet (share at least one edge in common), then there is
a corresponding maximal subpath which $f$ and $g$ share.
If the paths go in the same direction on the common subpath, we will say they are \emph{parallel} on that common subpath.
If $f$ and $g$ are parallel on a subpath, which is a proper subpath of both $f$ and $g$, then there are two possibilities, 
which are depicted in \textbf{Figure \ref{fig:paths_meet}}. In the figure, $f$ is the blue path, $g$ is the red path.
In the first situation (the left image in the figure), the paths are said to ``touch'', and in the second
situation, they are said to ``cross''. The other possibility is that one of the paths (either $f$ or $g$) is
a subpath of the other. In this case, if $f$ is a subpath of $g$, then we may write $g = fg'$, where $g'$ is the
rest of $g$. We extend the notions of ``touching'' and ``crossing'' to the situation where $g = fg'$ by saying $f$ and $g$
``touch'' on the subpath $f$ if $g'$ and $f$ touch, and similarly for ``crossing''.

\bigskip

We will denote by $f \cap g$ the set of all maximal common subpaths of $f$ and $g$. We will define a function $\varepsilon(f,g)$ on the set $f \cap g$,
and for $x \in f \cap g$, we denote the value by $\varepsilon_x(f,g)$. It depends on whether $f$ and $g$ touch or cross, and on
the colors of the vertices at the endpoints of the common subpath. The values are given in \textbf{Table \ref{table:def_epsilon}}.

\begin {table}[h!]
\caption {Definition of $\varepsilon(f,g)$ for parallel paths}
\label {table:def_epsilon}
\[
    \begin {array} {|c|c|c|c|} \hline
        \mathrm{type}  & \mathrm{left~endpt} & \mathrm{right~endpt} & \varepsilon_x(f,g) \\ \hline
        \mathrm{touch} & \bullet             & \bullet              & 0 \\
        \mathrm{touch} & \circ               & \circ                & 0 \\
        \mathrm{touch} & \circ               & \bullet              & A_\circ+A_\bullet \\
        \mathrm{touch} & \bullet             & \circ                & -(A_\circ+A_\bullet) \\
        \mathrm{cross} & \bullet             & \bullet              & -2A_\bullet \\
        \mathrm{cross} & \circ               & \circ                & 2A_\circ \\
        \mathrm{cross} & \circ               & \bullet              & A_\circ-A_\bullet \\
        \mathrm{cross} & \bullet             & \circ                & A_\circ-A_\bullet \\ \hline
    \end {array}
\]
\end {table}

We define $\varepsilon(f,g)$ to be skew-symmetric, so that $\varepsilon(g,f) = -\varepsilon(f,g)$ if the roles are switched. 
Also, there are other cases of intersecting paths which we have not considered. If we allow paths in $\overline{Q}$,
then two paths can meet with opposing orientations on a common subpath. This can be obtained from the local pictures above
by changing one of the paths to its inverse. Since Poisson brackets extend uniquely to localizations, this will not
be an issue.

We then have the following formula:

\begin {prop}
    Suppose $f$ and $g$ are two paths. Then
    \[ \{f,g\} = \sum_{* \in f \cap g} \varepsilon_*(f,g) \, fg \]
\end {prop}
\begin {proof}
    Since $f$ and $g$ are monomials in the coordinate variables, it is obvious that $\{f,g\} = c \, fg$ for some constant $c$.
    We must prove that $c = \sum \varepsilon_*(f,g)$. Each time $f$ and $g$ share a proper subpath on which they are parallel,
    it will look like one of the two pictures in \textbf{Figure \ref{fig:paths_meet}}. We can then write the paths as
    $f = f_0axwyb$ and $g = g_0cxwyd$, where $xwy$ is the common subpath, $a$,$x$, and $c$ are the edges incident to the
    ``left'' vertex (where the paths come together), and $y$,$b$, and $d$ are the edges incident to the ``right'' vertex
    (where the paths diverge at the end of the common subpath). Then after expanding $\{f,g\}$ by the Leibniz rule,
    we get contributions from $\{a,c\}$, $\{a,x\}$, and $\{x,c\}$ at the left vertex, and $\{y,d\}$, $\{b,y\}$, and $\{b,d\}$
    from the right vertex. Obviously the contributions from the common subpath cancel, since if $p$ and $q$ are two consecutive
    edges in the common subpath, we will get a contribution from both $\{p,q\}$ and $\{q,p\}$. There are of course other terms
    in the expansion coming from $f_0$ and $g_0$, but they will correspond to other common subpaths in the same way. It is
    therefore enough to consider the six terms mentioned above (3 each for the beginning and end vertices) for each common
    subpath.

    \medskip

    In general, we must consider 36 different cases, since each endpoint could be either white or black, and there are three
    orientations for the arrows incident to each of the endpoint vertices. A simple calculation shows that the orientations
    at the vertices are in fact irrelevant. Therefore we only have 4 cases, depending on the colors of the vertices. As mentioned above,
    by the Leibniz rule, we can add the contributions from the endpoint vertices separately.

    \medskip

    First consider the left endpoint. For all three choices of orientations, a black left endpoint will give a contribution of $-A_\bullet$
    to the coefficient, and a white left endpoint will give a contibution of $A_\circ$. For the right endpoint, a black vertex
    gives $A_\bullet$, and a white gives $-A_\circ$. Combining all the possibilities together gives the list of values in the table for $\varepsilon_x(f,g)$.

    \medskip

    We assumed so far that the common subpaths were parallel. If the paths are oriented in opposite directions on a common subpath,
    then it corresponds to reversing one of the paths in \textbf{Figure \ref{fig:paths_meet}}. Then we can extend the table for
    $\varepsilon_x(f,g)$ by the rule $\{f,g^{-1}\} = -g^{-2}\{f,g\}$.

    \medskip

    We also assumed so far that the subpaths were proper. It remains to consider when one of the paths is a subpath of the other.
    For instance, let's say $f$ is contained in $g$. Then we can write $g$ as $g=fg'$. Then by the Leibniz rule,
    $\{f,g\} = \{f,fg'\} = f\{f,g'\}$. It must be (since all vertices are trivalent) that the beginning of $g'$ is
    a common subpath with $f$, and so this is the term corresponding to the common subpath of $f$ and $g$ which is all of $f$.
    The rest of the common subpaths of $f$ and $g'$ are just usual proper common subpaths of $f$ and $g$, as discussed above.
\end {proof}

\bigskip

\begin {rmk}
    Note that this result depends only on the parameters $A_\circ$ and $A_\bullet$, and not on the actual values of $w_1$, $w_2$, $w_3$, $k_1$, $k_2$, $k_3$.
    In particular, we could choose $w_1=w_2=k_1=k_2=0$, and $w_3 = -A_\circ$, $k_3 = -A_\bullet$, and the formula from the proposition
    would be the same. From now on, we will assume this is the case (that all but $w_3$ and $k_3$ are zero). 
\end {rmk}

\bigskip

\section {Interpreting the Pentagram Map in Terms of Networks} \label{sec:model_via_networks}

\bigskip

We now use the constructions from the previous section to realize the pentagram map as a sequence of Postnikov
moves and gauge transformations of a particular network in an annulus. The properties of the boundary measurements
will give us invariants of the pentagram map, which also turn out to be involutive with respect to the
Poisson structures described in the previous section. This section summarizes the main points of \cite{gsvt_16}.

\bigskip

\subsection {The Quiver and Poisson Bracket}

\bigskip

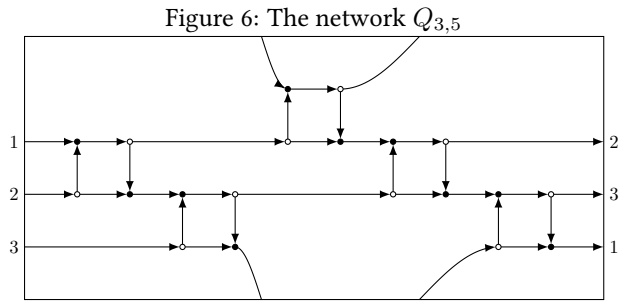
\begin {figure}[h!]
\centering
\caption {The network $Q_{3,5}$}
\label {fig:pentagram_quiver}
\begin {tikzpicture}[scale=0.7, every node/.style={scale=0.7}]
    % The rectangle
    \draw (0,0) -- (11,0) -- (11,5) -- (0,5) -- cycle;

    % The Labels
    \draw (-0.2,1) node {$3$};
    \draw (-0.2,2) node {$2$};
    \draw (-0.2,3) node {$1$};
    \draw (11.2,1) node {$1$};
    \draw (11.2,2) node {$3$};
    \draw (11.2,3) node {$2$};

    \draw [fill = black] (1,3) circle (0.05cm);
    \draw                (1,2) circle (0.05cm);
    \draw [fill = black] (2,2) circle (0.05cm);
    \draw                (2,3) circle (0.05cm);

    \draw [-latex] (0,3) -- (1-0.05,3);
    \draw [-latex] (0,2) -- (1-0.05,2);

    % Square 1

    \draw [-latex] (1,2+0.05) -- (1,3-0.05);
    \draw [-latex] (2,3-0.05) -- (2,2+0.05);
    \draw [-latex] (1+0.05,2) -- (2-0.05,2);
    \draw [-latex] (1+0.05,3) -- (2-0.05,3);

    \draw [fill = black] (3,2) circle (0.05cm);
    \draw                (3,1) circle (0.05cm);
    \draw [fill = black] (4,1) circle (0.05cm);
    \draw                (4,2) circle (0.05cm);

    \draw [-latex] (2+0.05,2) -- (3-0.05,2);
    \draw [-latex] (0,1) -- (3-0.05,1);

    \draw [-latex] (3+0.05,2) -- (4-0.05,2);
    \draw [-latex] (3+0.05,1) -- (4-0.05,1);
    \draw [-latex] (3,1+0.05) -- (3,2-0.05);
    \draw [-latex] (4,2-0.05) -- (4,1+0.05);

    \draw                (5,3) circle (0.05cm);
    \draw [fill = black] (6,3) circle (0.05cm);
    \draw [fill = black] (5,4) circle (0.05cm);
    \draw                (6,4) circle (0.05cm);

    \draw [-latex] (2+0.05,3) -- (5-0.05,3);
    
    \draw [-latex] (5+0.05,3) -- (6-0.05,3);
    \draw [-latex] (5+0.05,4) -- (6-0.05,4);
    \draw [-latex] (5,3+0.05) -- (5,4-0.05);
    \draw [-latex] (6,4-0.05) -- (6,3+0.05);

    \draw [domain = 4:4.5] plot (\x, {cos(180*(\x-4))});
    \draw [-latex, domain = 4.5:4.95] plot (\x, {5+cos(180*(\x-4))});

    \draw                (7,2) circle (0.05cm);
    \draw [fill = black] (8,2) circle (0.05cm);
    \draw [fill = black] (7,3) circle (0.05cm);
    \draw                (8,3) circle (0.05cm);

    \draw [-latex] (6+0.05,3) -- (7-0.05,3);
    \draw [-latex] (4+0.05,2) -- (7-0.05,2);

    \draw [-latex] (7+0.05,3) -- (8-0.05,3);
    \draw [-latex] (7+0.05,2) -- (8-0.05,2);
    \draw [-latex] (7,2+0.05) -- (7,3-0.05);
    \draw [-latex] (8,3-0.05) -- (8,2+0.05);

    \draw                (9,1) circle (0.05cm);
    \draw [fill = black] (10,1) circle (0.05cm);
    \draw [fill = black] (9,2) circle (0.05cm);
    \draw                (10,2) circle (0.05cm);

    \draw [-latex] (8+0.05,2) -- (9-0.05,2);

    \draw [domain = 6.05:7.5] plot (\x, {5-cos(180*((\x-6)/3))});
    \draw [-latex, domain = 7.5:8.95] plot (\x, {-cos(180*((\x-6)/3))});

    \draw [-latex] (9+0.05,2) -- (10-0.05,2);
    \draw [-latex] (9+0.05,1) -- (10-0.05,1);
    \draw [-latex] (9,1+0.05) -- (9,2-0.05);
    \draw [-latex] (10,2-0.05) -- (10,1+0.05);

    % Sinks
    \draw [-latex] (8+0.05,3) -- (11,3);
    \draw [-latex] (10+0.05,2) -- (11,2);
    \draw [-latex] (10+0.05,1) -- (11,1);

\end {tikzpicture}
\end {figure}

We now look at a very specific example of the quivers and Poisson structures discussed previously. The quiver $Q_{k,n}$ (or just $Q_n$ if $k=3$) is embedded in an annulus (or cylinder), 
and the case $k=3$, $n=5$ is shown below in \textbf{Figure \ref{fig:pentagram_quiver}}. In general, there are $k$ sources and $k$ sinks,
the number of square faces is $n$,
and they are connected as in the figure. More specifically, the bottom right of the $i^\mathrm{th}$ square face connects to the top left of
the $(i+1)^\mathrm{st}$, and the top right of the $i^\mathrm{th}$ square face is connected to the bottom left of the $(i+2)^\mathrm{nd}$,
with the indices understood cyclically. This is how the labels of the sinks on the right side are chosen.
The top and bottom edges of the boundary rectangle are identified, giving a cylinder. We take the cut to be the top and bottom edges
which we identify.
The sources and sinks are labeled on the left and right edges, which become the inner and outer boundary circles, respectively. As was mentioned
in the previous sections, we will sometimes consider the network on a torus by also identifying the left and right edges (up to a twist, indicated
by the labels).

As was mentioned in the previous section, the Poisson bracket depends only on the parameters $A_\circ$ and $A_\bullet$.
So we may assume that all $w_i$ and $k_i$ are zero except $w_3 = -A_\circ$ and $k_3 = -A_\bullet$.
From now on we will consider the specific choice $w_3 = k_3 = \frac{1}{2}$. Note that when two paths meet in this
network, they must come together at a black vertex, and separate at a white vertex.
This means for two paths $f$ and $g$, that $\varepsilon_x(f,g) = \pm 1$ when $f$ and $g$ touch, 
and $\varepsilon_x(f,g) = 0$ when $f$ and $g$ cross.

\bigskip

\subsection {The $x,y$ Coordinates}

\bigskip

Next, we will apply gauge transformations to this network and consider some other functions on $\mathscr{E}_Q$. We start by
defining our notation for the edge weights. We will consider the toric network obtained by glueing the boundary components
together. So the $n^\mathrm{th}$ square connects to the $1^\mathrm{st}$ square.
We label the edge weights around each square face as follows:
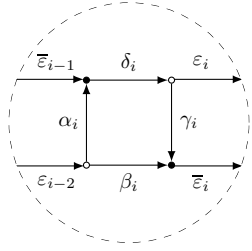
\begin {center}
\begin {tikzpicture}[scale=0.8, every node/.style={scale=0.8}]
    \clip (0,0) circle (2cm);
    \draw [dashed]     (0,0) circle (2cm);
    \draw (0.707,0.707) circle (0.05cm);
    \draw[fill=black] (0.707,-0.707) circle (0.05cm);
    \draw[fill=black] (-0.707,0.707) circle (0.05cm);
    \draw (-0.707,-0.707) circle (0.05cm);

    \draw [-latex] (-2,0.72) -- (-0.75,0.72);
    \draw [-latex] (0.75,0.72) -- (1.9,0.72);
    \draw [-latex] (0.75,-0.72) -- (1.9,-0.72);
    \draw [-latex] (-1.9,-0.72) -- (-0.74,-0.72);

    \draw [-latex] (-0.65,0.707) -- (0.65,0.707);
    \draw [-latex] (-0.65,-0.707) -- (0.65,-0.707);
    \draw [-latex] (-0.707,-0.65) -- (-0.707,0.65);
    \draw [-latex] (0.707,0.65) -- (0.707,-0.65);

    \draw (-1,0) node {$\alpha_i$};
    \draw (0,-1) node {$\beta_i$};
    \draw (1,0) node {$\gamma_i$};
    \draw (0,1) node {$\delta_i$};

    \draw (1.2,1) node {$\varepsilon_i$};
    \draw (1.2,-1) node {$\overline{\varepsilon}_i$};
    \draw (-1.2,1) node {$\overline{\varepsilon}_{i-1}$};
    \draw (-1.2,-1) node {$\varepsilon_{i-2}$};
\end {tikzpicture}
\end {center}

Because of our choice of coefficients $A_\circ=A_\bullet=-\frac{1}{2}$, the brackets between these coordinates are given by

\begin {alignat*}{2}
    \{\beta_i,\alpha_i\}                      &= \frac{1}{2} \beta_i \alpha_i                      &\hspace{1cm} \{\beta_i,\gamma_i\}       &= \frac{1}{2} \gamma_i \beta_i \\
    \{\alpha_i,\overline{\varepsilon}_{i-1}\} &= \frac{1}{2} \overline{\varepsilon}_{i-1} \alpha_i &\hspace{1cm} \{\gamma_i,\varepsilon_i\} &= \frac{1}{2} \gamma_i \varepsilon_i \\
\end {alignat*}

After applying guage transformations at the corners of each square face, we can obtain the following edge weights:

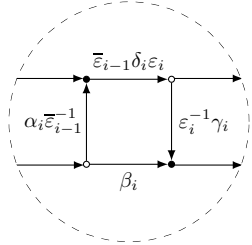
\begin {center}
\begin {tikzpicture}[scale=0.8, every node/.style={scale=0.8}]
    \clip (0,0) circle (2cm);
    \draw [dashed]     (0,0) circle (2cm);
    \draw (0.707,0.707) circle (0.05cm);
    \draw[fill=black] (0.707,-0.707) circle (0.05cm);
    \draw[fill=black] (-0.707,0.707) circle (0.05cm);
    \draw (-0.707,-0.707) circle (0.05cm);

    \draw [-latex] (-2,0.72) -- (-0.75,0.72);
    \draw [-latex] (0.75,0.72) -- (1.9,0.72);
    \draw [-latex] (0.75,-0.72) -- (1.9,-0.72);
    \draw [-latex] (-1.9,-0.72) -- (-0.74,-0.72);

    \draw [-latex] (-0.65,0.707) -- (0.65,0.707);
    \draw [-latex] (-0.65,-0.707) -- (0.65,-0.707);
    \draw [-latex] (-0.707,-0.65) -- (-0.707,0.65);
    \draw [-latex] (0.707,0.65) -- (0.707,-0.65);

    \draw (-1.25,0) node {$\alpha_i \overline{\varepsilon}_{i-1}^{-1}$};
    \draw (0,-1) node {$\beta_i$};
    \draw (1.25,0) node {$\varepsilon_i^{-1} \gamma_i$};
    \draw (0,1) node {$\overline{\varepsilon}_{i-1} \delta_i \varepsilon_i$};

\end {tikzpicture}
\end {center}

The effect of these gauge transformations is that all the edges other than those bounding the square faces have been set to $1$.
These new weights, which are monomials in the original weights, we will call by $a_i,b_i,c_i,d_i$:
\begin {alignat*}{2}
    a_i &= \alpha_i \overline{\varepsilon}_{i-1}^{-1} &\hspace{1cm} b_i &= \beta_i \\
    c_i &= \varepsilon_i^{-1} \gamma_i &\hspace{1cm} d_i &= \overline{\varepsilon}_{i-1} \delta_i \varepsilon_i
\end {alignat*}
It is easily checked that the Poisson brackets of these monomials are given by
\begin {alignat*}{2}
    \{b_i,a_i\} &= \frac{1}{2} b_ia_i &\hspace{1cm} \{b_i,c_i\} &= \frac{1}{2} b_ic_i \\
    \{a_i,d_i\} &= \frac{1}{2} a_id_i &\hspace{1cm} \{c_i,d_i\} &= \frac{1}{2} c_id_i
\end {alignat*}

\begin {figure}[h!]
\centering
\caption {The $n=5$ quiver after gauge transformations}
\label {fig:xy_quiver}
\begin {tikzpicture}[scale=0.9, every node/.style={scale=0.9}]
    % The rectangle
    \draw (0,0) -- (11,0) -- (11,5) -- (0,5) -- cycle;

    % The Labels
    \draw (-0.2,1) node {$3$};
    \draw (-0.2,2) node {$2$};
    \draw (-0.2,3) node {$1$};
    \draw (11.2,1) node {$1$};
    \draw (11.2,2) node {$3$};
    \draw (11.2,3) node {$2$};

    \draw [fill = black] (1,3) circle (0.05cm);
    \draw                (1,2) circle (0.05cm);
    \draw [fill = black] (2,2) circle (0.05cm);
    \draw                (2,3) circle (0.05cm);

    \draw [-latex] (0,3) -- (1-0.05,3);
    \draw [-latex] (0,2) -- (1-0.05,2);

    % Square 1

    \draw [-latex] (1,2+0.05) -- (1,3-0.05);
    \draw [-latex] (2,3-0.05) -- (2,2+0.05);
    \draw [-latex] (1+0.05,2) -- (2-0.05,2);
    \draw [-latex] (1+0.05,3) -- (2-0.05,3);

    \draw [fill = black] (3,2) circle (0.05cm);
    \draw                (3,1) circle (0.05cm);
    \draw [fill = black] (4,1) circle (0.05cm);
    \draw                (4,2) circle (0.05cm);

    \draw [-latex] (2+0.05,2) -- (3-0.05,2);
    \draw [-latex] (0,1) -- (3-0.05,1);

    \draw (0.8,2.5) node {$x_1$};
    \draw (1.5,1.8) node {$y_1$};

    % Square 2

    \draw [-latex] (3+0.05,2) -- (4-0.05,2);
    \draw [-latex] (3+0.05,1) -- (4-0.05,1);
    \draw [-latex] (3,1+0.05) -- (3,2-0.05);
    \draw [-latex] (4,2-0.05) -- (4,1+0.05);

    \draw                (5,3) circle (0.05cm);
    \draw [fill = black] (6,3) circle (0.05cm);
    \draw [fill = black] (5,4) circle (0.05cm);
    \draw                (6,4) circle (0.05cm);

    \draw [-latex] (2+0.05,3) -- (5-0.05,3);

    \draw (2.8,1.5) node {$x_2$};
    \draw (3.5,0.8) node {$y_2$};

    % Square 3
    
    \draw [-latex] (5+0.05,3) -- (6-0.05,3);
    \draw [-latex] (5+0.05,4) -- (6-0.05,4);
    \draw [-latex] (5,3+0.05) -- (5,4-0.05);
    \draw [-latex] (6,4-0.05) -- (6,3+0.05);

    \draw [domain = 4:4.5] plot (\x, {cos(180*(\x-4))});
    \draw [-latex, domain = 4.5:4.95] plot (\x, {5+cos(180*(\x-4))});

    \draw                (7,2) circle (0.05cm);
    \draw [fill = black] (8,2) circle (0.05cm);
    \draw [fill = black] (7,3) circle (0.05cm);
    \draw                (8,3) circle (0.05cm);

    \draw [-latex] (6+0.05,3) -- (7-0.05,3);
    \draw [-latex] (4+0.05,2) -- (7-0.05,2);

    \draw (4.8,3.5) node {$x_3$};
    \draw (5.5,2.8) node {$y_3$};

    % Square 4

    \draw [-latex] (7+0.05,3) -- (8-0.05,3);
    \draw [-latex] (7+0.05,2) -- (8-0.05,2);
    \draw [-latex] (7,2+0.05) -- (7,3-0.05);
    \draw [-latex] (8,3-0.05) -- (8,2+0.05);

    \draw                (9,1) circle (0.05cm);
    \draw [fill = black] (10,1) circle (0.05cm);
    \draw [fill = black] (9,2) circle (0.05cm);
    \draw                (10,2) circle (0.05cm);

    \draw [-latex] (8+0.05,2) -- (9-0.05,2);

    \draw [domain = 6.05:7.5] plot (\x, {5-cos(180*((\x-6)/3))});
    \draw [-latex, domain = 7.5:8.95] plot (\x, {-cos(180*((\x-6)/3))});

    \draw (6.8,2.5) node {$x_4$};
    \draw (7.5,1.8) node {$y_4$};

    % Square 5

    \draw [-latex] (9+0.05,2) -- (10-0.05,2);
    \draw [-latex] (9+0.05,1) -- (10-0.05,1);
    \draw [-latex] (9,1+0.05) -- (9,2-0.05);
    \draw [-latex] (10,2-0.05) -- (10,1+0.05);

    \draw (8.8,1.5) node {$x_5$};
    \draw (9.5,0.8) node {$y_5$};

    % Sinks
    \draw [-latex] (8+0.05,3) -- (11,3);
    \draw [-latex] (10+0.05,2) -- (11,2);
    \draw [-latex] (10+0.05,1) -- (11,1);

    \draw (10.5,3.2) node {$z$};
    \draw (10.5,2.2) node {$z$};
    \draw (10.5,1.2) node {$z$};

\end {tikzpicture}
\end {figure}
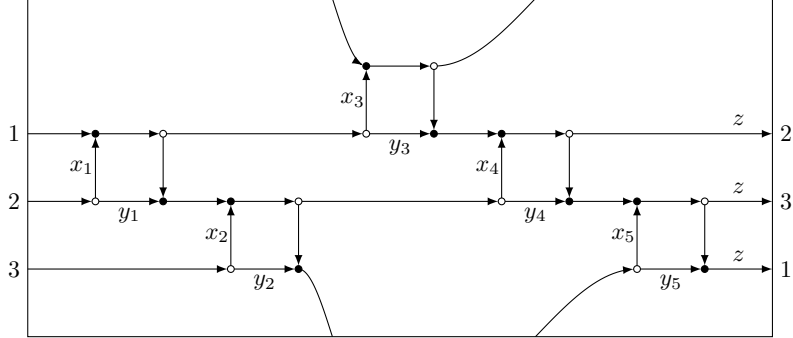

We will now apply further gauge transformations, in order to set as many edge weights as possible equal to $1$. 
It is possible, after further gauge transformations,
to set most of the weights equal to 1, except the bottom and left edges of each square face, and a few other edges. The result (again for $n=5$) 
is pictured in \textbf{Figure \ref{fig:xy_quiver}}.

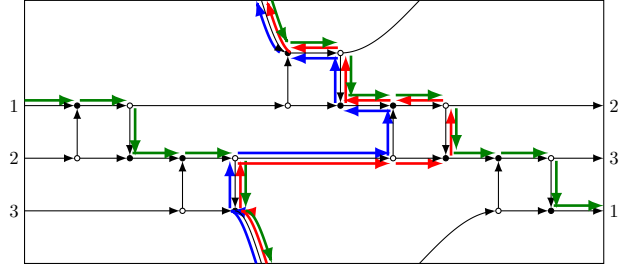
\begin {figure}[h!]
\centering
\caption {The weights $x_i,y_i,z$ represented as cycles}
\label {fig:xyz_cycles}
\begin {tikzpicture}[scale=0.7, every node/.style={scale=0.7}]
    % The rectangle
    \draw (0,0) -- (11,0) -- (11,5) -- (0,5) -- cycle;

    % The Labels
    \draw (-0.2,1) node {$3$};
    \draw (-0.2,2) node {$2$};
    \draw (-0.2,3) node {$1$};
    \draw (11.2,1) node {$1$};
    \draw (11.2,2) node {$3$};
    \draw (11.2,3) node {$2$};

    \draw [fill = black] (1,3) circle (0.05cm);
    \draw                (1,2) circle (0.05cm);
    \draw [fill = black] (2,2) circle (0.05cm);
    \draw                (2,3) circle (0.05cm);

    \draw [-latex, black!50!green, line width=1pt] (0,3.1) -- (1-0.05,3.1);

    \draw [-latex] (0,3) -- (1-0.05,3);
    \draw [-latex] (0,2) -- (1-0.05,2);

    % Square 1

    \draw [-latex] (1,2+0.05) -- (1,3-0.05);
    \draw [-latex] (2,3-0.05) -- (2,2+0.05);
    \draw [-latex] (1+0.05,2) -- (2-0.05,2);
    \draw [-latex] (1+0.05,3) -- (2-0.05,3);

    \draw [-latex, black!50!green, line width=1pt] (1+0.05,3.1) -- (2-0.05,3.1);
    \draw [-latex, black!50!green, line width=1pt] (2.1,3-0.05) -- (2.1,2+0.05);

    \draw [fill = black] (3,2) circle (0.05cm);
    \draw                (3,1) circle (0.05cm);
    \draw [fill = black] (4,1) circle (0.05cm);
    \draw                (4,2) circle (0.05cm);

    \draw [-latex] (2+0.05,2) -- (3-0.05,2);
    \draw [-latex] (0,1) -- (3-0.05,1);

    \draw [-latex, black!50!green, line width=1pt] (2+0.05,2.1) -- (3-0.05,2.1);

    % Square 2

    \draw [-latex] (3+0.05,2) -- (4-0.05,2);
    \draw [-latex] (3+0.05,1) -- (4-0.05,1);
    \draw [-latex, blue, line width=1pt] (4-0.1,1+0.05) -- (4-0.1,2-0.05);
    \draw [-latex, red, line width=1pt] (4+0.1,1+0.05) -- (4+0.1,2-0.05);
    \draw [-latex, black] (4,2-0.05) -- (4,1+0.05);
    \draw [-latex, black] (3,1+0.05) -- (3,2-0.05);

    \draw [-latex, black!50!green, line width=1pt] (3+0.05,2.1) -- (4-0.05,2.1);
    \draw [-latex, black!50!green, line width=1pt] (4+0.2,2-0.05) -- (4+0.2,1+0.05);

    \draw                (5,3) circle (0.05cm);
    \draw [fill = black] (6,3) circle (0.05cm);
    \draw [fill = black] (5,4) circle (0.05cm);
    \draw                (6,4) circle (0.05cm);

    \draw [-latex] (2+0.05,3) -- (5-0.05,3);

    % Square 3
    
    \draw [-latex] (5+0.05,3) -- (6-0.05,3);
    \draw [-latex, blue, line width=1pt] (6-0.05,4-0.1) -- (5+0.05,4-0.1);
    \draw [-latex, red, line width=1pt] (6-0.05,4+0.1) -- (5+0.05,4+0.1);
    \draw [-latex, black] (5+0.05,4) -- (6-0.05,4);
    \draw [-latex, blue, line width=1pt] (6-0.1,3+0.05) -- (6-0.1,4-0.05);
    \draw [-latex, red, line width=1pt] (6+0.1,3+0.05) -- (6+0.1,4-0.05);
    \draw [-latex, black] (6,4-0.05) -- (6,3+0.05);
    \draw [-latex, black] (5,3+0.05) -- (5,4-0.05);

    \draw [-latex, blue, line width=1pt, domain = 4.4:3.9] plot (\x, {cos(180*(\x-4+0.1))});
    \draw [-latex, red, line width=1pt, domain = 4.6:4.1] plot (\x, {cos(180*(\x-4-0.1))});
    \draw [black, domain = 4:4.5] plot (\x, {cos(180*(\x-4))});
    \draw [-latex, blue, line width=1pt, domain = 4.85:4.4] plot (\x, {5+cos(180*(\x-4+0.1))});
    \draw [-latex, red, line width=1pt, domain = 5.05:4.6] plot (\x, {5+cos(180*(\x-4-0.1))});
    \draw [-latex, black, domain = 4.5:4.95] plot (\x, {5+cos(180*(\x-4))});

    \draw [-latex, black!50!green, line width=1pt, domain = 4.2:4.7] plot (\x, {cos(180*(\x-4-0.2))});
    \draw [-latex, black!50!green, line width=1pt, domain = 4.7:5] plot (\x, {5+cos(180*(\x-4-0.2))});

    \draw                (7,2) circle (0.05cm);
    \draw [fill = black] (8,2) circle (0.05cm);
    \draw [fill = black] (7,3) circle (0.05cm);
    \draw                (8,3) circle (0.05cm);

    \draw [-latex, blue, line width=1pt] (7-0.05,3-0.1) -- (6+0.05,3-0.1);
    \draw [-latex, red, line width=1pt] (7-0.05,3+0.1) -- (6+0.05,3+0.1);
    \draw [-latex, black] (6+0.05,3) -- (7-0.05,3);
    \draw [-latex, blue, line width=1pt] (4+0.05,2+0.1) -- (7-0.05,2+0.1);
    \draw [-latex, red, line width=1pt] (4+0.05,2-0.1) -- (7-0.05,2-0.1);
    \draw [-latex, black] (4+0.05,2) -- (7-0.05,2);

    \draw [-latex, black!50!green, line width=1pt] (5+0.05,4.2) -- (6-0.05,4.2);
    \draw [-latex, black!50!green, line width=1pt] (6.2,4-0.05) -- (6.2,3+0.15);
    \draw [-latex, black!50!green, line width=1pt] (6+0.15,3+0.2) -- (7-0.05,3+0.2);

    % Square 4

    \draw [-latex, red, line width=1pt] (8-0.05,3+0.1) -- (7+0.05,3+0.1);
    \draw [-latex] (7+0.05,3) -- (8-0.05,3);
    \draw [-latex, red, line width=1pt] (7+0.05,2-0.1) -- (8-0.05,2-0.1);
    \draw [-latex] (7+0.05,2) -- (8-0.05,2);
    \draw [-latex, blue, line width=1pt] (7-0.1,2+0.05) -- (7-0.1,3-0.05);
    \draw [-latex, black] (7,2+0.05) -- (7,3-0.05);
    \draw [-latex, red, line width=1pt] (8+0.1,2+0.05) -- (8+0.1,3-0.05);
    \draw [-latex, black] (8,3-0.05) -- (8,2+0.05);

    \draw [-latex, black!50!green, line width=1pt] (7+0.05,3.2) -- (8-0.05,3.2);
    \draw [-latex, black!50!green, line width=1pt] (8.2,3-0.05) -- (8.2,2+0.15);

    \draw                (9,1) circle (0.05cm);
    \draw [fill = black] (10,1) circle (0.05cm);
    \draw [fill = black] (9,2) circle (0.05cm);
    \draw                (10,2) circle (0.05cm);

    \draw [-latex] (8+0.05,2) -- (9-0.05,2);

    \draw [-latex, black!50!green, line width=1pt] (8+0.15,2.1) -- (9-0.05,2.1);

    \draw [domain = 6.05:7.5] plot (\x, {5-cos(180*((\x-6)/3))});
    \draw [-latex, domain = 7.5:8.95] plot (\x, {-cos(180*((\x-6)/3))});

    % Square 5

    \draw [-latex] (9+0.05,2) -- (10-0.05,2);
    \draw [-latex] (9+0.05,1) -- (10-0.05,1);
    \draw [-latex] (9,1+0.05) -- (9,2-0.05);
    \draw [-latex] (10,2-0.05) -- (10,1+0.05);

    \draw [-latex, black!50!green, line width=1pt] (9+0.05,2.1) -- (10-0.05,2.1);
    \draw [-latex, black!50!green, line width=1pt] (10.1,2-0.05) -- (10.1,1+0.05);

    % Sinks
    \draw [-latex] (8+0.05,3) -- (11,3);
    \draw [-latex] (10+0.05,2) -- (11,2);
    \draw [-latex] (10+0.05,1) -- (11,1);

    \draw [-latex, black!50!green, line width=1pt] (10+0.05,1.1) -- (11,1.1);

\end {tikzpicture}
\end {figure}

The weights in the new quiver after gauge transformations are given as follows
\begin {align*}
    x_i &= \frac{a_i}{c_{i-1}d_{i-1}c_{i-2}} \\
    y_i &= \frac{b_i}{c_id_ic_{i-1}d_{i-1}c_{i-2}} \\
    z   &= \prod_{k=1}^n d_kc_k
\end {align*}
Thinking of the quiver as being on a torus, all of these monomials represent loops. They are
depicted in \textbf{Figure \ref{fig:xyz_cycles}}.
The blue loop is $x_4$, the red loop is $y_4$, and the green loop is $z$.

\bigskip

It is an easy calculation to see that $z$ is a Casimir of the Poisson bracket, and so the bracket descends
to the level surface in $\mathscr{E}_Q$ defined by $z=1$. We will now consider just edge weights lying on this
hypersurface. In this case, we have that all edge weights except those labeled by $x_i$ and $y_i$ are equal to $1$.
Therefore, this level surface has coordinates given by $x_i$ and $y_i$, and so it is $2n$-dimensional.
We will associate these coordinates with the $x_i,y_i$ coordinates on $\mathcal{P}_n$ introduced in section 1.
Their brackets are given by

\begin {alignat*}{2}
    \{x_{i+1},x_i\} &= x_{i+1} x_i &\hspace{1cm} \{y_i,x_i\}     &= y_i x_i \\
    \{y_{i+1},y_i\} &= y_{i+1} y_i &\hspace{1cm} \{y_i,x_{i-1}\} &= y_i x_{i-1} \\
    \{y_{i+2},y_i\} &= y_{i+2} y_i &\hspace{1cm} \{x_i,y_{i-1}\} &= x_i y_{i-1} \\
                    &\phantom{=}   &\hspace{1cm} \{x_i,y_{i-2}\} &= x_i y_{i-2}
\end {alignat*}

\bigskip

\subsection {The $p,q$ Coordinates}

It will be convenient to also consider a different set of functions on $\mathscr{E}_Q$. We return to the $a,b,c,d$ coordinates
before the guage transformations which gave the $x,y$ coordinates. We define the \emph{face weights} $p_i$ and $q_i$ as the 
counter-clockwise paths around the faces of the quiver, taking inverses when the orientation disagrees:

\begin {alignat*}{2}
    p_i &= \frac{b_i}{a_ic_id_i}  &\hspace{1cm} q_i &= \frac{c_{i-2}d_{i-1}a_{i+1}}{b_i} \\
\end {alignat*}

These are related to the $x,y$ coordinates by

\begin {alignat*}{2}
    p_i &= \frac{y_i}{x_i}  &\hspace{1cm} q_i &= \frac{x_{i+1}}{y_i} \\
\end {alignat*}

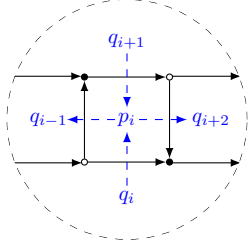
\begin {figure}[h!]
\centering
\caption {The dual quiver $Q'$}
\label {fig:dual_quiver}
\begin {tikzpicture}[scale=0.8, every node/.style={scale=0.8}]
    \clip (0,0) circle (2cm);
    \draw [dashed]     (0,0) circle (2cm);
    \draw (0.707,0.707) circle (0.05cm);
    \draw[fill=black] (0.707,-0.707) circle (0.05cm);
    \draw[fill=black] (-0.707,0.707) circle (0.05cm);
    \draw (-0.707,-0.707) circle (0.05cm);

    \draw [-latex] (-2,0.72) -- (-0.75,0.72);
    \draw [-latex] (0.75,0.72) -- (1.9,0.72);
    \draw [-latex] (0.75,-0.72) -- (1.9,-0.72);
    \draw [-latex] (-1.9,-0.72) -- (-0.74,-0.72);

    \draw [-latex] (-0.65,0.707) -- (0.65,0.707);
    \draw [-latex] (-0.65,-0.707) -- (0.65,-0.707);
    \draw [-latex] (-0.707,-0.65) -- (-0.707,0.65);
    \draw [-latex] (0.707,0.65) -- (0.707,-0.65);

    \draw (-1.3,0) node {$\color{blue} q_{i-1}$};
    \draw (0,-1.3) node {$\color{blue} q_i$};
    \draw (1.4,0)  node {$\color{blue} q_{i+2}$};
    \draw (0,1.3)  node {$\color{blue} q_{i+1}$};
    \draw (0,0)    node {$\color{blue} p_i$};

    \draw [-latex, dashed, color=blue] (0,-1.1) -- (0,-0.2);
    \draw [-latex, dashed, color=blue] (0,1.1) -- (0,0.2);
    \draw [-latex, dashed, color=blue] (0.2,0) -- (1,0);
    \draw [-latex, dashed, color=blue] (-0.2,0) -- (-1,0);
\end {tikzpicture}
\end {figure}

Define the quiver $Q'$, dual to $Q$, to have vertices corresponding to the faces of $Q$, and arrows corresponding to
arrows of $Q$ connecting vertices of different colors, directed so that white vertices are on the left and black vertices on the right.
The vertices of $Q'$ are labeled by $p_i$ and $q_i$, with arrows $p_i \to q_{i-1}$, $p_i \to q_{i+2}$, $q_i \to p_i$, and $q_{i+1} \to p_i$.
This is pictured in \textbf{Figure \ref{fig:dual_quiver}}, with the dual quiver drawn in blue.

The brackets of the $p$'s and $q$'s is ``log-canonical'', with the skew-symmetric coefficient matrix 
given by the adjacency matrix of the dual quiver $Q'$. That is,

\begin {alignat*}{2}
    \{q_i,p_i\} &= q_i p_i  &\hspace{1cm} \{q_{i+1},p_i\} &= q_{i+1} p_i \\
    \{p_i,q_{i-1}\} &= p_i q_{i-1} &\hspace{1cm} \{p_i,q_{i+2}\} &= p_i q_{i+2}
\end {alignat*}

\subsection {The Postnikov Moves and the Invariants}

\bigskip

Now we will describe a sequence of Postnikov moves which will transform this quiver (considered as being on the torus)
into an isomorphic quiver. The new edge weights obtained after this sequence
will be the expressions for the pentagram map on $\mathcal{P}_n$ in the $x,y$ coordinates given in section 1. 
We start with the quiver as in \textbf{Figure \ref{fig:xy_quiver}}, 
where all weights are $1$ except the $x_i$ and $y_i$ weights on
the bottom and sides of the square faces. We then apply the following moves, in order:
\begin {enumerate}
    \item Perform the ``square move'' at each of the $n$ square faces
    \item Perform the ``white-swap'' at each white-white edge
    \item Perform the ``black-swap'' at each black-black edge
\end {enumerate}
After this sequence of moves, the underlying directed graph is isomorphic to the one we started with. However, the
edge weights will not be of the same form. It remains to perform gauge transformations (as we did in the previous section)
so that again all weights are $1$ except the bottom and left of each square face. After these gauge transformations,
and after choosing a particular choice of graph isomorphism, the $x,y$ weights transform as
\begin {align*}
    x_i &\mapsto x_i \, \frac{x_{i+2}+y_{i+2}}{x_i+y_i} \\
    y_i &\mapsto y_{i+1} \, \frac{x_{i+3}+y_{i+3}}{x_{i+1}+y_{i+1}}
\end {align*}
This is almost (but not quite) the same as the expression for the pentagram map on $\mathcal{P}_n$ given in section 1.
But after making the change of variables $y_i \mapsto y_{i-1}$, the formulas agree. So it only differs by a shift of indices
in the $y$-variables.

\bigskip

As discussed before, the Postnikov moves and gauge transformations performed above do not change the boundary measurements. However, since the quiver
was considered on a torus, we may have had some of the vertices or edges ``wrap around'' from the right side to the left side
when doing the ``white-swap'' and ``black-swap'' moves. This corresponds to cutting a piece off the right side of the picture,
and glueing it back onto the left side. As discussed earlier, this changes the boundary measurement matrix by a re-factorization $B_1 B_2 \mapsto B_2 B_1$. 
The boundary measurement matrix $\mathscr{B}_Q(\lambda)$ is therefore changed only up to conjugation, since $B_2 B_1 = B_2(B_1B_2)B_2^{-1}$.

\bigskip

Thus the components of the characteristic polynomial $\chi(t) = \mathrm{det}(\mathrm{Id}+tB)$ are unchanged by the sequence of moves described above.
Recall that we consider the elements of the boundary measurement matrix to be the modified edge weights, which are Laurent
polynomials in the variable $\lambda$. So $\chi$ is a function of both $\lambda$ and $t$, which is polynomial in $t$, and Laurent in $\lambda$.
We denote the coefficients (which are functions of the edge weights of the quiver) by $I_{ij}$:
\[ \chi(\lambda,t) = \sum_{i,j} I_{ij} \lambda^i t^j \]
The discussion above implies that $I_{ij}$ are invariants of the pentagram map. Furthermore, we have:

\bigskip

\begin {thm} \cite{gsvt_16} \label{thm:commutative_integrability}
    Let $I_{ij}$ be as defined above. Then \\
    \begin {tabular}{cl}
        $(a)$ & $I_{ij}$ are invariant under the pentagram map \\[1.2ex]
        $(b)$ & $\{I_{ij},I_{k\ell}\} = 0$ for all $i,j,k,\ell$ \\[1.2ex]
        $(c)$ & The pentagram map is completely integrable
    \end {tabular}
\end {thm}

\bigskip

We will also consider the dynamics in the $p,q$ coordinates described before. 
A simple calculation using the relations between $p,q$ and $x,y$ shows that the pentagram map
transforms the $p,q$ coordinates by

\begin {align*}
    q_i &\mapsto \frac{1}{p_{i+1}} \\
    p_i &\mapsto q_{i+2} \frac{(1+p_i)(1+p_{i+3})}{(1+p_{i+1}^{-1})(1+p_{i+2}^{-1})}
\end {align*}

In these coordinates we can interpret the Postnikov moves as cluster transformations. To see this,
consider the dual quiver $Q'$ described before.
Interpret the $p,q$ variables as the initial cluster of a cluster algebra whose exchange matrix $B = (b_{ij})$
is defined by the dual quiver $Q'$. We use the mutation formula 
\[ 
    \mu_k(\tau_i) = \begin{cases}
                     \frac{1}{\tau_i} & \mathrm{if}~ i=k \\[1.2ex]
                     \tau_i(1+\tau_k)^{b_{ik}} & \mathrm{if}~ b_{ik}>0 \\[1.2ex]
                     \tau_i(1+\tau_k^{-1})^{b_{ik}} & \mathrm{if}~ b_{ik}<0
                 \end{cases}
\]
These are referred to as cluster $\mathcal{X}$-variables (as opposed to cluster $\mathcal{A}$-variables) by Fock and Goncharov \cite{fock_goncharov}.
They are also called ``$\tau$-coordinates'' by Gekhtman, Shapiro, and Vainshtein \cite{gsv_book}. This type of cluster algebra is called
a \emph{cluster Poisson algebra}, since it inherits a natural Poisson bracket which is ``compatible'' with mutation (see \cite{gsv_book}), given by
\[ \{\tau_i,\tau_j\} = b_{ij} \tau_i \tau_j \]
Using the dual quiver $Q'$ for the exchange matrix $B$, this cluster Poisson bracket in the $p,q$ variables coincides with the bracket on $\mathscr{E}_Q$.
Additionally, the square moves in the sequence giving the pentagram map coincide with mutations at the $p$-vertices. The sequence of mutations, $\mu$,
that mutates at each $p_i$ once gives
\begin {alignat*}{2}
    \widetilde{p}_i &= \mu(p_i) = \frac{1}{p_i}  &\hspace{1cm} 
    \widetilde{q}_i &= \mu(q_i) = q_i \, \frac{(1+p_{i+1})(1+p_{i-2})}{(1+p_i^{-1})(1+p_{i-1}^{-1})}
\end {alignat*}

The white-swap and black-swap Postnikov moves give a graph isomorphic to the original, which exchanges the $\widetilde{p}$ and $\widetilde{q}$ face weights.
That is, the new square face weights $p_i^*$ and the new octagonal face weights $q_i^*$ are
\begin {align*}
    p_i^* &= \widetilde{q}_i = q_i \, \frac{(1+p_{i+1})(1+p_{i-2})}{(1+p_i^{-1})(1+p_{i-1}^{-1})} \\
    q_i^* &= \widetilde{p}_i = \frac{1}{p_i} \\
\end {align*}

The pentagram map coincides with this formula up to the permutation of the variables which shifts the indices by $p^*_i \mapsto p^*_{i+2}$ and $q^*_i \mapsto q^*_{i+1}$,
More technically, we have

\begin {alignat*}{2}
    T(p_i) &= \widetilde{q}_{i+1}  &\hspace{1cm} 
    T(q_i) &= \widetilde{p}_{i+2}
\end {alignat*}

It is a simple calculation to verify that the brackets $\{T(p_i),T(q_j)\}$ have exactly the same form as $\{p_i,q_j\}$.

\bigskip

\section {The Grassmann Pentagram Map}

\bigskip

\subsection {Twisted Grassmann Polygons and the Pentagram Map}

\bigskip

Mari Beffa and Felipe studied a generalization of the pentagram map \cite{bf15} to the Grassmann manifold. The exposition and notation in this section is
largely borrowed from that paper, with minor variations. The projective plane $\Bbb{P}^2$ coincides with the Grassmannian $\mathrm{Gr}(1,3)$ of
$1$-dimensional subspaces of $\Bbb{R}^3$. A natural generalization would be to consider $\mathrm{Gr}(N,3N)$, in which the previous case is
just when $N=1$. In actuality, Mari Beffa and Felipe considered more generally $\mathrm{Gr}(N,kN)$, which generalizes $\Bbb{P}^{k-1}$, 
but we will focus here on the case $k=3$.

\bigskip

Consider the set $\mathrm{Mat}_{3N \times N}$ of $3N$-by-$N$ real matrices. There are two natural multiplication actions: on the left by
$\mathrm{GL}_{3N}$, and on the right by $\mathrm{GL}_N$. If we call $\mathcal{M}_N \subset \mathrm{Mat}_{3N \times N}$ the subset of rank-$N$ matrices, 
then we will identify $\mathrm{Gr}(N,3N)$ with the orbit space $\mathcal{M}_N / \mathrm{GL}_N$
by this right action, since two matrices will be equivalent if they have the same column-span. Then $\mathrm{Gr}(N,3N)$ carries a natural
left action by $\mathrm{GL}_{3N}$, induced by the action on $\mathrm{Mat}_{3N \times N}$.

\bigskip

\begin {defn}
    A \emph{twisted Grassmann $n$-gon} (or just ``polygon'' if $n$ is understood) is a
    bi-infinite sequence $(p_i)_{i \in \Bbb{Z}}$ of points in $\mathrm{Gr}(N,3N)$, with the property that  
    $p_{i+n} = M p_i$ for some projective transformation $M \in \mathrm{PGL}_{3N}$ called the \emph{monodromy}. 
\end {defn}

\medskip

As in the classical case, we will only consider twisted polygons
satsifying some nondegeneracy condition. More specifically, choose a lift $(V_i)$ of $(p_i)$ to $\mathcal{M}_N$. Define the block column
matrices $\mathbf{V}_i = (V_i V_{i+1} V_{i+2})$. We require for any $i$
that $\mathbf{V}_i$ is nonsingular, or equivalently, that the combined columns of $V_i,V_{i+1},V_{i+2}$
form a basis of $\Bbb{R}^{3N}$. If this nondegeneracy condition is satisfied, then we obtain something analogous to the linear dependence
relations given in \textbf{Equation \ref{eq:lin_dep}}. For each $i$, there are matrices $A_i,B_i,C_i \in \mathrm{GL}_N$ so that
\begin {equation} \label{eq:grassmann_lin_dep}
    V_{i+3} = V_i A_i + V_{i+1} B_i + V_{i+2} C_i
\end {equation}
For the remainder of the paper, we will use the abbreviated phrase \emph{twisted polygon} to mean a twisted Grassmann polygon. 
We will use the notation $\mathcal{GP}_{n,N}$ to denote the moduli space of twisted Grassmann $n$-gons in $\mathrm{Gr}(N,3N)$, up to the action of $\mathrm{PGL}_{3N}$.
We will use the adjective ``classical'' when we wish to distinguish the $\Bbb{P}^2$ case ($N=1$). Next we will describe the pentagram map
on Grassmann polygons.

\bigskip

Let $(p_i)$ be a twisted polygon in $\mathrm{Gr}(N,3N)$. Since $p_i$ and $p_{i+2}$ are $N$-dimensional subspaces with trivial intersection,
they span a $2N$-dimenional (codimension $N$) subspace, which we call $\mathscr{L}_i$. 
Then the intersection $\mathscr{L}_i \cap \mathscr{L}_{i+1}$ is a codimension $2N$
(dimension $N$) subspace, which is again an element of $\mathrm{Gr}(N,3N)$. We define the pentagram map to be the map that sends $(p_i)$ to $(q_i)$, where
$q_i = \mathscr{L}_i \cap \mathscr{L}_{i+1}$. We will also refer to the pentagram map as $T \colon \mathcal{GP}_{n,N} \to \mathcal{GP}_{n,N}$.
Abusing notation as in the classical case, we also write the map as if it were defined on individual vertices $T(p_i) = q_i$,
and also write the map as if it is defined on lifts $T(V_i) = W_i$ if $(V_i)$ and $(W_i)$ are lifts of $(p_i)$ and $(q_i)$.

\bigskip

\subsection {Corrugated Grassmann Polygons and Higher Pentagram Maps}

\bigskip

Analogous to the classical case, we can also define \emph{generalized higher pentagram maps} in the Grassmann case.
This is a generalization from $\mathrm{Gr}(N,3N)$ to $\mathrm{Gr}(N,kN)$. A twisted Grassmann polygon in $\mathrm{Gr}(N,kN)$
is a sequence $(p_i)$ with monodromy $M \in \mathrm{PGL}_{kN}$ so that $p_{i+n} = M p_i$ for all $i$. We denote the
set of all equivalence classes of twisted polygons (up to the action of $\mathrm{PGL}_{kN}$) by $\mathcal{GP}_{k,n,N}$.
In keeping with the notation of the previous section, we write just $\mathcal{GP}_{n,N}$ for $k=3$.

\medskip

Given a polygon $P = (p_i)$, we define the subspaces $\mathscr{L}_i := p_i + p_{i+k-1}$.
We say that a twisted Grassmann polygon is \emph{corrugated} if $\mathscr{L}_i + \mathscr{L}_{i+1}$ is a
$3N$-dimensional subspace of $\Bbb{R}^{kN}$ for all $i$. Let $\mathcal{GP}^0_{k,n,N}$ denote the set of classes of corrugated
polygons with the following additional properties, all of which hold generically: 
\begin {itemize}
    \item Any $3$ of the $4$ subspaces $p_i$, $p_{i+1}$, $p_{i+k-1}$, $p_{i+k}$ span a $3N$-dimensional subspace.
    \item $\mathscr{L}_i \cap \mathscr{L}_{i+1} \cap \mathscr{L}_{i+2} = 0$
    \item $\dim (\mathscr{L}_i \cap \mathscr{L}_{i+k-1}) = N$
\end {itemize}

\medskip

We define the generalized higher pentagram map $T \colon \mathcal{GP}^0_{k,n,N} \to \mathcal{GP}_{k,n,N}$ as follows.
The ``corrugated'' property guarantees that
\[ \dim(\mathscr{L}_i \cap \mathscr{L}_{i+1}) = \dim \mathscr{L}_i + \dim \mathscr{L}_{i+1} - \dim(\mathscr{L}_i + \mathscr{L}_{i+1}) = 2N + 2N - 3N = N \]
So $q_i := \mathscr{L}_i \cap \mathscr{L}_{i+1}$ is again an element of $\mathrm{Gr}(N,kN)$, and we define
the pentagram map to be $(p_i) \mapsto (q_i)$.

\medskip

\begin {prop}
    If $P = (p_i) \in \mathcal{GP}^0_{k,n,N}$, then its image under $T$ is corrugated.
\end {prop}
\begin {proof}
    We need to show that for the image $T(P) = Q = (q_i)$, the span $q_i + q_{i+1} + q_{i+k-1} + q_{i+k}$ is $3N$-dimensional.
    Define subspaces $\mathscr{K}_i := q_i + q_{i+1}$. Then we are trying to show that $\dim (\mathscr{K}_i + \mathscr{K}_{i+k-1}) = 3N$.
    Note that by definition, $q_i$ and $q_{i+1}$ are both subspaces of $\mathscr{L}_{i+1}$,
    and so $\mathscr{K}_i \subseteq \mathscr{L}_{i+1}$. The assumption that $\mathscr{L}_i \cap \mathscr{L}_{i+1} \cap \mathscr{L}_{i+2} = 0$
    ensures that $\dim \mathscr{K}_i = 2N$, and so in fact $\mathscr{K}_i = \mathscr{L}_{i+1}$. Then by shifting indices
    we also get that $\mathscr{K}_{i+k-1} = \mathscr{L}_{i+k}$. It is always true that $\mathscr{L}_{i+1} \cap \mathscr{L}_{i+k}$ contains
    $p_{i+k}$, and so its dimension is at least $N$. The third assumption, however, gaurantees that the intersection is exactly $p_{i+k}$.
    We then have that

    \[ \dim(\mathscr{K}_i + \mathscr{K}_{i+k-1}) = \dim \mathscr{K}_i + \dim \mathscr{K}_{i+k-1} - \dim(\mathscr{K}_i \cap \mathscr{K}_{i+k-1}) = 3N \]

\end {proof}

\bigskip

\begin {rmk}
    It is worth pointing out that the image of the map is not necessarily in the set $\mathcal{GP}^0_{k,n,N}$. The above proof shows
    that the image is corrugated, but it does not necessarily satisfy the genericity assumptions. The set on which the map may be iterated
    is thus the complement of countably many subsets of codimension 1 (and hence non-empty).
\end {rmk}

\bigskip

\subsection {Description of the Moduli Space by Networks}

\bigskip

For simplicity of presentation, we will for the remainder of the paper restrict to considering the case $k=3$,
where all polygons are automatically corrugated.
In this section, we model the moduli space $\mathcal{GP}_{n,N}$ using the networks $Q_n$ from \textbf{Section \ref{sec:model_via_networks}}.
We will use variations of this network on three surfaces: the cylinder $S^1 \times [0,1]$ (with two boundary components), the infinite cylinder $S^1 \times (0,1)$,
and the torus $S^1 \times S^1$. 

\bigskip

First, we establish some notation. As before, we let $\mathcal{M}_{N}$ be the space of all 
$3N$-by-$N$ matrices of full rank. It is an open subset of $\mathrm{Mat}_{3N \times N}$. 
In particular, for a twisted polygon $P \in \mathcal{GP}_{n,N}$, a lift of $P$ will be a point
in $\mathcal{M}_{N}^{\Bbb{Z}}$. We denote by $\mathcal{L}_{n,N} \subset \mathcal{M}_{N}^{\Bbb{Z}}$ the subset of lifts of twisted $n$-gons,
and we let $\mathcal{TL}_{n,N}$ denote the set of lifts which are ``twisted'' (the lifts satisfy $V_{i+n} = M V_i$).
Since a twisted lift is determined by the first $n$ of the $V_i$ and the monodromy, we may view
$\mathcal{TL}_{n,N}$ as an open subset of $\mathcal{M}_N^n \times \mathrm{GL}_{3N}$. The dimension is thus $3N^2(n+3)$.

\bigskip

Recall the networks
$Q_n$ defined in \textbf{Section \ref{sec:model_via_networks}}, which are embedded on the cylinder. We will also consider the infinite
networks $\widetilde{Q}_n$, embedded on the infinite cylinder, which are just infinitely many copies of $Q_n$ concatenated
together, according to the labels of the sources and sinks. Finally, we also consider the networks $\widehat{Q}_n$, on the torus,
which are obtained by glueing the boundary components of the clyinder to each other according to the labels of the sources and sinks
in $Q_n$.

\medskip

Let $\mathcal{E}_{n,N}$ be the space of all possible choices of $\mathrm{GL}_N$ weights on the edges which are the left, bottom, and right
of each square face of the network $Q_n$. Clearly there is a bijection $\mathcal{E}_{n,N} \cong \mathrm{GL}_N^{3n}$.
Similarly, let $\widetilde{\mathcal{E}}_{n,N}$ be space of $\mathrm{GL}_N$-weights on the same edges, but on the infinite network $\widetilde{Q}_n$.
Then of course $\widetilde{\mathcal{E}}_{n,N} \cong \mathrm{GL}_N^\Bbb{Z}$.

\bigskip

We now define \emph{gauge transformations} on the networks $\widetilde{Q}_n$. 
At any vertex, we may choose some $A \in \mathrm{GL}_N$, and multiply all incoming edge weights at that vertex by $A$ on the right,
and all outgoing edge weights by $A^{-1}$ on the left. Clearly, the weight of any path passing through this vertex is unchanged
by this action. We let $\mathcal{G}$ be the group generated by these transformations. If $S$ is the set of vertices, 
then clearly $\mathcal{G} \cong \mathrm{GL}_N^S$.

\medskip

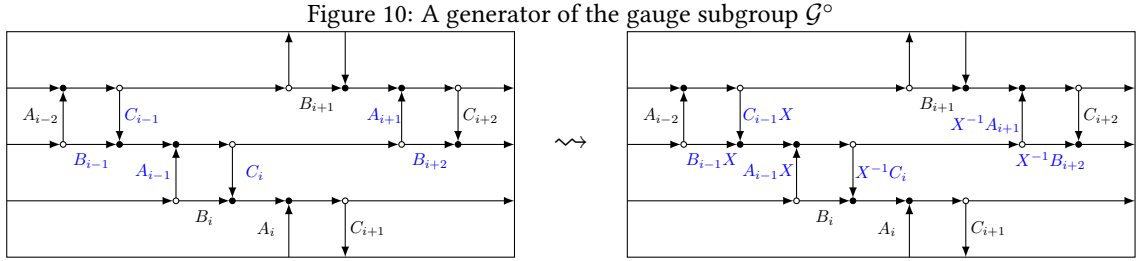
\begin {figure}[h!]
\centering
\caption{A generator of the gauge subgroup $\mathcal{G}^\circ$}
\label{fig:gauge_subgroup}
\begin {tikzpicture} [scale=0.75, every node/.style={scale=0.65}]
    \draw (0,0) -- (9,0) -- (9,4) -- (0,4) -- cycle;

    % Row 1
    \draw [-latex] (0,3) -- (0.95,3);
    \draw [fill=black] (1,3) circle (0.05cm);
    \draw [-latex] (1.05,3) -- (1.95,3);
    \draw [fill=white] (2,3) circle (0.05cm);
    \draw [-latex] (2.05,3) -- (4.95,3);
    \draw [fill=white] (5,3) circle (0.05cm);
    \draw [-latex] (5.05,3) -- (5.95,3);
    \draw [fill=black] (6,3) circle (0.05cm);
    \draw [-latex] (6.05,3) -- (6.95,3);
    \draw [fill=black] (7,3) circle (0.05cm);
    \draw [-latex] (7.05,3) -- (7.95,3);
    \draw [fill=white] (8,3) circle (0.05cm);
    \draw [-latex] (8.05,3) -- (9,3);

    % In-between 1-2
    \draw [-latex] (1,2.05) -- (1,2.95);
    \draw [-latex] (2,2.95) -- (2,2.05);
    \draw [-latex] (7,2.05) -- (7,2.95);
    \draw [-latex] (8,2.95) -- (8,2.05);

    % Row 2
    \draw [-latex] (0,2) -- (0.95,2);
    \draw [fill=white] (1,2) circle (0.05cm);
    \draw [-latex] (1.05,2) -- (1.95,2);
    \draw [fill=black] (2,2) circle (0.05cm);
    \draw [-latex] (2.05,2) -- (2.95,2);
    \draw [fill=black] (3,2) circle (0.05cm);
    \draw [-latex] (3.05,2) -- (3.95,2);
    \draw [fill=white] (4,2) circle (0.05cm);
    \draw [-latex] (4.05,2) -- (6.95,2);
    \draw [fill=white] (7,2) circle (0.05cm);
    \draw [-latex] (7.05,2) -- (7.95,2);
    \draw [fill=black] (8,2) circle (0.05cm);
    \draw [-latex] (8.05,2) -- (9,2);

    % In between 2-3
    \draw [-latex] (3,1.05) -- (3,1.95);
    \draw [-latex] (4,1.95) -- (4,1.05);

    % Row 3
    \draw [-latex] (0,1) -- (2.95,1);
    \draw [fill=white] (3,1) circle (0.05cm);
    \draw [-latex] (3.05,1) -- (3.95,1);
    \draw [fill=black] (4,1) circle (0.05cm);
    \draw [-latex] (4.05,1) -- (4.95,1);
    \draw [fill=black] (5,1) circle (0.05cm);
    \draw [-latex] (5.05,1) -- (5.95,1);
    \draw [fill=white] (6,1) circle (0.05cm);
    \draw [-latex] (6.05,1) -- (9,1);

    % In between 3-1
    \draw [-latex] (5,3.05) -- (5,4);
    \draw [-latex] (5,0) -- (5,0.95);
    \draw [-latex] (6,0.95) -- (6,0);
    \draw [-latex] (6,4) -- (6,3.05);

    % Labels
    ----------------------------------------------

    % A's
    \draw (0.6,2.5) node {$A_{i-2}$};
    \draw (2.6,1.5) node {$\color{blue}A_{i-1}$};
    \draw (4.6,0.5) node {$A_{i}$};
    \draw (6.7,2.5) node {$\color{blue}A_{i+1}$};

    % B's
    \draw (1.5,1.7) node {$\color{blue}B_{i-1}$};
    \draw (3.5,0.7) node {$B_{i}$};
    \draw (5.5,2.7) node {$B_{i+1}$};
    \draw (7.5,1.7) node {$\color{blue}B_{i+2}$};

    % C's
    \draw (2.4,2.5) node {$\color{blue}C_{i-1}$};
    \draw (4.4,1.5) node {$\color{blue}C_{i}$};
    \draw (6.4,0.5) node {$C_{i+1}$};
    \draw (8.4,2.5) node {$C_{i+2}$};

    % Arrow
    \draw (10,2) node {\huge $\rightsquigarrow$};

    % After
    \draw (11,0) -- (11+9,0) -- (11+9,4) -- (11+0,4) -- cycle;

    % Row 1
    \draw [-latex] (11+0,3)    -- (11+0.95,3);
    \draw [-latex] (11+1.05,3) -- (11+1.95,3);
    \draw [-latex] (11+2.05,3) -- (11+4.95,3);
    \draw [-latex] (11+5.05,3) -- (11+5.95,3);
    \draw [-latex] (11+6.05,3) -- (11+6.95,3);
    \draw [-latex] (11+7.05,3) -- (11+7.95,3);
    \draw [-latex] (11+8.05,3) -- (11+9,3);

    \draw [fill=black] (11+6,3) circle (0.05cm);
    \draw [fill=white] (11+5,3) circle (0.05cm);
    \draw [fill=black] (11+1,3) circle (0.05cm);
    \draw [fill=white] (11+2,3) circle (0.05cm);
    \draw [fill=black] (11+7,3) circle (0.05cm);
    \draw [fill=white] (11+8,3) circle (0.05cm);

    % In-between 1-2
    \draw [-latex] (11+1,2.05) -- (11+1,2.95);
    \draw [-latex] (11+2,2.95) -- (11+2,2.05);
    \draw [-latex] (11+7,2.05) -- (11+7,2.95);
    \draw [-latex] (11+8,2.95) -- (11+8,2.05);

    % Row 2
    \draw [-latex] (11+0,2)    -- (11+0.95,2);
    \draw [-latex] (11+1.05,2) -- (11+1.95,2);
    \draw [-latex] (11+2.05,2) -- (11+2.95,2);
    \draw [-latex] (11+3.05,2) -- (11+3.95,2);
    \draw [-latex] (11+4.05,2) -- (11+6.95,2);
    \draw [-latex] (11+7.05,2) -- (11+7.95,2);
    \draw [-latex] (11+8.05,2) -- (11+9,2);

    \draw [fill=white] (11+1,2) circle (0.05cm);
    \draw [fill=black] (11+2,2) circle (0.05cm);
    \draw [fill=black] (11+3,2) circle (0.05cm);
    \draw [fill=white] (11+4,2) circle (0.05cm);
    \draw [fill=white] (11+7,2) circle (0.05cm);
    \draw [fill=black] (11+8,2) circle (0.05cm);

    % In between 2-3
    \draw [-latex] (11+3,1.05) -- (11+3,1.95);
    \draw [-latex] (11+4,1.95) -- (11+4,1.05);

    % Row 3
    \draw [-latex] (11+0,1)    -- (11+2.95,1);
    \draw [-latex] (11+3.05,1) -- (11+3.95,1);
    \draw [-latex] (11+4.05,1) -- (11+4.95,1);
    \draw [-latex] (11+5.05,1) -- (11+5.95,1);
    \draw [-latex] (11+6.05,1) -- (11+9,1);

    \draw [fill=white] (11+3,1) circle (0.05cm);
    \draw [fill=black] (11+4,1) circle (0.05cm);
    \draw [fill=black] (11+5,1) circle (0.05cm);
    \draw [fill=white] (11+6,1) circle (0.05cm);

    % In between 3-1
    \draw [-latex] (11+5,3.05) -- (11+5,4);
    \draw [-latex] (11+5,0)    -- (11+5,0.95);
    \draw [-latex] (11+6,0.95) -- (11+6,0);
    \draw [-latex] (11+6,4)    -- (11+6,3.05);

    % Labels
    ----------------------------------------------

    % A's
    \draw (11+0.6,2.5) node {$A_{i-2}$};
    \draw (11+2.5,1.5) node {$\color{blue}A_{i-1}X$};
    \draw (11+4.6,0.5) node {$A_{i}$};
    \draw (11+6.35,2.35) node {$\color{blue}X^{-1}A_{i+1}$};

    % B's
    \draw (11+1.5,1.7) node {$\color{blue}B_{i-1}X$};
    \draw (11+3.5,0.7) node {$B_{i}$};
    \draw (11+5.5,2.7) node {$B_{i+1}$};
    \draw (11+7.5,1.7) node {$\color{blue}X^{-1}B_{i+2}$};

    % C's
    \draw (11+2.5,2.5) node {$\color{blue}C_{i-1}X$};
    \draw (11+4.5,1.5) node {$\color{blue}X^{-1}C_{i}$};
    \draw (11+6.4,0.5) node {$C_{i+1}$};
    \draw (11+8.4,2.5) node {$C_{i+2}$};
\end {tikzpicture}
\end {figure}

We would like to think of the gauge group acting on the space $\widetilde{\mathcal{E}}_{n,N}$. However, this space consists only of those
assignments of weights in which all edges other than the left, right, and bottom of each square face have weight equal to $\mathrm{Id}_N$.
The gauge transformations described above do not, in general, preserve this property. So we will instead restrict to the subgroup $\mathcal{G}^\circ \subset \mathcal{G}$
which preserves $\widetilde{\mathcal{E}}_{n,N}$. We will now describe a system of generators for $\mathcal{G}^\circ$.
The top edge of each square face has six ``nearby'' nontrivial weights -- three to the left, and three to the right.
For a matrix $X \in \mathrm{GL}_N$, the corresponding generator of $\mathcal{G}^\circ$ will multiply the three weights
to the left by $X$ on the right, and multiply the three weights to the right by $X^{-1}$ on the left. This is pictured in
\textbf{Figure \ref{fig:gauge_subgroup}}.

\medskip

Also, the subgroup $\mathcal{G}^\bullet \subset \mathcal{G}^\circ$ consisting of gauge transformations which are $n$-periodic
can be considered to act on $\mathcal{E}_{n,N}$, thought of as the space of weights for $\widehat{Q}_n$ on the torus,
since $\mathcal{E}_{n,N}$ is naturally identified with the subset of $\widetilde{\mathcal{E}}_{n,N}$ with $n$-periodic weights.

\medskip

We will now consider several quotients of $\mathcal{TL}_{n,N}$ by different group actions. First, we consider the action of $\Bbb{R}^*$
on $\mathcal{TL}_{n,N}$ by dilations. For $\lambda \in \Bbb{R}^*$, the action is given by $V_i \mapsto \lambda V_i$, and the monodromy
matrix is unchanged. This obviously does not change the polygon which the lift represents. We let $\widetilde{\mathcal{TL}}_{n,N}$ denote
the quotient by this action. Since the action is free and proper, we have
\[ \dim \widetilde{\mathcal{TL}}_{n,N} = \dim \mathcal{TL}_{n,N} - 1 = 3N^2(n+3) - 1 \]

\medskip

Now we may define an action of $\mathrm{PGL}_{3N}$ on $\widetilde{\mathcal{TL}}_{n,N}$. Let $A \in \mathrm{GL}_{3N}$ be a representative
of $\overline{A} \in \mathrm{PGL}_{3N}$, and let $(V_1,\dots,V_n,M) \in \mathcal{TL}_{n,N}$. Then $A$ acts in the usual way by $V_i \mapsto A V_i$
and $M \mapsto AMA^{-1}$. The action by conjugation is well-defined on $\mathrm{PGL}_{3N}$, since scalar matrices act trivially.
The left action $V_i \mapsto A V_i$ is also well-defined, since $\overline{A} = \overline{\lambda A}$ and $\overline{V_i} = \overline{\lambda^{-1} V_i}$.
This action is also free and proper, so
\[ \dim \widetilde{\mathcal{TL}}_{n,N} / \mathrm{PGL}_{3N} = \dim \widetilde{\mathcal{TL}}_{n,N} - \dim \mathrm{PGL}_{3N} = 3nN^2 \]

\bigskip

\begin {thm} \label{thm:lift_weight_correspondence}
    The spaces $\mathcal{TL}_{n,N}$ and $\mathcal{E}_{n,N}$ are related in the following ways: \\ \\
    \begin {tabular}{cl}
        $(a)$ & There is a bijective map $\widetilde{\mathcal{TL}}_{n,N} / \mathrm{PGL}_{3N} \to \mathcal{E}_{n,N}$ \\[1.5ex]
        $(b)$ & There is a bijective map $\widetilde{\mathcal{L}}_{n,N} / \mathrm{PGL}_{3N}  \to \widetilde{\mathcal{E}}_{n,N}$. \\[1.5ex]
        $(c)$ & Under these identifications, the actions of $\mathcal{G}^\circ$ on $\widetilde{\mathcal{E}}_{n,N}$ 
                and $\mathcal{G}^\bullet$ on $\mathcal{E}_{n,N}$ correspond to changing the lift of a fixed polygon. \\[1.5ex]
        $(d)$ & $\mathcal{GP}_{n,N} \cong \mathcal{E}_{n,N} / \mathcal{G}^\bullet \cong \widetilde{\mathcal{E}}_{n,N} / \mathcal{G}^\circ$.
    \end {tabular}
\end {thm}
\begin {proof} \ \relax \\
    $(a)$ Given a twisted lift $(V_i)$ with monodromy $M$, we have for each $i$ matrices $A_i$, $B_i$, and $C_i$ in $\mathrm{GL}_N$
    so that 
    \begin {equation} \label{eq:ABC_coefs}
        V_{i+3} = V_i A_i + V_{i+1} B_i + V_{i+2} C_i
    \end {equation}
    The fact that the lift is twisted guarantees that the $A_i$, $B_i$, and $C_i$
    are periodic. Placing $A_{i-1}$ on the left edge of each square face, $B_i$ on the bottom, and $C_i$ on the right
    defines a point in $\mathcal{E}_{n,N}$. This gives a map $\mathcal{TL}_{n,N} \to \mathcal{E}_{n,N}$. 

    \medskip
    
    It is clear that two equivalent
    polygons have the same $A_i$, $B_i$, $C_i$ coefficients for both the $\Bbb{R}^*$ and the $\mathrm{PGL}_{3N}$ actions.
    Therefore the proposed map $\widetilde{\mathcal{TL}}_{n,N}/\mathrm{PGL}_{3N} \to \mathcal{E}_{n,N}$ is well-defined.
    Now, supposing that $\overline{V}$ and $\overline{W}$ in $\widetilde{\mathcal{TL}}_{n,N}$ have the same $A_i$,
    $B_i$, $C_i$ coefficients,
    we wish to show that there is some $\overline{g} \in \mathrm{PGL}_{3N}$ so that $\overline{W} = \overline{g} \overline{V}$. 
    This will show that the map is injective.
    We first define the ``transfer matrices'' $L_i$:
    \[
        L_i := \left( \begin{array}{ccc}
                   0 & 0 & A_i \\
                   \mathrm{Id}_N & 0 & B_i \\
                   0 & \mathrm{Id}_N & C_i
               \end {array} \right)
    \]
    Now, if we form the block-column matrices $\mathbf{V}_i$, whose columns are the combined columns of $V_i$, $V_{i+1}$, and $V_{i+2}$,
    then \textbf{Equation \ref{eq:ABC_coefs}} can be written succinctly as $\mathbf{V}_{i+1} = \mathbf{V}_i L_i$.
    By the nondegeneracy assumption on lifts, $\mathrm{det}(\mathbf{V}_i) \neq 0$ for all $i$. For any two lifts $V$ and $W$
    (even of two different polygons), there are matrices $g_i \in \mathrm{GL}_{3N}$ so that $W_i = g_i V_i$ for all $i$.
    We claim that if $V$ and $W$ have the same $A_i$,$B_i$,$C_i$ coefficients, then all $g_i$ are the same. To see this,
    note that $V$ and $W$ will have the same transfer matrices $L_i$. That is, $\mathbf{V}_{i+1} = \mathbf{V}_i L_i$ 
    and $\mathbf{W}_{i+1} = \mathbf{W}_i L_i$. Let $G_i$ be the block diagonal matrices with $g_i$, $g_{i+1}$, $g_{i+2}$
    as the diagonal blocks. Then the equations $W_i = g_i V_i$ can be written as $\mathbf{W}_i = G_i \mathbf{V}_i$. Then we have
    \[ \mathbf{W}_{i+1} = \mathbf{W}_i L_i = G_i \mathbf{V}_i L_i = G_i \mathbf{V}_{i+1} \]
    But also, by definition of the $G_i$, we have $\mathbf{W}_{i+1} = G_{i+1} \mathbf{V}_{i+1}$. So it must be that
    $g_{i+1} = g_i$ for all $i$. We have shown that for two lifts $V$ and $W$ in the same fiber of the map
    $\mathcal{TL}_{n,N} \to \mathcal{E}_{n,N}$, there is some $g \in \mathrm{GL}_{3N}$ so that $W = gV$. 
    In terms of the elements of $\widetilde{\mathcal{TL}}_{n,N}$, we can say that $\overline{W} = \overline{g} \overline{V}$
    for the corresponding element $\overline{g} \in \mathrm{PGL}_{3N}$. This shows the map is injective.

    \medskip

    To see that the map is surjective, let $A_i$, $B_i$, $C_i$ be an arbitrary choice of weights in $\mathcal{E}_{n,N}$.
    We will construct a lift which maps to this choice of weights. We will define $V_1$, $V_2$, and $V_3$ so that
    the matrix $\mathbf{V}_1 = (V_1 V_2 V_3)$ is the identity matrix $\mathrm{Id}_{3N}$. We then define the rest
    of the $V_i$ by \textbf{Equation \ref{eq:ABC_coefs}}. Any other equivalent polygon is given by a different choice
    of $\mathbf{V}_1$, which differs only by the $\mathrm{PGL}_{3N}$ action, giving the same point in $\widetilde{\mathcal{TL}}_{n,N} / \mathrm{PGL}_{3N}$.
    This gives an inverse to the map described above.

    \bigskip

    $(b)$ There is the obvious map $\mathcal{L}_{n,N} \to \widetilde{\mathcal{E}}_{n,N}$ which sends a lifted polygon $V$ to the
    $A_i,B_i,C_i$ defined above. If the lift is not twisted, these are not necessarily periodic. Two lifted polygons $V$ and $W$
    have the same $A_i$, $B_i$, $C_i$ coefficients if and only if $W_i = g V_i$ for all $i$ and some fixed $g \in \mathrm{GL}_{3N}$, 
    as described above. The proof is the same as in part $(a)$. This shows that the induced map
    $\widetilde{\mathcal{L}}_{n,N} / \mathrm{PGL}_{3N} \to \widetilde{\mathcal{E}}_{n,N}$ is bijective. 

    \bigskip

    $(c)$ To see the equivalence of the gauge action of $\mathcal{G}^\circ$ and changing lifts of a polygon, we will use the following
    visualization technique. We will associate to the top edge of each square face the lifted vertex $V_i$. The left edge will
    be labelled by $A_{i-3}$, the bottom with $B_{i-2}$, and the right by $C_{i-2}$. Then the relation $V_{i+3} = V_i B_i + V_{i+1} A_i + V_{i+2} C_i$
    can be visualized in the following way. We sum over all paths from $V_j$ to $V_i$ which do not pass over another $V_k$ in between.
    We see that there are only 3 such paths, which start at $V_i$, $V_{i+1}$, and $V_{i+2}$. The weights of the paths are the matrix
    coefficients which we multiply by on the right. 

    \medskip

    The claim is that the action of the generators of $\mathcal{G}^\circ$ and $\mathcal{G}^\bullet$ described above correspond precisely to changing to a
    different lift of the same twisted polygon. So consider changing one lifted vertex $V_i$ to $\widehat{V}_i = V_i G$ for some $G \in \mathrm{GL}_N$.
    Each $V_i$ is involved in exactly four versions of \textbf{Equation \ref{eq:ABC_coefs}}. Since $V_i = \widehat{V}_i G^{-1}$,
    these four relations for the new lift become:
    \begin {align*} 
        V_{i+3} &= \widehat{V}_i (G^{-1}B_i) + V_{i+1} A_i + V_{i+2} C_i \\
        V_{i+2} &= V_{i-1} B_{i-1} + \widehat{V}_{i} (G^{-1}A_{i-1}) + V_{i+1} C_{i-1} \\
        V_{i+1} &= V_{i-2} B_{i-2} + V_{i-1} A_{i-2} + \widehat{V}_i (G^{-1}C_{i-2}) \\
        \widehat{V}_i &= V_{i-3} (B_{i-3}G) + V_{i-2} (A_{i-3}G) + V_{i-1} (C_{i-3}G)
    \end {align*}
    These are precisely how the weights change in \textbf{Figure \ref{fig:gauge_subgroup}} by taking $X=G$.

    \bigskip

    $(d)$ Recall that $\mathcal{GP}_{n,N}$ is the quotient of $\mathcal{TL}_{n,N}$
    by the actions of $\mathrm{GL}_{3N}$ and of changing the lift.
    By part $(a)$, we identify $\mathcal{E}_{n,N}$ with the space of twisted lifts, up to the actions of
    $\Bbb{R}^*$ and $\mathrm{PGL}_{3N}$. The $\Bbb{R}^*$ action is one way of changing the lift,
    and corresponds to a subgroup of $\mathcal{G}^\bullet$. This subgroup is the diagonal embedding of the scalar
    matrices $\{\lambda \mathrm{Id}_N\}$ into $\mathrm{GL}_N^n$.
    By part $(c)$, the action of $\mathcal{G}^\bullet$ on $\mathcal{E}_{n,N}$ is equivalent to chaging the lift.
    We may thus identify $\mathcal{GP}_{n,N}$ with $\mathcal{E}_{n,N} / \mathcal{G}^\bullet$.
    Similarly, we may describe
    it also as $\mathcal{GP}_{n,N} \cong \widetilde{\mathcal{E}}_{n,N}/\mathcal{G}^\circ$.
\end {proof}

\bigskip

We now use the gauge actoin to choose a convenient normalization, which we will use to coordinatize
the moduli space. The following is an analogue of \textbf{Proposition \ref{prop:normal_lift}}.

\bigskip

\begin {prop} \label{prop:grassmann_normal_lift}
    The lift $(V_i)$ of $(p_i)$ can be chosen so that $C_i = \mathrm{Id}_N$ for all $i$.
\end {prop}
\begin {proof}
    We start by choosing a twisted lift $V_i$. This guarantees that the $A_i$, $B_i$, and $C_i$ are periodic.
    Any other lift $\widehat{V}_i$ differs from $V_i$ by a change of basis. That is, there are $G_i \in \mathrm{GL}_N$ so that $\widehat{V}_i = V_i G_i$.
    Equivalently, $V_i = \widehat{V}_i G_i^{-1}$. Substituting this into \textbf{Equation \ref{eq:grassmann_lin_dep}}, we obtain
    \[ \widehat{V}_{i+3} = \widehat{V}_i \left(G_i^{-1} A_i G_{i+3}\right) + \widehat{V}_{i+1} \left(G_{i+1}^{-1}B_iG_{i+3}\right) + \widehat{V}_{i+2} \left(G_{i+2}^{-1}C_iG_{i+3}\right) \]
    We want to prove that we can always ensure that $G_{i+2}^{-1} C_i G_{i+3} = \mathrm{Id}$. In other words, we desire $G_{i+3} = C_i^{-1} G_{i+2}$.
    Re-indexing gives $G_{i+1} = C_{i-2}^{-1} G_i$. Now we may choose for instance $G_0 = \mathrm{Id}$, and determine the rest by this recurrence.

    \medskip

    Equivalently, we may think of this in terms of the weights in $\widetilde{\mathcal{E}}_{n,N}$ and the gauge action
    of $\mathcal{G}^\circ$. Solving the recurrence above means first choosing a position (a square face) in the network,
    and then applying elementary gauge transformations to the left and right to cancel the weights on the right
    edges of each square face.
\end {proof}

\medskip

The lift obtained in the above proposition is identified with a choice of weights in $\widetilde{\mathcal{E}}_{n,N}$
for which the only non-trivial weights are on the left and bottom edges of each square face. We call the weights
on the left edges $X_i$, and the ones on the bottom $Y_i$, 
so that \textbf{Equation \ref{eq:grassmann_lin_dep}} becomes
\begin {equation} \label{eq:grassmann:lin_dep_xy}
    V_{i+3} = V_i Y_i + V_{i+1} X_i + V_{i+2} 
\end {equation}

\medskip

\begin {rmk}
    It is worth pointing out that although in the proof of the above proposition, we started with a twisted lift,
    the resulting weights in $\widetilde{\mathcal{E}}_{n,N}$ are \emph{not} periodic, meaning the lift guaranteed
    by the proposition is not twisted. However, the $X_i$, $Y_i$ weights are \emph{almost} periodic in the
    following sense. Define the matrix which is the cyclic product of the $C_i$'s:
    \[ Z := C_{n-1} C_{n} C_1 \cdots C_{n-2} \]
    Then shifting the indices by $n$ corresponds to conjugation by $Z$. That is,
    \begin {align*}
        X_{i+n} &= Z^{-1} \, X_i \, Z \\
        Y_{i+n} &= Z^{-1} \, Y_i \, Z
    \end {align*}
    So although the weights themselves are not periodic, their conjugacy classes \emph{are} periodic.
\end {rmk}

\medskip

As a consequence of this normalization, we have a new model for the moduli space.

\bigskip

\begin {thm} \label{thm:moduli_space_identification}
    There is a bijection between the moduli space $\mathcal{GP}_{n,N}$ and $\mathrm{GL}_N^{2n+1}/\mathrm{Ad} \, \mathrm{GL}_N$,
    the space of $(2n+1)$-tuples of matrices up to simultaneous conjugation. Therefore we have
    that $\dim \mathcal{GP}_{n,N} = 2nN^2 + 1$.
\end {thm}
\begin {proof}
    We would like to say that the collection of matrices $X_i,Y_i,Z$ for $1 \leq i \leq n$ defines a map $\mathcal{GP}_{n,N} \to \mathrm{GL}_N^{2n+1}$.
    However, this is not quite well-defined. This depended on our choice of lift. The gauge action does not always preserve
    this collection of $X_i,Y_i,Z$ matrices, but in general changes them all by simultaneous conjugation. Therefore the induced
    map $\mathcal{GP}_{n,N} \to \mathrm{GL}_N^{2n+1}/\mathrm{GL}_N$ is well-defined.

    \medskip

    The map is injective, since the $X_i$'s, $Y_i$'s and $Z$ determine the class of a polygon. To see this, 
    choose any elements for $V_1$, $V_2$, and $V_3$. We may for instance choose them so that $\mathbf{V}_1 = \mathrm{Id}_{3N}$.
    Then the rest of the $V_i$ are determined by the relations $V_{i+3} = V_i Y_i + V_{i+1}X_i + V_{i+2}$.
    However, the $X_i$'s, $Y_i$'s and $Z$ are only defined up to simultaneous conjugation. Choosing some matrix $A$
    and conjugating all the $X_i$, $Y_i$ and $Z$ by $A$, however, changes the lift by $V_i \mapsto V_i A$ for all $i$.
    This just changes to a different lift of the same polygon, and so the choice of representatives of the $X_i$, $Y_i$,
    and $Z$ do not matter. Also, any other choice of the initial $V_1$, $V_2$, $V_3$ differs by the left action
    of $\mathrm{PGL}_{3N}$, and so gives an equivalent polygon. This gives an inverse to the map 
    $\mathcal{GP}_{n,N} \to \mathrm{GL}_{N}^{2n+1}/\mathrm{Ad} \, \mathrm{GL}_N$. 
\end {proof}

\bigskip

\begin {rmk}
    The space $\mathrm{GL}_N^{2n+1}/\mathrm{Ad} \, \mathrm{GL}_N$ parameterizes isomorphism classes of $N$-dimensional
    representations of the free associative algebra on $2n+1$ generators. Later, we will think of the matrices
    $X_i$, $Y_i$, and $Z$ as formal non-commutative variables (generators of a free algebra), and we will identify
    our moduli space $\mathcal{GP}_{n,N}$ with the moduli space of $N$-dimensional representations of this algebra.
\end {rmk}

\bigskip

With this description of the moduli space, we can now give some convenient coordinates.
There is the following classical theorem of Procesi:

\medskip

\begin {thm} \cite{procesi_76} 
    Any polynomial invariant of an $n$-tuple of matrices $(A_1,\dots,A_n)$,
    under the action of simultaneous conjugation, is a polynomial in the functions
    $\mathrm{tr}(A_{i_1} \cdots A_{i_k})$, where $A_{i_1} \cdots A_{i_k}$ ranges over
    all non-commutative monomials in the matrices $A_1,\dots,A_n$.
\end {thm}

\medskip

As mentioned above, we identify $\mathcal{GP}_{n,N}$ with the quotient $\mathrm{GL}_N^{2n+1}/\mathrm{Ad} \, \mathrm{GL}_N$.
We consider the coordinate ring to be the subring of invariants of polynomial functions on $\mathrm{GL}_N^{2n+1}$.
Procesi's theorem says that this coordinate ring is generated by the traces of monomials in $2n+1$ matrices.
As coordinates
for the moduli space $\mathcal{GP}_{n,N}$, we may therefore take any algebraically independent family of traces 
of size $2nN^2 + 1$.

\bigskip

\begin {ex}
    Consider the case of $\mathrm{GL}_2^3 / \mathrm{Ad} \, \mathrm{GL}_2$, the space of triples $(X,Y,Z)$ of $2 \times 2$
    matrices up to simultaneous conjugation. A theorem of Sibiirski \cite{sibirskii_68}
    says that the ring of invariants is minimally generated by the ten functions
    \[ \begin {array}{ccc}
            \mathrm{tr}(X)   & \mathrm{tr}(Y)   & \mathrm{tr}(Z) \\
            \mathrm{tr}(X^2) & \mathrm{tr}(Y^2) & \mathrm{tr}(Z^2) \\
            \mathrm{tr}(XY)  & \mathrm{tr}(XZ)  & \mathrm{tr}(YZ) \\
                             & \mathrm{tr}(XYZ) &
        \end {array}
    \]
    However, the dimension of the space is $9$. We claim that the first $9$ (all but $\mathrm{tr}(XYZ)$)
    can be taken as a system of local coordinates. The tangent space to $\mathrm{GL}_2^3$ may be identified
    with $\mathrm{Mat}_2^3$. The gradients of the ten functions above are given by
    \[
        \begin {array}{ccc}
            (\mathrm{Id}_2,0,0) & (0,\mathrm{Id}_2,0) & (0,0,\mathrm{Id}_2) \\
            (X^T,0,0)           & (0,Y^T,0)           & (0,0,Z^T) \\
            (Y^T,X^T,0)         & (Z^T,0,X^T)         & (0,Z^T,Y^T) \\
                                & ((YZ)^T,(ZX)^T,(XY)^T) &
        \end {array}
    \]
    The first $9$ are clearly linearly independent for generic choices of $X$, $Y$, and $Z$.
    The fact that the last is a linear combination of the first $9$ is equivalent to the statement that
    there are coefficients $c_1,\dots,c_9$ for which
    \begin {align*}
        ZY &= c_1 \mathrm{Id}_2 + c_4 X + c_7 Y + c_8 Z \\
        XZ &= c_2 \mathrm{Id}_2 + c_7 X + c_5 Y + c_9 Z \\
        YX &= c_3 \mathrm{Id}_2 + c_8 X + c_9 Y + c_6 Z
    \end {align*}
    Of course, if $X$, $Y$, and $Z$ are linearly independent, then the set $\{\mathrm{Id}_2,X,Y,Z\}$
    is a linear basis of the set of $2 \times 2$ matrices, and we can expand the products $YX$, $ZY$,
    and $XZ$ in terms of this basis. The claim above is that when we do this, three pairs of coefficients
    will be equal --- the $X$-coefficient of $YX$ is equal to the $Z$-coefficient of $ZY$, etc.
    The fact that these pairs of coefficients will be equal follows from the Jacobi identity:
    \[ [X,[Y,Z]] + [Z,[X,Y]] + [Y,[Z,X]] = 0 \]
\end {ex}

\bigskip

\subsection {Expression for the Pentagram Map}

\bigskip

In this section, we see how the pentagram map transforms the $X_i$, $Y_i$, and $Z$ matrices.
The map is only defined generically, since we need to assume certain matrices are invertible.

\begin {prop}
    The pentagram map transforms the $X_i$, $Y_i$, and $Z$ by
    \begin {align*}
        X_i &\mapsto (X_i+Y_{i+1})^{-1}X_i(X_{i+2}+Y_{i+3}) \\
        Y_i &\mapsto (X_i+Y_{i+1})^{-1}Y_{i+1}(X_{i+2}+Y_{i+3}) \\
        Z   &\mapsto Z
    \end {align*}
\end {prop}
\begin {proof}
    The proof is essentially the same as that of \textbf{Proposition \ref{prop:pentagram_in_coords}}.

    By re-arranging \textbf{Equation \ref{eq:grassmann:lin_dep_xy}}, we see that for a lift of the image, we may take for $T(V_i)$ the subspace spanned by the columns of
    \begin {equation} \label{eq:grassmann_pent_lifts} 
        V_{i+3} - V_{i+1} X_i = V_i Y_i + V_{i+2} 
    \end {equation}

    Since the lift of the image under the pentagram map is again a twisted polygon, it satisfies \textbf{Equation \ref{eq:grassmann_lin_dep}}:
    \[ T(V_{i+3}) = T(V_i) A_i + T(V_{i+1}) B_i + T(V_{i+2}) C_i \]
    We substitute the expressions from \textbf{Equation \ref{eq:grassmann_pent_lifts}} into the equation above, using either the left-hand or right-hand
    side, according to the same convention used in the proof of \textbf{Proposition \ref{prop:pentagram_in_coords}}. Doing so, we obtain
    \[ V_{i+5} + V_{i+3}Y_{i+3} = V_{i+1} \left( Y_{i+1}B_i - X_iA_i \right) + V_{i+3} \left( A_i + B_i - X_{i+2}C_i \right) + V_{i+5}C_i \]
    Comparing coefficients of each $V_i$, we conclude that $C_i = \mathrm{Id}$, and that
    \begin {align*}
        X_i A_i &= Y_{i+1}B_i \\
        Y_{i+3} + X_{i+2} &= A_i + B_i
    \end {align*}
    Assuming that each $X_i+Y_{i+1}$ is invertible, these equations can be solved for $A_i$ and $B_i$ to give the desired result.
    To see that $Z$ is unchanged, it is a simple check that a shift of indices by $n$ still corresponds to conjugation by the same $Z$.
\end {proof}

Notice that the expressions in the proposition are a non-commutative version of the expressions obtained in the classical case, in \textbf{Proposition \ref{prop:pentagram_in_coords}}.
For the remainder of the paper, we attempt to generalize much of what was done in the previous expository sections to noncommutative variables,
using elements of a noncommutative ring in place of the commutative coordinates on $\mathscr{E}_Q$, and we introduce a noncommutative
Poisson structure which mimics the classical counterpart in \cite{gsvt_16}.

\bigskip

\section {Non-Commutative Poisson Structures}

\bigskip

\subsection {Double Brackets}

\bigskip

In this section, we discuss the basic constructions and definitions introduced by Van den Bergh \cite{bergh_08}. Throughout, $\Bbb{K}$ is a field,
and $A$ is an associative $\Bbb{K}$-algebra. Unadorned tensor products are assumed to be over $\Bbb{K}$.
We will assume the characteristic of $\Bbb{K}$ is zero, but we do not necessarily assume 
it is algebraically closed.

\bigskip

If $A$ is an associative algebra, then $A^{\otimes n}$ is an $A$-bimodule in the obvious way:
\[ x \cdot (a_1 \otimes a_2 \otimes \cdots \otimes a_{n-1} \otimes a_n) \cdot y = xa_1 \otimes a_2 \otimes \cdots \otimes a_{n-1} \otimes a_ny \]
We refer to this as the \emph{outer bimodule structure} on $A^{\otimes n}$.

\bigskip

\begin {defn}
    A \emph{double bracket} on $A$ is a $\Bbb{K}$-bilinear map
    \[ \db{-,-} \colon A \times A \to A \otimes A \]
    which satisfies:
    \begin {enumerate}
        \item[$(1)$] $\db{-,-}$ is a derivation in the second argument with respect to the outer bimodule structure:
              \[ \db{a,bc} = \db{a,b}(1 \otimes c) + (b \otimes 1)\db{a,c} \]
        \item[$(2)$] $\db{b,a} = -\db{a,b}^\tau$, where $(x \otimes y)^\tau := y \otimes x$
    \end {enumerate}
\end {defn}

\bigskip

Note that these properties imply that $\db{-,-}$ also satisfies
\[ \db{ab,c} = (1 \otimes a)\db{b,c} + \db{a,c}(b \otimes 1) \]

\bigskip

\begin {defn}
    A double bracket $\db{-,-}$ is called a \emph{double Poisson bracket} if it additionally satisfies a version
    of the Jacobi identity:
    \[ 0 = \sum_{k=0}^2 \sigma^k \circ (\db{-,-} \otimes \mathrm{Id}) \circ (\mathrm{Id} \otimes \db{-,-}) \circ \sigma^{-k} \]
    The right-hand side is an operator on $A \otimes A \otimes A$, and $\sigma$ is the permutation operator which sends
    $x \otimes y \otimes z$ to $z \otimes x \otimes y$.
\end {defn}

\bigskip

For an algebra $A$, let $\mu \colon A \otimes A \to A$ be the multilpication map. If $A$ has a double bracket $\db{-,-}$,
then we define another operation $\{-,-\} \colon A \times A \to A$ by composing with $\mu$:
\[ \{a,b\} := \mu(\db{a,b}) \]

\bigskip

\begin {prop} \label {prop:single_bracket_properties} \cite{bergh_08}
    Suppose $\db{-,-}$ is a double bracket on $A$ (not necessarily a double Poisson bracket).
    The induced bracket $\{-,-\}$ has the following properties
    \begin {enumerate}
        \item $\{a,bc\} = \{a,b\}c + b\{a,c\}$ \hfill (Leibniz in the $2^\mathrm{nd}$ argument)
        \item $\{ab,c\} = \{ba,c\}$ \hfill (cyclic in the $1^\mathrm{st}$ argument)
        \item $\{a,b\} \equiv -\{b,a\} \; \mathrm{mod} \; [A,A]$ \hfill (skew mod commutators)
    \end {enumerate}
\end {prop}
\begin {proof}
    Suppose that $\db{a,b} = \sum_i \omega_i \otimes \bar{\omega}_i$ and $\db{a,c} = \sum_i \eta_i \otimes \bar{\eta}_i$.
    Then using the Leibniz rule for $\db{-,-}$ in the second argument, and composing with multiplication, we get
    \[ \{a,bc\} = \sum_i \omega \bar{\omega}_i c +  \sum_i b \eta_i \bar{\eta}_i =  \{a,b\}c + b \{a,c\} \]
    This proves the first identity. The third identity follows from the fact that
    \[ \{a,b\} + \{b,a\} = \sum_i [\omega_i,\bar{\omega_i}] \]
    For the second identity, suppose that $\db{b,c} = \sum_i \zeta_i \otimes \bar{\zeta}_i$.
    Then on the one hand, we have
    \[ \db{ab,c} = \sum_i \eta_i b \otimes \bar{\eta}_i + \sum_i \zeta_i \otimes a \bar{\zeta}_i \]
    On the other hand, we have
    \[ \db{ba,c} = \sum_i \eta_i \otimes b \bar{\eta}_i + \sum_i \zeta_i a \otimes \bar{\zeta}_i \]
    Obviously both will be the same after composing with the multiplication map.
\end {proof}

\bigskip

\begin {defn}
For an associative algebra $A$, define the \emph{cyclic space}, denoted $A^\natural$, 
to be the vector space quotient $A/[A,A]$, by the linear span of all commutators.
\end {defn}

\bigskip

Note that since $[A,A]$ is not in general an ideal, this is not necessarily an algebra.
For a vector subspace $V \leq A$, we let $V^\natural$ denote its image under the projection.
Similarly, for an element $a \in A$, we use the notation $a^\natural$ to denote
its image under the projection. Note that conjugate elements are equivalent, since $xyx^{-1} - y = [xy,x^{-1}]$.
From this it follows that in $A^\natural$, monomials are equivalent up to cyclic permutation,
since for any $x_1,\dots,x_n \in A$, we have $x_nx_1\cdots x_{n-1} = x_n(x_1\cdots x_n)x_n^{-1}$.

\medskip

The following follows easily from the properties above.

\bigskip

\begin {prop} \cite{bergh_08}
    Suppose $\db{-,-}$ is a double bracket on $A$. Then there is 
    a well-defined bilinear skew-symmetric map 
    $\left<-,-\right> \colon A^\natural \times A^\natural \to A^\natural$, given by
    \[ \left<a^\natural,b^\natural\right> := \{a,b\}^\natural \]
\end {prop}

The next few results indicate the extra structure which $\{-,-\}$ and $\left<-,-\right>$ inherit 
if $\db{-,-}$ additionally satisfies the double Jacobi identity.

\bigskip

\begin {prop} \cite{bergh_08}
    \label{prop:double_poisson_is_H_0}
    Suppose that $\db{-,-}$ is a double Poisson bracket. Then
    \begin {enumerate}
    \item $\{-,-\}$ is a ``Loday bracket'' on $A$. That is, it satisfies a version of the Jacobi identity:
          \[ \{a,\{b,c\}\} = \{\{a,b\},c\} + \{b,\{a,c\}\} \]
    \item $\left<-,-\right>$ is a Lie bracket on $A^\natural$.
    \end {enumerate}
\end {prop}

\bigskip

A more general notion is given by what William Crawley-Boevey calls an $H_0$-Poisson structure:

\medskip

\begin {defn} \cite{crawley-boevey_11}
    An \emph{$H_0$-Poisson strcuture} on an associative algebra $A$ is a Lie bracket $[-,-]$ on $A^\natural$
    such that each $[a,-]$ is induced by a derivation of $A$.
\end {defn}

\medskip

In particular, the induced bracket $\left<-,-\right>$ of a double Poisson bracket $\db{-,-}$ is an $H_0$-Poisson 
structure by \textbf{Proposition \ref{prop:single_bracket_properties}}
and \textbf{Proposition \ref{prop:double_poisson_is_H_0}}.
However, there are $H_0$-Poisson
structures which do not arise from double Poisson brackets. Later, we will define a double bracket which is \emph{not}
double Poisson, but still has the property that the induced bracket $\left<-,-\right>$ on $A^\natural$ is a Lie bracket.
This will be an example of an $H_0$-Poisson structure which comes from a double bracket, but not a double Poisson bracket.

\bigskip

There is a theorem of Crawley-Boevey which says that $H_0$-Poisson structures on $A$ induce usual ``commutative'' Poisson
structures on the representation space $\mathrm{Rep}_N(A) := \mathrm{Hom}(A,\mathrm{Mat}_N)/\mathrm{Ad} \, \mathrm{GL}_N$.

\medskip

\begin {thm} \cite{crawley-boevey_11}
    If $A$ is an associative algebra, and the bracket $\left<-,-\right>$ on $A^\natural$ defines an $H_0$-Poisson 
    structure, then there is a unique Poisson structure on $\mathrm{Rep}_N(A)$ such that the trace map
    is a Lie algebra homomorphism. That is,
    \[ \{\mathrm{tr}(a),\mathrm{tr}(b)\} = \mathrm{tr} \left<a^\natural,b^\natural\right> \]
\end {thm}

\bigskip

With this theorem in mind, we now describe the plan of the remainder of the paper. Recall that by
\textbf{Theorem \ref{thm:moduli_space_identification}}, our moduli space $\mathcal{GP}_{n,N}$ is identified
with $\mathrm{GL}_N^{2n+1}/\mathrm{Ad} \, \mathrm{GL}_N$. We will consider our algebra $A$ to be a
free algebra with $2n+1$ generators, corresponding to the $X_i$, $Y_i$, and $Z$. Since any homomorphism
with domain $A$ is determined by the images of the generators, we see that 
$\mathrm{Rep}_N(A) \cong \mathrm{Mat}_N^{2n+1}/\mathrm{Ad} \, \mathrm{GL}_N$. We may therefore interpret $\mathrm{tr}(X_i)$,
$\mathrm{tr}(Y_i)$, and $\mathrm{tr}(Z)$ as functions on $\mathcal{GP}_{n,N}$. We will spend the next couple
sections defining non-commutative Poisson structures on this free algebra. The theorem above will then induce
a Poisson structure on our moduli space.

\bigskip

\section {Non-Commutative Networks and Double Brackets}

\bigskip

\subsection {Weighted Directed Fat Graphs (Revisited)}

\bigskip

Let $Q=(Q_0,Q_1)$ be a quiver as in section 2, which is embedded in an annulus, with all sources on the inner boundary circle,
all sinks on the outer boundary circle, and all internal vertices trivalent. Define the algebra $A = A_Q = \Bbb{Q}\left<\alpha ~|~ \alpha \in Q_1\right>$
to be the free associative algebra generated by the arrows of the quiver. 
Assuming, as before, that $Q$ is acyclic, we can define the boundary measurement matrix $\mathscr{B}_Q(\lambda) = (b_{ij}(\lambda))$,
where $b_{ij}(\lambda)$ is the sum of the weights of all paths from source $i$ to sink $j$. Here, the weight of a path is the
product of the weights, in order from left to right. We again assign an indeterminate $\lambda$ to the cut, and consider the boundary
measurements to be Laurent polynomials in $A[\lambda^\pm]$, where $\lambda$ commutes with elements in $A$.

\medskip

We can define \emph{non-commutative gauge transformations}. For any Laurent monomial $t$ in the generators of $A$,
we can multiply all incoming arrows at one vertex by $t$ on the right, and multiply all outgoing arrows at that vertex by $t^{-1}$ on the left.
This will obviously not change the boundary measurements.

\medskip

\begin {figure}[h!]
\centering
\caption {non-commutative Postnikov moves}
\label {fig:noncomm_postnikov_moves}
\begin {tikzpicture}[scale=0.8, every node/.style={scale=0.8}]
    % Label (I)
    \draw (-4,0) node {( I )};

    % Before
    \draw [dashed]     (0,0) circle (2cm);
    \draw (0.707,0.707) circle (0.05cm);
    \draw[fill=black] (0.707,-0.707) circle (0.05cm);
    \draw[fill=black] (-0.707,0.707) circle (0.05cm);
    \draw (-0.707,-0.707) circle (0.05cm);

    \draw [-latex] (-1.414,1.414) -- (-0.75,0.75);
    \draw [-latex] (0.75,0.75) -- (1.414,1.414);
    \draw [-latex] (0.75,-0.75) -- (1.414,-1.414);
    \draw [-latex] (-1.414,-1.414) -- (-0.74,-0.74);

    \draw [-latex] (-0.65,0.707) -- (0.65,0.707);
    \draw [-latex] (-0.65,-0.707) -- (0.65,-0.707);
    \draw [-latex] (-0.707,-0.65) -- (-0.707,0.65);
    \draw [-latex] (0.707,0.65) -- (0.707,-0.65);

    \draw (-1,0) node {$a$};
    \draw (0,-1) node {$b$};
    \draw (1,0) node {$c$};
    \draw (0,1) node {$d$};

    \draw (-0.9,1.3) node {$e_1$};
    \draw (-0.9,-1.3) node {$e_2$};
    \draw (0.9,-1.3) node {$e_3$};
    \draw (0.9,1.3) node {$e_4$};

    % Arrow

    \draw [-latex] (2.5,0) -- (3.5,0);

    % After
    \draw [dashed]    (6+0,0)           circle (2cm);
    \draw[fill=black] (6+0.707,0.707)   circle (0.05cm);
    \draw[fill=white] (6+0.707,-0.707)  circle (0.05cm);
    \draw[fill=white] (6+-0.707,0.707)  circle (0.05cm);
    \draw[fill=black] (6+-0.707,-0.707) circle (0.05cm);

    \draw [-latex] (6+-1.414,1.414)  -- (6+-0.75,0.75);
    \draw [-latex] (6+0.75,0.75)     -- (6+1.414,1.414);
    \draw [-latex] (6+0.75,-0.75)    -- (6+1.414,-1.414);
    \draw [-latex] (6+-1.414,-1.414) -- (6+-0.74,-0.74);

    \draw [-latex] (6+-0.65,0.707)  -- (6+0.65,0.707);
    \draw [-latex] (6+-0.65,-0.707) -- (6+0.65,-0.707);
    \draw [-latex] (6+-0.707,0.65) -- (6+-0.707,-0.65);
    \draw [-latex] (6+0.707,-0.65)   -- (6+0.707,0.65);

    \draw (5-0.25,0) node {$d c \Delta^{-1}$};
    \draw (6+0,-1) node {$\Delta$};
    \draw (7+0.25,0)  node {$\Delta^{-1} a d$};
    \draw (6+0.08,1.05)  node {$d c \Delta^{-1} b c^{-1}$};

    \draw (6+-0.9,1.45) node {$e_1$};
    \draw (6+-0.9,-1.3) node {$e_2$};
    \draw (6+0.9,-1.3) node {$e_3$};
    \draw (6+0.9,1.45) node {$e_4$};

    % Type (II)

    \draw (-4,-5) node {( II )};

    % Before
    \draw [dashed]     (0,-5)      circle (2cm);
    \draw [fill=white] (-0.667,-5) circle (0.05cm);
    \draw [fill=white] (0.667,-5)  circle (0.05cm);

    \draw [-latex] (-2,-5) -- (-0.667-0.05,-5);
    \draw [-latex] (-0.667+0.05,-5) -- (0.667-0.05,-5);
    \draw [-latex] (0.667+0.05,-5) -- (2,-5);
    \draw [-latex] (-0.667,-5+0.05) -- (-0.667,-5+1.886);
    \draw [-latex] (0.667,-5-0.05) -- (0.667,-5-1.886);

    \draw (-1.33,-5.2) node {$a$};
    \draw (0,-5.2)     node {$b$};
    \draw (1.33,-5.2)  node {$c$};
    \draw (-0.85,-4)   node {$x$};
    \draw (0.85,-6)    node {$y$};

    % Arrow
    \draw [-latex] (2.5,-5) -- (3.5,-5);

    % After
    \draw [dashed]     (6+0,-5)      circle (2cm);
    \draw [fill=white] (6+-0.667,-5) circle (0.05cm);
    \draw [fill=white] (6+0.667,-5)  circle (0.05cm);

    \draw [-latex] (6+-2,-5)          -- (6+-0.667-0.05,-5);
    \draw [-latex] (6+-0.667+0.05,-5) -- (6+0.667-0.05,-5);
    \draw [-latex] (6+0.667+0.05,-5)  -- (6+2,-5);
    \draw [-latex] (6+0.667,-5+0.05) -- (6+-0.667,-5+1.886);
    \draw [-latex] (6-0.667,-5-0.05)  -- (6+0.667,-5-1.886);

    \draw (6+-1.33,-5.2) node {$a$};
    \draw (6+0,-5.2)     node {$b$};
    \draw (6+1.33,-5.2)  node {$c$};
    \draw (6+0.5,-5+1)   node {$b^{-1}x$};
    \draw (6-0.3,-5-1)   node {$by$};

    % Type (III)

    \draw (-4,-10) node {( III )};

    % Before
    \draw [dashed]     (0,-10)      circle (2cm);
    \draw [fill=black] (-0.667,-10) circle (0.05cm);
    \draw [fill=black] (0.667,-10)  circle (0.05cm);

    \draw [-latex] (-2,-10) -- (-0.667-0.05,-10);
    \draw [-latex] (-0.667+0.05,-10) -- (0.667-0.05,-10);
    \draw [-latex] (0.667+0.05,-10) -- (2,-10);
    \draw [-latex] (-0.667,-10+1.886) -- (-0.667,-10+0.05);
    \draw [-latex] (0.667,-10-1.886) -- (0.667,-10-0.05);

    \draw (-1.33,-10.2) node {$a$};
    \draw (0,-10.2)     node {$b$};
    \draw (1.33,-10.2)  node {$c$};
    \draw (-0.85,-9)    node {$x$};
    \draw (0.85,-11)    node {$y$};

    % Arrow
    \draw [-latex] (2.5,-10) -- (3.5,-10);

    % After
    \draw [dashed]     (6+0,-10)      circle (2cm);
    \draw [fill=black] (6+-0.667,-10) circle (0.05cm);
    \draw [fill=black] (6+0.667,-10)  circle (0.05cm);

    \draw [-latex] (6+-2,-10)          -- (6+-0.667-0.05,-10);
    \draw [-latex] (6+-0.667+0.05,-10) -- (6+0.667-0.05,-10);
    \draw [-latex] (6+0.667+0.05,-10)  -- (6+2,-10);
    \draw [-latex] (6+-0.667,-10+1.886) -- (6+0.667,-10+0.05);
    \draw [-latex] (6+0.667,-10-1.886) -- (6-0.667,-10-0.05);

    \draw (6+-1.33,-10.2) node {$a$};
    \draw (6+0,-10.2)     node {$b$};
    \draw (6+1.33,-10.2)  node {$c$};
    \draw (6+0.3,-10+1)   node {$xb$};
    \draw (6-0.5,-10-1)   node {$yb^{-1}$};

\end {tikzpicture}
\end {figure}
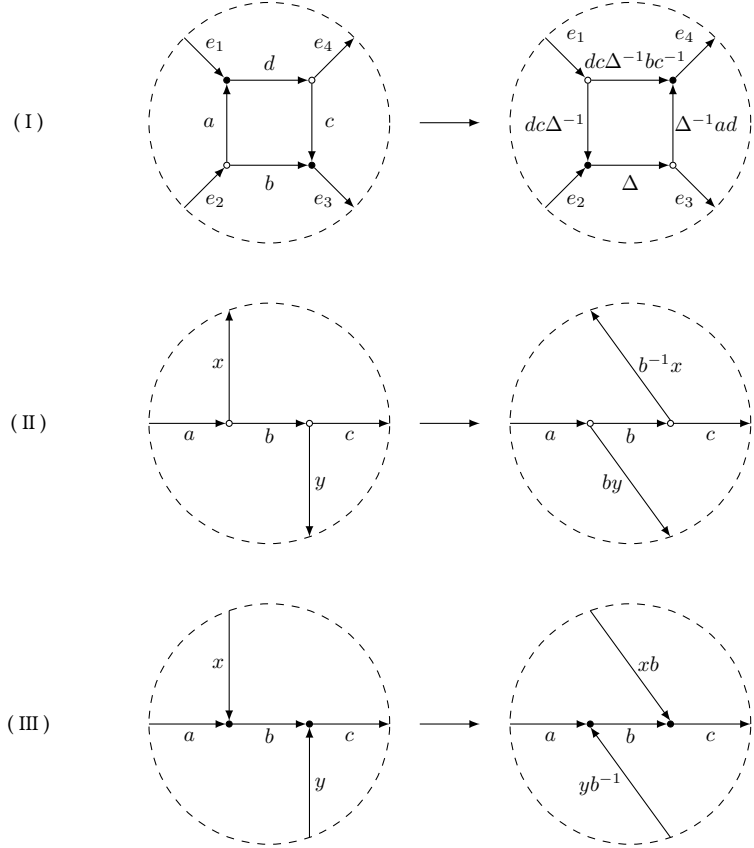

We can also define non-commutative versions of the the Postnikov moves, which preserve the boundary measurements.
They are pictured in \textbf{Figure \ref{fig:noncomm_postnikov_moves}}. As in the commutative case, $\Delta := b + adc$.
We must take special care at this point to say what we mean by expressions such as $(b+acd)^{-1}$. We do so now.

\bigskip

\subsection {The Free Skew Field and Mal'cev Neumann Series}

\bigskip

In order to make sense of expressions of the form $(x+y)^{-1}$, we need an appropriate
notion of ``noncommutative rational functions''. This will be the \emph{free skew field},
which we will define below. Then, we will discuss how the free skew field can be identified
with a certain subset of noncommutative formal power series.

\medskip

For a set $X = (x_1,\dots,x_n)$ of formal noncommuting variables, the \emph{free skew field} or \emph{universal field of fractions},
which we will denote by $\mathcal{F}_{\Bbb{Q}}(X)$, or just $\mathcal{F}(X)$, is a division algebra over $\Bbb{Q}$
characterized by the following universal property: 

\medskip

There is an injective homomorphism $i \colon \Bbb{Q}\left<X\right> \hookrightarrow \mathcal{F}_{\Bbb{Q}}(X)$
from the free associative algebra on $X$ into $\mathcal{F}_{\Bbb{Q}}(X)$ such that for any homomorphism $\varphi \colon \Bbb{Q}\left<X\right> \to D$
into a division ring $D$, there is a unique subring $\Bbb{Q}\left< X \right> \subset \mathcal{R}_\varphi \subset \mathcal{F}_{\Bbb{Q}}(X)$ and
a homomorphism $\psi \colon \mathcal{R} \to D$ with $\varphi = \psi \circ i$, and such that if $0 \neq a \in \mathcal{R}$ and $\psi(a) \neq 0$, then $a^{-1} \in \mathcal{R}$.

\medskip

The explicit construction of $\mathcal{F}(X)$ is a bit complicated (see \cite{grw02} or \cite{cohn77} for details), 
but informally it consists of noncommutative rational expressions
in the variables $x_1,\dots,x_n$, under some suitable notion of equivalence. For example, $\mathcal{F}(X)$ contains expressions such as
$(1+x)^{-1}$ and $w(x+y)^{-1}z$.

\medskip

Another way to embed the free algebra into a division ring is by so-called \emph{Mal'cev Neumann series} (also called Hahn-Mal'cev-Neumann series).
Suppose that the free group generated by $X$ is given some order relation compatible with multiplication. This means that if $f \leq g$,
then $hf \leq hg$ for any $h$. Given an ordering of the free group, the ring of Mal'cev Nemuann series is the subset of all formal
series over the free group which have well-ordered support. This is a division ring which contains $\Bbb{Q}\left<X\right>$.
Jaques Lewin proved \cite{lewin74} that the free skew field $\mathcal{F}(X)$ is isomorphic to the subfield of the ring of
Mal'cev Neumann series generated by the variables $x_1,\dots,x_n$. In the next section, we show an example of expanding
an element of the free skew field as a series.

\medskip

Since $A = A_Q = \Bbb{Q}\left<Q_1\right>$ is a free algebra, it is a subalgebra of $\mathcal{F}(Q_1)$.
So, when we consider the noncommutative Postnikov ``square move'', we interpret the expressions $(b+adc)^{-1}$
as elements of the free skew field $\mathcal{F}(Q_1)$.
If we choose an order relation
compatible with multiplication, then we may also view these expressions as noncommutative Mal'cev Neumann series.
Thus, from now on, given a quiver, we will work more generally with the skew field $\mathcal{F}(Q_1)$ rather than
the free algebra $\Bbb{Q}\left<Q_1\right>$.

\subsection {An Example of a Series Expansion}

Consider the free skew field on two variables $x$ and $y$. We will write the element $(x+y)^{-1}$ as a noncommutative power series.
As mentioned above, we need to put an order on the free group generated by $x$ and $y$ to determine the ring of Mal'cev Neumann series.
We will use the order induced by the \emph{Magnus embedding} of the free group into the ring $\Bbb{Z}[[x,y]]$ of formal power series in 
non-commuting variables $x$ and $y$ \cite{magnus_76}.

\medskip 

More specifically, we embed the free group into $\Bbb{Z}[[x,y]]$ by the map $\alpha \mapsto 1+\alpha$ and $\alpha^{-1} \mapsto \sum_{i \geq 0} (-1)^i \alpha^i$
where $\alpha$ is either $x$ or $y$. This embeds the free group as a subgroup of the multiplicative group of all power series with constant term $1$.
Choosing an order of the variables, say $x < y$, this determines a graded lexicographic order on monomials/words in $\Bbb{Z}[[x,y]]$. Then
we define an order on $\Bbb{Z}[[x,y]]$ where $f < g$ if the coefficient of $f$ at the first place they disagree is smaller. This induces an order
on the free group. We then define the ring of Mal'cev Neumann series to be those series which have well-ordered support with respect to
this ordering of the free group.

\medskip

There are a couple different ways we could try to expand $(x+y)^{-1}$ as a series. We could factor out either of the variables, and then expand as a geometric series.
For example, we could first factor out $x$ to get
\[ (x+y)^{-1} = (x(1+x^{-1}y))^{-1} = (1+x^{-1}y)^{-1}x^{-1} \]
Then we expand $(1+x^{-1}y)^{-1}$ as a geometric series to get
\[ (x+y)^{-1} = (1-x^{-1}y + x^{-1}yx^{-1}y + \cdots)x^{-1} = x^{-1} - x^{-1}yx^{-1} + x^{-1}yx^{-1}yx^{-1} + \cdots \]
To see if the support is well-ordered, we need to go through the above construction. Under the Magnus embedding, we have
\begin {align*}
    x^{-1} &\mapsto 1 - x + x^2 - x^3 + \cdots \\
    x^{-1}yx^{-1} &\mapsto 1 - 2x + y + 3x^2 - xy - yx + \cdots \\
    x^{-1}yx^{-1}yx^{-1} &\mapsto 1 - 3x + 2y + 6x^2 + y^2 - 3xy - 3yx + \cdots \\
    &\vdots \\
    (x^{-1}y)^n x^{-1} &\mapsto 1 - (n+1)x + n y + \cdots
\end {align*}
The sequence of terms in our series expansion for $(x+y)^{-1}$ all differ at the coefficient for $x$. We see that the coefficients
are decreasing in the sequence $-1,-2,-3,\cdots,-(n+1),\cdots$. There is thus no lowest term in this sequence, and so it is \emph{not}
well-ordered. This expansion is then \emph{not} a Mal'cev Neumann series for our chosen ordering.

\medskip

If, however, we factor out $y$ instead of $x$, things will work out. Then we get
\[ (x+y)^{-1} = y^{-1} (1+xy^{-1})^{-1} = y^{-1} - y^{-1}xy^{-1} + y^{-1}xy^{-1}xy^{-1} + \cdots \]
Then the general term, under the Magnus embedding, will be
\[ y^{-1}(xy^{-1})^n \mapsto 1 + nx - (n+1)y + \cdots \]
Now the sequence of $x$ coefficients is the increasing sequence $0,1,2,3,\cdots,n,\cdots$. Therefore the first term $y^{-1}$ is the lowest in
the sequence, and so the support is well-ordered. This expansion is thus a Mal'cev Neumann series.

\bigskip

\subsection {Double Brackets Associated to a Quiver}

\bigskip

We now define a family of double brackets on the skew field $\mathcal{F}_Q := \mathcal{F}(Q_1)$ which generalize the Poisson structures
on the space of edge weights $\mathcal{E}_Q$ described in section 2. 
We again define it locally, and then describe the concatenation procedure
by which they can be glued together. We start with the two local pictures of white and black vertices:

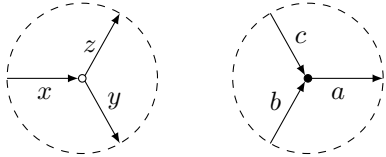
\begin {center}
\begin {tikzpicture}
    \draw [dashed]     (0,0) circle (1cm);
    \draw (0,0) circle (0.05cm);
    \draw [-latex]     (-1,0) -- (-0.05,0);
    \draw [-latex]     (0.035,0.035) -- (0.5,0.866);
    \draw [-latex]     (0.035,-0.035) -- (0.5,-0.866);

    \draw (-0.5,-0.2) node {$x$};
    \draw (0.1,0.43) node {$z$};
    \draw (0.43,-0.3) node {$y$};

    \draw [dashed]     (3,0) circle (1cm);
    \draw [fill=black] (3,0) circle (0.05cm);
    \draw [-latex]     (3.05,0) -- (4,0);
    \draw [-latex]     (2.5,0.866) -- (2.965,0.035);
    \draw [-latex]     (2.5,-0.866) -- (2.965,-0.035);

    \draw (3.4,-0.2) node {$a$};
    \draw (2.57,-0.3) node {$b$};
    \draw (2.9,0.55) node {$c$};
\end {tikzpicture}
\end {center}

Let $\mathcal{F}_\circ = \mathcal{F}(x,y,z)$ be the free skew field (over $\Bbb{Q}$) on 3 generators. We will think of the variables $x,y,z$
as representing the edge weights on the white vertex picture above. Similarly, let $\mathcal{F}_\bullet = \mathcal{F}(a,b,c)$
correspond to the black vertex. Of course they are isomorphic as associative $\Bbb{Q}$-algebras, but we will define different brackets on them.
Choose any $w_1,w_2,w_3 \in \Bbb{Q}$, and define a double bracket $\db{-,-}_\circ$ on $\mathcal{F}_\circ$ by

\begin {align*}
    \db{x,y}_\circ &= w_1 (1 \otimes xy) \\
    \db{x,z}_\circ &= w_2 (1 \otimes xz)  \\
    \db{y,z}_\circ &= w_3 (y \otimes z)  
\end {align*}

Similarly, for scalars $k_1,k_2,k_3 \in \Bbb{Q}$, define $\db{-,-}_\bullet$ on $\mathcal{F}_\bullet$ by

\begin {align*}
    \db{a,b}_\bullet &= k_1 (ba \otimes 1) \\
    \db{a,c}_\bullet &= k_2 (ca \otimes 1) \\
    \db{b,c}_\bullet &= k_3 (c \otimes b)
\end {align*}

Also define $\mathcal{F}_\partial = \mathcal{F}(x) = \Bbb{Q}(x)$ to be the field of rational functions in one variable,
corresponding to a univalent boundary vertex, with trivial double bracket.

\medskip

These double brackets are again ``compatible'' with concatenating/glueing in a similar sense as in the commutative case. 
Let $H_Q$ be the free skew field generated by the half-edges of $Q$, denoted again by $\alpha_s$ and $\alpha_t$ for the
source and target ends of $\alpha \in Q_1$. Then $H_Q$ is the free product of the algebras defined above:
\[ H_Q \cong \left( \freeprod_{\circ \in V_\circ} \mathcal{F}_\circ \right) \, \ast \, 
             \left( \freeprod_{\bullet \in V_\bullet} \mathcal{F}_\bullet \right) \, \ast \,
             \left( \freeprod_{\star \in V_\partial} \mathcal{F}_\partial \right) \]
As is mentioned in \cite{bergh_08}, for algebras $A$ and $B$ with double brackets, there is a uniquely defined double bracket
on the free product $A \ast B$ such that $\db{a,b} = 0$ for $a \in A$ and $b \in B$. Since the free product is the coproduct
in the category of associative algebras (analogously the tensor product is the coproduct for commutative algebras), this is
analogous to the product bracket on $\mathscr{H}_Q$ in the commutative case.

\medskip

Again we define a glueing map $g^* \colon \mathcal{F}_Q \to H_Q$ by $g^*(\alpha) = \alpha_s \alpha_t$ for $\alpha \in Q_1$. Since the
algebras are non-commutative, this is not actually the pull-back of a geometric/topological map, but we keep the notation
for the sake of analogy. We have the following result, which mimics the commutative case:

\bigskip

\begin {prop}
    There is a unique double bracket on $\mathcal{F}_Q$ so that the glueing homomorphism $g^* \colon \mathcal{F}_Q \to H_Q$ satisfies
    \[ \db{g^*(\alpha),g^*(\beta)}_{H_Q} = (g^* \otimes g^*) \left( \db{\alpha,\beta}_{\mathcal{F}_Q} \right) \]
\end {prop}
\begin {proof}
    As in the commutative case, there are many cases to consider, but all of them are similar
    and are simple calculations. We show one as an example to illustrate the idea. We again look
    at \textbf{Figure \ref{fig:concat_local}}, and try to define $\db{w,y}_{\mathcal{F}_Q}$. If it is to
    satisfy the desired property, we must have
    \begin {align*}
        (g^* \otimes g^*) (\db{w,y}_{\mathcal{F}_Q}) &= \db{g^*(w), g^*(y)}_{H_Q} \\
            &= \db{ax, y_s y_t}_{H_Q} \\
            &= (1 \otimes a) \db{x,y_s}_{H_Q} (1 \otimes y_t) \\
            &= w_1 (1 \otimes a) (1 \otimes xy_s) (1 \otimes y_t) \\
            &= w_1 (1 \otimes axy_s y_t) \\
            &= w_1 (1 \otimes g^*(w)g^*(y)) \\
            &= (g^* \otimes g^*)(w_1 (1 \otimes wy))
    \end {align*}
    This suggests that $\db{w,y}_{\mathcal{F}_Q}$ must be defined to be $w_1 (1 \otimes wy)$. As mentioned
    above, all other cases are similar. Just as in the commutative case, the double bracket on $\mathcal{F}_Q$
    is given by the same expressions as in $H_Q$, treating edges as the corresponding half-edges
    which meet at a common vertex.
\end {proof}

\bigskip

We may also consider the doubled quiver $\overline{Q}$, as before. 
We then associate to each opposite arrow $\alpha^*$ the element $\alpha^{-1}$ in the free skew field.
The double bracket defined above on generators extends in a unique way to the free skew field by the formulas
\begin {align*}
    \db{\beta,\alpha^{-1}} &= -(\alpha^{-1} \otimes 1) \db{\beta,\alpha} (1 \otimes \alpha^{-1}) \\
    \db{\alpha^{-1},\beta} &= -(1 \otimes \alpha^{-1}) \db{\alpha,\beta} (\alpha^{-1} \otimes 1)
\end {align*}
We may then interpret any path in $\overline{Q}$ as a non-commutative Laurent monomial in $\mathcal{F}_Q$. 

\bigskip

\subsection {A Formula for the Bracket}

We will primarily be concerned with paths that are closed loops in $\overline{Q}$. Let $\mathscr{L} \subset \mathcal{F}_Q$
be the vector subspace spanned by all monomials which represent closed loops,
and let $f,g \in \mathscr{L}$ be two elements.

\bigskip

As before, we let $f \cap g$ denote the set of all maximal common subpaths.
We will prove that, as in the commutative case, we only need to consider the ends of 
common subpaths in order to compute $\left<f,g\right>$. 
First we give a result about the local structure of the double bracket.

\bigskip

\begin {lem} \label{lem:local_bracket}
    Let $x$ and $y$ be edges in $\overline{Q}$. That is, $x$ and $y$ can either be arrows in the quiver,
    or ``reverse'' arrows. Then \\[0.5ex] \\[0.1ex]
    \begin {tabular}{cp{5in}}
        ($a$) & If $y$ follows $x$ (i.e. $s(y) = t(x)$), then $\db{x,y} = \lambda (1 \otimes xy)$ for some $\lambda \in \Bbb{Q}$. \\[1.5ex]
        ($b$) & If $x$ and $y$ have the same source, then $\db{x,y} = \lambda (x \otimes y)$ for some $\lambda \in \Bbb{Q}$. \\[1.5ex]
        ($c$) & If $x$ and $y$ have the same target, then $\db{x,y} = \lambda (y \otimes x)$ for some $\lambda \in \Bbb{Q}$.
    \end {tabular}
\end {lem}
\begin {proof}
    The pictures below show the possibilities at a white vertex.
    The red are for part $(a)$, the blue for part $(b)$, and green for part $(c)$.
    Note that there are 3 additional possibilities for part $(a)$, given by the reversal/inverse of the pairs shown
    in red.
    
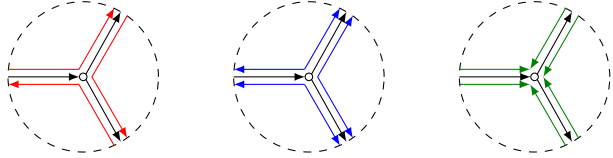
\begin {center}
    \begin {tikzpicture}
        % Left : consecutive
        \draw [dashed]     (0,0) circle (1cm);
        \draw (0,0) circle (0.05cm);

        \draw [-latex]     (-1,0) -- (-0.05,0);
        \draw [-latex]     (0.035,0.035) -- (0.5,0.866);
        \draw [-latex]     (0.035,-0.035) -- (0.5,-0.866);

        \draw [-latex, red] (-1,0.1) -- (-0.0577,0.1) -- (0.415,0.918);
        \draw [-latex, red] (0.413,-0.916) -- (-0.0578,-0.1) -- (-1,-0.1);
        \draw [-latex, red] (0.587,0.816) -- (0.116,0) -- (0.587,-0.816);

        % Middle : same source
        \draw [dashed]     (3,0) circle (1cm);
        \draw (3,0) circle (0.05cm);

        \draw [-latex]     (3+-1,0)         -- (3+-0.05,0);
        \draw [-latex]     (3+0.035,0.035)  -- (3+0.5,0.866);
        \draw [-latex]     (3+0.035,-0.035) -- (3+0.5,-0.866);

        \draw [-latex, blue] (3+-0.0577,0.1)  -- (3+0.415,0.918);
        \draw [-latex, blue] (3+-0.0577,0.1)  -- (3+-1,0.1);
        \draw [-latex, blue] (3+-0.0578,-0.1) -- (3+0.413,-0.916);
        \draw [-latex, blue] (3+-0.0578,-0.1) -- (3+-1,-0.1);
        \draw [-latex, blue] (3+0.116,0)      -- (3+0.587,-0.816);
        \draw [-latex, blue] (3+0.116,0)      -- (3+0.587,0.816);

        % Right : same target
        \draw [dashed]     (6,0) circle (1cm);
        \draw (6,0) circle (0.05cm);

        \draw [-latex]     (6+-1,0)         -- (6+-0.05,0);
        \draw [-latex]     (6+0.035,0.035)  -- (6+0.5,0.866);
        \draw [-latex]     (6+0.035,-0.035) -- (6+0.5,-0.866);

        \draw [-latex, black!50!green] (6+0.415,0.918) -- (6+-0.0577,0.1);
        \draw [-latex, black!50!green] (6+-1,0.1)      -- (6+-0.0577,0.1);
        \draw [-latex, black!50!green] (6+0.413,-0.916)-- (6+-0.0578,-0.1); 
        \draw [-latex, black!50!green] (6+-1,-0.1)     -- (6+-0.0578,-0.1); 
        \draw [-latex, black!50!green] (6+0.587,-0.816)-- (6+0.116,0);      
        \draw [-latex, black!50!green] (6+0.587,0.816) -- (6+0.116,0);      
    \end {tikzpicture}
    \end {center}
    As in the picture from the previous section, call the unique incoming arrow $x$, and the outgoing arrows $y$ and $z$,
    in counter-clockwise order from $x$. Then the three pairs of edges in the picture above for $(a)$ give
    \begin {align*}
        \db{x,z} &= w_2 (1 \otimes xz) \\
        \db{z^{-1},y} &= w_3 (1 \otimes z^{-1}) (z \otimes y) (z^{-1} \otimes 1) = w_3 (1 \otimes z^{-1} y) \\
        \db{y^{-1},x^{-1}} &= -w_1 (x^{-1} \otimes y^{-1}) (xy \otimes 1) (y^{-1} \otimes x^{-1}) = -w_1 (1 \otimes y^{-1}x^{-1})
    \end {align*}
    For part $(b)$, the pictured pairs give
    \begin {align*}
        \db{y,z} &= w_3 (y \otimes z) \\
        \db{x^{-1},y} &= -w_1 (1 \otimes x^{-1})(1 \otimes xy)(x^{-1} \otimes 1) = -w_1 (x^{-1} \otimes y) \\
        \db{x^{-1},z} &= -w_2 (1 \otimes x^{-1})(1 \otimes xz)(x^{-1} \otimes 1) = -w_2 (x^{-1} \otimes z)
    \end {align*}
    For part $(c)$, the pairs pictured above give
    \begin {align*}
        \db{x,y^{-1}}      &= -w_1 (y^{-1} \otimes 1)(1 \otimes xy)(1 \otimes y^{-1})         = -w_1 (y^{-1} \otimes x) \\
        \db{x,z^{-1}}      &= -w_2 (z^{-1} \otimes 1)(1 \otimes xz)(1 \otimes z^{-1})         = -w_2 (z^{-1} \otimes x) \\
        \db{y^{-1},z^{-1}} &= w_3 (z^{-1} \otimes y^{-1})(y \otimes z)(y^{-1} \otimes z^{-1}) = w_3 (z^{-1} \otimes y^{-1})
    \end {align*}
    All other possibilities at a black vertex are verified by similar calculations.
\end {proof}

\bigskip

\begin {lem}
    Let $f,g \in \mathscr{L}$. The induced bracket $\left< f,g \right>$ depends only on the endpoints of maximal common
    subpaths of $f$ and $g$.
\end {lem}
\begin {proof}
    The result can be stated more technically as follows. Suppose that $p$ is an edge in $f$ and $q$ is an edge in $g$.
    Expanding $\db{f,g}$ with the Leibniz rule, there will be a term involving $\db{p,q}$. The result says that unless
    $p$ and $q$ are incident to a common vertex which is an endpoint of a maximal common subpath, then this term is either
    zero, or it cancels with another term in the Leibniz expansion.

    \medskip

    First of all, if $p$ and $q$ are not incident to the same vertex, then $\db{p,q} = 0$. Since each vertex is trivalent,
    if two paths go through a common vertex, then they must have at least one edge in common at that vertex. So every
    non-zero contribution $\db{p,q}$ comes from when $p$ and $q$ belong to a common subpath of $f$ and $g$. We now argue
    that even in this case, $\db{p,q}$ is either zero or cancels unless $p$ and $q$ occur at the ends of a common subpath.
    For the remainder of the proof, let $w$ be a maximal common subpath of $f$ and $g$.

    \medskip

    We first consider the simplest case, which is when $p$ and $q$ are consecutive edges in $w$.
    So suppose that $q$ immediately follows $p$ in $w$. Then we can write $w =  w_1 pq w_2$, and $f$ and $g$
    can be written as $f = f_0 w_1 pq w_2$ and $g = g_0 w_1 pq w_2$. 

    \medskip

    Then after expanding using the Leibniz rules, we get two terms in $\db{f,g}$ involving $\db{p,q}$ and $\db{q,p}$.
    By the previous lemma, $\db{p,q} = \lambda (1 \otimes pq)$ for some $\lambda$, and $\db{q,p} = -\lambda(pq \otimes 1)$.
    Then these two terms are given by
    \[ \lambda (g_0 w_1 p \otimes f_0 w_1) (1 \otimes pq) (q w_2 \otimes w_2) = \lambda g_0 w_1 pq w_2 \otimes f_0 w_1 pq w_2 \]
    \[ \mathrm{and} \]
    \[ -\lambda (g_0 w_1 \otimes f_0 w_1 p) (pq \otimes 1) (w_2 \otimes q w_2) = -\lambda g_0 w_1 pq w_2 \otimes f_0 w_1 pq w_2 \]
    Obviously, these terms cancel.

    \medskip

    There is also the case that the path $w$ intersects itself, so it may be possible that $p$ and $q$ share a common vertex,
    but do not occur consecutively in the path. There are three cases,
    corresponding to parts $(a)$, $(b)$, and $(c)$ in the previous lemma.

    \medskip

    First we consider part $(a)$ of the previous lemma. So suppose $p$ and $q$ are consecutive arrows in the quiver, but
    they occur non-consecutively in the path. Let $r$ be the third edge incident to the common vertex of $p$ and $q$,
    oriented outward. There are three cases:
    \begin {enumerate}
        \item $p$ and $q$ occur consecutively twice, so $w = w_1 pq w_2 pq w_3$.
        \item The path follows $r$ after $p$, and $p$ before $q$, so $w = w_1 pr w_2 pq w_3$.
        \item The path follows $r$ after $p$, and $r^{-1}$ before $q$, so $w = w_1 pr w_2 r^{-1} q w_3$.
    \end {enumerate}
    In cases $(2)$ and $(3)$, the vertex in question is the end of another common subpath, and so we don't
    consider these cases.
    In case $(1)$, the brackets of the first $p$ and first $q$ cancel by the previous discussion, since they are
    consecutive edges in the path. The terms coming
    from the bracket of the first $p$ with the second $q$ and vice versa are given by
    \[ \lambda (g_0 w_1 pq w_2 pq w_2 pq w_3 \otimes f_0 w_1 pq w_3) ~ \mathrm{and} ~ -\lambda (g_0 w_1 pq w_2 pq w_2 pq w_3 \otimes f_0 w_1 pq w_3) \]

    \medskip

    All other cases correspond to cases $(b)$ and $(c)$ from the previous lemma. All these cases imply that the vertex
    which $p$ and $q$ share is the end of some other common subpath, and so we don't need to consider these. We have
    thus covered all cases, and the lemma is proved.
\end {proof}

\bigskip

Write $f$ and $g$ as products of edges: $f = f_1 \cdots f_n$ and $g = g_1 \cdots g_m$. Then by the Leibniz rule:
\[ \db{f,g} = \sum_{i,j} (g_1 \cdots g_{j-1} \otimes f_1 \cdots f_{i-1}) \db{f_i,g_j} (f_{i+1} \cdots f_n \otimes g_{j+1} \cdots g_m) \]
The result of the preceding lemma says that this sum is only over the pairs $f_i,g_j$ which are incident to a common
vertex which is the end of a maximal common subpath of $f$ and $g$. So we will now compute what happens at these
vertices.

\medskip

As in the commutative case, there appear to be 36 cases to consider, since the two endpoints of a common subpath could
each have one of two colors and one of three orientations. It happens again, however, that the orientations do not affect
the outcome, which we formulate now as a lemma. 

\bigskip

\begin {lem}
    The contribution to $\left<f,g\right>$ coming from a common subpath $w \in f \cap g$ depends only on the colors
    of the endpoint vertices, and not on the orientations.
\end {lem}
\begin {proof}
    Let us first consider the beginning of a common subpath $w$, where the paths $f$ and $g$ first come together.
    Then we may write $f = f_1 \cdots f_i w_1 \cdots f_n$ and $g = g_1 \cdots g_j w_1 \cdots g_m$, where $w_1 = f_{i+1} = g_{j+1}$
    is the first edge in the common subpath $w$, and $f_i$ and $g_j$ are incident to the same vertex as $w_1$.
    There are three terms in the expansion of $\db{f,g}$ coming from this vertex, corresponding to $w$, given by
    \[ (g_1 \cdots g_{j-1} \otimes f_1 \cdots f_{i-1}) \db{f_i,g_j} (w_1f_{i+2} \cdots f_n \otimes w_1 g_{j+2} \cdots g_m) \]
    \[ (g_1 \cdots g_{j} \otimes f_1 \cdots f_{i-1}) \db{f_i,w_1} (w_1f_{i+2} \cdots f_n \otimes g_{j+2} \cdots g_m) \]
    \[ (g_1 \cdots g_{j-1} \otimes f_1 \cdots f_{i}) \db{w_1,g_j} (f_{i+2} \cdots f_n \otimes w_1 g_{j+2} \cdots g_m) \]
    By \textbf{Lemma \ref{lem:local_bracket}}, all three of these terms give the simple tensor
    \[ g_1 \cdots g_j w_1 f_{i+2} \cdots f_n \otimes f_1 \cdots f_i w_1 g_{j+2} \cdots g_m, \]
    but with different coefficients. Which coefficient goes with which term depends on the orientation of the vertex.
    The following picture illustrates the three possibilities at a black vertex, with $f_i w_1$ in blue and $g_j w_1$ in red.

    
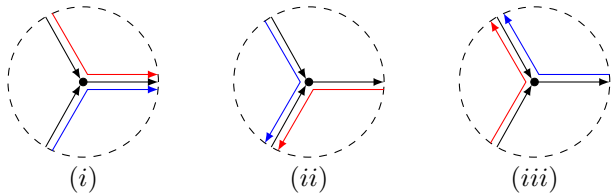
\begin {center}
    \begin {tikzpicture}
        % Left
        \draw [dashed]     (0,0) circle (1cm);
        \draw [fill=black] (0,0) circle (0.05cm);

        \draw [-latex]     (0.05,0)      -- (1,0);
        \draw [-latex]     (-0.5,0.866)  -- (-0.035,0.035);
        \draw [-latex]     (-0.5,-0.866) -- (-0.035,-0.035);

        \draw [-latex, red]  (-0.4134,0.916) -- (0.0577,0.1) -- (1,0.1);
        \draw [-latex, blue] (-0.4134,-0.916) -- (0.0577,-0.1) -- (1,-0.1);

        \draw (0,-1.3) node {$(i)$};

        % Middle
        \draw [dashed]     (3,0) circle (1cm);
        \draw [fill=black] (3,0) circle (0.05cm);

        \draw [-latex, black]     (3+0.05,0)      -- (3+1,0);
        \draw [-latex, black]     (3+-0.5,0.866)  -- (3+-0.035,0.035);
        \draw [-latex, black]     (3+-0.5,-0.866) -- (3+-0.035,-0.035);

        \draw [-latex, red]  (3+1,-0.1) --  (3+0.0577,-0.1) --  (3-0.4134,-0.916);
        \draw [-latex, blue] (3+-0.587,0.816) -- (3-0.115,0) -- (3-0.587,-0.816);

        \draw (3,-1.3) node {$(ii)$};

        % Right
        \draw [dashed]     (6,0) circle (1cm);
        \draw [fill=black] (6,0) circle (0.05cm);

        \draw [-latex, black]     (6+0.05,0)      -- (6+1,0);
        \draw [-latex, black]     (6+-0.5,0.866)  -- (6+-0.035,0.035);
        \draw [-latex, black]     (6+-0.5,-0.866) -- (6+-0.035,-0.035);

        \draw [-latex, blue]  (6+1,0.1) -- (6+0.0577,0.1) -- (6+-0.4134,0.916);
        \draw [-latex, red]   (6-0.587,-0.816) -- (6-0.115,0) -- (6+-0.587,0.816);

        \draw (6,-1.3) node {$(iii)$};

    \end {tikzpicture}
    \end {center}

    The coefficients for the three pictures above are given in the following table:
    \[
        \begin {array}{|c|c|c|c|} \hline
            \mathrm{Case} & \db{f_i,g_j}           & \db{f_i,w_1}             & \db{w_1,g_j} \\ \hline
            (i)           & k_3  (g_j \otimes f_i) & -k_1 (1 \otimes f_i w_1) & k_2  (g_j w_1 \otimes 1) \\ \hline
            (ii)          & k_2  (g_j \otimes f_i) & k_3  (1 \otimes f_i w_1) & -k_1 (g_j w_1 \otimes 1) \\ \hline
            (iii)         & -k_1 (g_j \otimes f_i) & k_2  (1 \otimes f_i w_1) & k_3  (g_j w_1 \otimes 1 ) \\ \hline
        \end {array}
    \]
    In all three cases, the three terms are the same simple tensor, so the coefficients add up. As we can see from the table,
    the combined coefficient is always $-A_\bullet = k_2 + k_3 - k_1$. Similarly, for all three orientations at a white
    vertex, we will get the same simple tensor with a coefficient of $A_\circ = w_1 - w_2 - w_3$. If we consider the endpoint
    vertex of the common subpath, we will again get the same simple tensor for all three terms. If the vertex is black, we will
    get a coefficient of $A_\bullet$, and if it is white, we will get $-A_\circ$. 
    These are exactly the same as the values for $\varepsilon_w(f,g)$ in the commutative case.
\end {proof}

\bigskip

Let $f_wg_w$ denote the loop which follows $f$ starting with the path $w$, followed by $g$ starting at $w$. 
Then putting the previous lemmas together gives

\bigskip

\begin {thm} \label{thm:bracket_formula}
    Let $f,g \in \mathscr{L}$. Then the induced bracket in $\mathcal{F}_Q^\natural$ is given by
    \[ \left<f,g\right> = \sum_{x \in f \cap g} \varepsilon_x(f,g) f_xg_x \]
    In particular, $\mathscr{L}^\natural$ is closed under $\left<-,-\right>$.
\end {thm}

\bigskip

Moreover, we have the following important observation.

\bigskip

\begin {thm} \label{thm:jacobi_identity}
    The induced bracket $\left<-,-\right>$ makes $\mathscr{L}^\natural$ a Lie algebra.
\end {thm}
\begin {proof}
    We only need to verify the Jacobi identity.
    So let $f,g,h \in \mathscr{L}$. Then
    \[ \left< f, \left<g,h\right> \right> = \sum_{i \in g \cap h} \varepsilon_i(g,h) \left< f, g_ih_i \right>
       = \sum_{i \in g \cap h} \varepsilon_i(g,h) \left( \sum_{j \in f \cap g} \varepsilon_j(f,g) f_j (g_ih_i)_j 
                                                  + \sum_{k \in f \cap h} \varepsilon_k(f,h) f_k (g_ih_i)_k \right) 
    \]

    \[
        \left< \left<f,g\right>, h \right> = \sum_{j \in f \cap g} \varepsilon_j(f,g) \left< f_jg_j, h \right>
        = \sum_{j \in f \cap g} \varepsilon_j(f,g) \left( \sum_{k \in f \cap h} \varepsilon_k(f,h) (f_jg_j)_k h_k
                                                   + \sum_{i \in g \cap h} \varepsilon_i(g,h) (f_jg_j)_i h_i \right)
    \]

    \[
        \left< g, \left<f,h\right> \right> = \sum_{k \in f \cap h} \varepsilon_k(f,h) \left< g, f_k h_k \right>
        = \sum_{k \in f \cap h} \varepsilon_k(f,h) \left( \sum_{i \in g \cap h} \varepsilon_i(g,h) g_i(f_k h_k)_i
                                                   - \sum_{j \in f \cap g} \varepsilon_j(f,g) g_j(f_k h_k)_j \right)
    \]

    The Jacobi identity will hold if the first equation is equal to the sum of the second and third.
    Comparing the right-hand side of the first with the sum of the right-hand sides of the second and third, 
    we see that the Jacobi identity will hold if the following identities are true:
    \[ f_j(g_ih_i)_j = (f_jg_j)_i h_i = g_i(f_jh_j)_i \]
    It is an easy check that these identities are indeed true for all cycles.
\end {proof}

\bigskip

\begin {rmk}
    This double bracket is \emph{not} a double Poisson bracket (it does not satisfy
    the double Jacobi identity), so it does not immediately guarantee the Jacobi identity for $\left<-,-\right>$. 
    Since 
    we do have the Jacobi identity, though, for $\left<-,-\right>$, this makes it into an $H_0$-Poisson structure
    as defined by Crawley-Boevey \cite{crawley-boevey_11}.
\end {rmk}

\bigskip

\subsection {Goldman's Bracket and the Twisted Ribbon Surface}

\bigskip

In this section, we give a geometric interpretation of the bracket just described. We first recall some preliminaries.

\medskip

Let $G$ be a connected Lie group and $S$ a smooth oriented surface with fundamental group $\pi := \pi_1(S)$.
We consider the space of representations of $\pi$ in $G$, modulo conjugations, which we call $\mathrm{Rep}_G(\pi)$:
\[ \mathrm{Rep}_G(\pi) := \mathrm{Hom}(\pi,G)/G \]
Let $f \colon G \to \Bbb{R}$ be any invariant function on $G$, with respect to conjugation. Then for $\alpha \in \pi$,
we can define the function $f_\alpha \colon \mathrm{Rep}_G(\pi) \to \Bbb{R}$ by the formula $f_\alpha([\varphi]) := f(\varphi(\alpha))$.
In particular, if $G$ is a group of matrices, we may take $f = \mathrm{tr}$. In this case we will write $\mathrm{tr}(\alpha)$
for $f_\alpha$. Also, if $G = \mathrm{GL}_n(\Bbb{R})$, we will write $\mathrm{Rep}_n(\pi)$ instead of $\mathrm{Rep}_{\mathrm{GL}_n(\Bbb{R})}(\pi)$.

\medskip

In 1984, William Goldman described a symplectic structure on $\mathrm{Rep}_G(\pi)$ which generalizes the Weil-Petersson
symplectic structure on Teichm\"{u}ller space in the case that $G = \mathrm{PSL}_2(\Bbb{R})$ \cite{goldman_84}. Then, in 1986,
he studied the Poisson bracket induced by this symplectic structure, in terms of the functions $f_\alpha$ \cite{goldman_86}.
Goldman gives explicit formulas for $\{f_\alpha,f_\beta\}$ for various choices of the group $G$,
in terms of the topology of the surface and the intersection of curves representing $\alpha$ and $\beta$. In particular,
when $G = \mathrm{GL}_n(\Bbb{R})$ and $f = \mathrm{tr}$, we have the following.

\begin {thm} \cite{goldman_86}
    The Poisson bracket of the functions $\mathrm{tr}(\alpha)$ on $\mathrm{Rep}_n(\pi)$ is given by
    \[ \{\mathrm{tr}(\alpha), \mathrm{tr}(\beta)\} = \sum_{p \in \alpha \cap \beta} \varepsilon_p(\alpha,\beta) \, \mathrm{tr}(\alpha_p \beta_p) \]
    Here, $\alpha_p \beta_p$ denotes the loop which traverses first $\alpha$, and then $\beta$, both based at the point $p$,
    and $\varepsilon_p(\alpha,\beta)$ is the oriented intersection number of the curves at $p$.
\end {thm}

Note the obvious similarity with our formula from \textbf{Theorem \ref{thm:bracket_formula}}. We will now formulate a geometric interpretation
of our double bracket from a quiver so that a special choice of constants $A_\bullet$ and $A_\circ$ realizes this Goldman bracket. First we define $\hat{\pi}$
to be the set of conjugacy classes in $\pi$ (or free homotopy classes of loops). Goldman observes that the bracket above induces a Lie bracket on the vector space spanned by $\hat{\pi}$.
Formally, we simply remove the ``$\mathrm{tr}$'' in the formula above. So
for $\alpha,\beta \in \hat{\pi}$, we have
\[ [\alpha,\beta] = \sum_{p \in \alpha \cap \beta} \varepsilon_p(\alpha,\beta) \, \alpha_p \beta_p \]
With this setup, we have the following

\begin {thm} \cite{goldman_86}
    $\Bbb{Z}\hat{\pi}$ is a Lie algebra with the bracket shown above, and the map $\mathrm{tr} \colon \Bbb{Z}\hat{\pi} \to C^\infty(\mathrm{Rep}_n(\pi))$
    given by $\alpha \mapsto \mathrm{tr}(\alpha)$ is a Lie algebra homomorphism.
\end {thm}

\medskip

In fact, our proof of \textbf{Theorem \ref{thm:jacobi_identity}} is essentially the same as Goldman's original proof (\cite{goldman_86}, \textbf{Theorem 5.3}).

\medskip

In the literature, fat graphs are also commonly called ``ribbon graphs''. From a fat graph $\Gamma$, one can construct an oriented surface with boundary, $S_\Gamma$, 
by replacing the edges with rectangular strips (ribbons) and the vertices with discs, where the ribbons are glued to the discs according to the cyclic ordering
prescribed by the fat graph structure. We call this surface the \emph{ribbon surface} associated to the fat graph. It is clear that the original graph $\Gamma$
is a deformation retract of $S_\Gamma$. We say that $\Gamma$ is a \emph{spine} of the surface $S = S_\Gamma$.

\medskip

Given a quiver $Q$ (oriented fat graph) with underlying unoriented graph $\Gamma$, we consider, as before, the subspace $\mathscr{L} \subset \mathcal{F}_Q$ of loops. 
Then in a natural way we identify the cyclic space $\mathscr{L}^\natural$ with $\Bbb{Q}\hat{\pi}$,
the space generated by free homotopy classes of loops on $S_\Gamma$. 
In the commutative case, we considered the case $A_\bullet=A_\circ=-\frac{1}{2}$, which gives $\varepsilon_p(f,g) = 1$
when $f$ and $g$ touch and $\varepsilon_p(f,g) = 0$ when $f$ and $g$ cross.
We will again be primarily concerned with this specific choice of coefficients in the non-commutative case.
This is opposite from what $\varepsilon_p(\alpha,\beta)$ means for Goldman's bracket, however. This is because if paths $f$ and $g$ touch in the quiver, then
on $S_\Gamma$, they are homologous to paths which do not touch at all, and so their Goldman bracket should be zero. The way around this is to consider the ``dual'' surface.

\medskip

\begin {defn}
    We define the \emph{dual surface} to $S_\Gamma$, which we denote $\widetilde{S}_\Gamma$, to be glued out of ribbons like $S_\Gamma$, except that whenever
    an edge joins vertices of different colors, the corresponding ribbon is given a half-twist. 
\end {defn}

\bigskip

\begin {thm}
    For the choice of coefficients $A_\bullet = A_\circ = -\frac{1}{2}$, the induced bracket $\left<-,-\right>$ on $\mathscr{L}^\natural$
    coincides with Goldman's bracket under the identification of $\mathscr{L}^\natural$ with $\Bbb{Q} \hat{\pi}$,
    where $\pi = \pi_1(\widetilde{S}_\Gamma)$.
\end {thm}
\begin {proof}
We think of the quiver as a planar projection picture
of the ribbon surface $\widetilde{S}_\Gamma$, where one of the colors is the ``top'' of the ribbon, and the other color is the ``bottom''. 
If we were to ``untwist'' the surface $\widetilde{S}_\Gamma$,
and try to view it without twists, we would see that paths which touch in the quiver end up crossing in $\widetilde{S}_\Gamma$, and paths which cross in
the quiver end up not touching in $\widetilde{S}_\Gamma$. Now that touching and crossing have been interchanged, we see that the Goldman bracket on $\widetilde{S}_\Gamma$
coincides with the induced bracket $\left<-,-\right>$ on $\mathscr{L}^\natural$ when $A_\bullet = A_\circ = -\frac{1}{2}$.
\end {proof}

\bigskip

We will sometimes want to choose a vertex in the quiver $Q$, and only consider loops in the fundamental
group of $S_\Gamma$ based at that point. After choosing a vertex, we let $\mathscr{L}_\bullet$ denote
the subspace of loops which start and end at this point. Then $\mathscr{L}_\bullet$ is naturally
identified with the group algebra of $\pi = \pi_1(S_\Gamma,\bullet)$. The induced bracket is independent of this
choice of basepoint, since conjugate elements are equivalent in $\mathscr{L}_\bullet^\natural$. Choosing a basepoint
gives some nice properties. For example, $\mathscr{L}_\bullet$ is a subalgebra, unlike $\mathscr{L}$, since
the common basepoint allows us to compose paths. Also, we have the following property for the double bracket.

\bigskip

\begin {prop}
    The subalgebra $\mathscr{L}_\bullet$ is closed under the double bracket.
    That is, $\db{\mathscr{L}_\bullet,\mathscr{L}_\bullet} \subseteq \mathscr{L}_\bullet \otimes \mathscr{L}_\bullet$. 
\end {prop}
\begin {proof}
    Let $f,g \in \mathscr{L}_\bullet$. Write $f$ and $g$ as monomials, $f = f_1 \cdots f_k$ and $g = g_1 \cdots g_\ell$,
    where each $f_i$ and $g_j$ are arrows in the quiver. Note that since they are both based at the point $\bullet$,
    we have $s(f_1) = s(g_1) = t(f_k) = t(g_\ell) = \bullet$. Now, using the Leibniz rule:
    \[ \db{f,g} = \sum_{i,j} (g_1 \cdots g_{j-1} \otimes f_1 \cdots f_{i-1}) \db{f_i,g_j} (f_{i+1} \cdots f_k \otimes g_{j+1} \cdots g_\ell) \]
    Using the formulas from \textbf{Lemma \ref{lem:local_bracket}}, and examining the three cases, we see that each
    term in the sum above is of the form $\alpha \otimes \beta$, where $\alpha$ and $\beta$ are both in $\mathscr{L}_\bullet$.
\end {proof}

\bigskip

We now point out some similar and related work which was brought to the author's attention while
working on the paper.

\bigskip

\begin {rmk}
    In \cite{turaev_12}, Turaev and Massuyeau construct a \emph{quasi-Poisson bracket} on the character variety
    $\mathrm{Rep}_n(\pi)$, which is induced by a \emph{double quasi-Poisson bracket} on the group algebra of $\pi$.
    As was mentioned above, we may identify the group algebra of $\pi$ with $\mathscr{L}_\bullet$, and so they
    are both double brackets on $\mathscr{L}_\bullet$.
    However, the double bracket considered in this paper is \emph{not} a double quasi-Poisson bracket, and
    so the two constructions are not exactly the same. However, the bracket of Turaev and Massuyeau also
    projects to the Goldman bracket in $\mathscr{L}_\bullet^\natural$.
    Also, Semeon Artamonov, in his recent thesis \cite{artamonov_18}, constructed a much more abstract and categorical
    version of the quasi-Poisson structure of Turaev and Massuyeau.
\end {rmk}

\bigskip

\subsection {The $X,Y$ Variables}

\bigskip

Just as in the commutative case, we will use gauge transformations to obtain new weights on the graph.
Start with the quiver $Q_n$ just as in \textbf{Figure \ref{fig:pentagram_quiver}}. Give the edge weights the same names
as in the commutative case, but now they are formal noncommutative variables in $\mathcal{F}_Q$. Perform guage transformations to
obtain variables $a_i,b_i,c_i,d_i$ just as before. Note that the pictures in \textbf{Section 3.2} are actually
noncommutative gauge transformations, as they take into account whether multiplication happens on the left or right.
In this case, the double bracket induced on the $a,b,c,d$ variables is given by:

\begin {alignat*}{2}
    \db{b_i,a_i} &= \frac{1}{2} b_i \otimes a_i  &\hspace{1cm} \db{b_i,c_i} &= \frac{1}{2} c_i \otimes b_i \\
    \db{a_i,d_i} &= \frac{1}{2} 1 \otimes a_id_i &\hspace{1cm} \db{c_i,d_i} &= \frac{1}{2} d_ic_i \otimes 1
\end {alignat*}

Define the following monomials in the $a,b,c,d$ variables:

\begin {align*}
    x_i &= a_i c_{i-1}^{-1} d_{i-1}^{-1} c_{i-2}^{-1} \\
    y_i &= b_i c_i^{-1} d_i^{-1} c_{i-1}^{-1} d_{i-1}^{-1} c_{i-2}^{-1} \\
    z_k &= d_1c_1 \cdots d_{k-1}c_{k-1} d_k
\end {align*}

These $x_i$ and $y_i$ are noncommutative versions of the same monomials given in the commutative case,
and $z_k$ are paths connecting the upper-left corner of the first square face to the upper-right corner
of the $k^\mathrm{th}$ square face. 
We may again perform the exact same sequence of gauge transformations as in the commutative case to
arrive at the weights depicted in \textbf{Figure \ref{fig:xy_quiver}}. However, these weights 
will now be noncommutative versions of the same monomials. 
In terms of the monomials defined above, the noncommutative weights we obtain
in \textbf{Figure \ref{fig:xy_quiver}} are given by

\begin {align*}
    X_i &= z_{i-2} \, x_i \, z_{i-2}^{-1} \\
    Y_i &= z_{i-2} \, y_i \, z_{i-2}^{-1} \\
    Z   &= z_n c_n
\end {align*}

Just as in the commutative case, these can be interpreted as closed loops in the quiver. The difference now
in the non-commutative case is that all these loops share a common basepoint. This common basepoint is
the upper-left corner of the first square face. As in the previous section, let $\mathscr{L}_\bullet \subset \mathscr{L}$ denote the subspace
of loops based at this point.
In fact, it is the group algebra of the fundamental group of the ribbon surface of the graph,
and it is generated by the $X_i$, $Y_i$, and $Z$.

\bigskip

A simple calculation shows that $Z$ is a Casimir of the double bracket.
The bracket is then well-defined on the quotient $\mathscr{L}_\bullet^{(1)} := \mathscr{L}_\bullet/\left<Z-1\right>$, where $Z$ is set to $1$.
The induced brackets in $\mathscr{L}_\bullet^{(1)}$ are given by:

\begin {alignat*}{2}
    \left<X_{i+1},X_i\right> &= X_{i+1}X_i &\hspace{1cm} \left<Y_{i+1},Y_i\right> &= Y_{i+1}Y_i \\
    \left<Y_{i+2},Y_i\right> &= Y_{i+2}Y_i &\hspace{1cm} \left<Y_i,X_i\right> &= Y_iX_i \\
    \left<Y_{i+1},X_i\right> &= Y_{i+1}X_i &\hspace{1cm} \left<X_{i+1},Y_i\right> &= X_{i+1}Y_i \\
    \left<X_{i+2},Y_i\right> &= X_{i+2}Y_i &\hspace{1cm} &
\end {alignat*}

\medskip

If we don't set $Z=1$, then the bracket relations in $\mathscr{L}_\bullet$ are mostly the same, except some exceptions when $i=1$ or $i=2$:

\begin {alignat*}{2}
    \left<X_3,X_2\right> &= X_3 Z^{-1} X_2 Z &\hspace{1cm} \left<Y_3,Y_2\right> &= Y_3 Z^{-1} Y_2 Z \\
    \left<Y_3,X_2\right> &= Y_3 Z^{-1} X_2 Z &\hspace{1cm} \left<X_3,Y_2\right> &= X_3 Z^{-1} Y_2 Z \\
    \left<Y_3,Y_1\right> &= Y_3 Z^{-1} Y_1 Z &\hspace{1cm} \left<Y_4,Y_2\right> &= Y_4 Z^{-1} Y_2 Z \\
    \left<X_3,Y_1\right> &= X_3 Z^{-1} Y_1 Z &\hspace{1cm} \left<X_4,Y_2\right> &= X_4 Z^{-1} Y_2 Z
\end {alignat*}

\medskip

If we perform the same sequence of Postnikov moves as in the commutative case, followed by gauge transformations, we get a graph isomorphism,
with the same edges having weight 1 as before. This gives the transformation

\begin {align*}
    X_i &\mapsto (X_i+Y_i)^{-1} X_i (X_{i+2}+Y_{i+2}) \\
    Y_i &\mapsto (X_{i+1}+Y_{i+1})^{-1} Y_{i+1} (X_{i+3}+Y_{i+3})
\end {align*}

Just as in the commutative case, this is almost the expression derived earlier for the pentagram map. It differs only by a shift in the $Y$-indices.

\medskip

\begin {rmk}
    In the formula above, the indices are \emph{not} read cyclically. If one of the indices is greater than $n$, we must conjugate by $Z$. For instance,
    if $i=n-2$, and $i+3=n+1$, then by $X_{i+3}$ we really mean $Z X_1 Z^{-1}$, and similarly for $Y$.
\end {rmk}

\bigskip

\subsection {The $P,Q$ Variables}

We now define a non-commutative version of the $p,q$ variables from the classical case. They are given by

\begin {alignat*}{2}
    p_i &= b_ic_i^{-1}d_i^{-1}a_i^{-1} &\hspace{1cm} q_i &= c_{i-2}d_{i-1}a_{i+1}b_i^{-1} \\
\end {alignat*}
Conjugating by the same $z_k$ paths as for the $X,Y$ variables, we obtain based versions:
\begin {alignat*}{2}
    P_i &= z_{i-2} \, p_i \, z_{i-2}^{-1} &\hspace{1cm} Q_i &= z_{i-2} \, q_i \, z_{i-2}^{-1} \\
\end {alignat*}
Similar to the commutative case, we have the following relation with the $X,Y$ variables:
\begin {alignat*}{2}
    P_i &= Y_iX_i^{-1} &\hspace{1cm} Q_i &= X_{i+1}Y_i^{-1} \\
\end {alignat*}
Their induced brackets are given by 
\begin {alignat*}{2}
    \left< Q_i, P_i \right> &= Q_i P_i &\hspace{1cm} \left< P_i, Q_{i-1} \right> &= P_i Q_{i-1} \\
    \left< Q_{i+1}, P_i \right> &= Q_{i+1} Y_i^{-1} P_i Y_i &\hspace{1cm} \left< P_i, Q_{i+2} \right> &= Y_i^{-1} P_i Y_i Q_{i+2}
\end {alignat*}

From the formulas given in the previous section for the pentagram map in the $X,Y$ variables, we get that the pentagram map transorms the $P,Q$ variables by

\begin {align*}
    P_i &\mapsto X_{i+1}^{-1} (1+P_{i+1}^{-1})^{-1} X_{i+1} \cdot (1+P_{i+3}) Q_{i+2} (1+P_{i+2}) \cdot Y_i^{-1} (1+P_i) Y_i \\
    Q_i &\mapsto X_{i+1}^{-1} P_{i+1}^{-1} X_{i+1}
\end {align*}

\section {Invariance, Invariants, and Integrability}

\bigskip

\subsection {Invariance of the Induced Bracket}

We will show in this section that the induced bracket $\left< -,- \right>$ on $\mathscr{L}_{\bullet}^\natural$
is invariant under the pentagram map. To do so,
we will consider step-by-step how the weights and double bracket change under the Postnikov moves. 

\bigskip

Recall that, starting with the $X,Y$ variables, the sequence of Postnikov moves which gives the pentagram map is as follows:
\begin {enumerate}
    \item Perform the square move at each square face
    \item Perform the white-swap move at each edge connecting two white vertices
    \item Perform the black-swap move at each edge connecting two black vertices
    \item Perform gauge transformations so all weights are 1 except the bottom and left of each square face
\end {enumerate}

\bigskip

We will actually start with
the $a,b,c,d$ weights instead of the $X,Y$ weights. The result will be the same, as we will show below.
Application of a square move gives

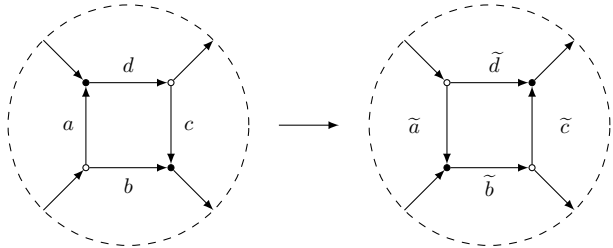
\begin {center}
\begin {tikzpicture}[scale=0.8, every node/.style={scale=0.8}]
    % Before
    \draw [dashed]     (0,0) circle (2cm);
    \draw (0.707,0.707) circle (0.05cm);
    \draw[fill=black] (0.707,-0.707) circle (0.05cm);
    \draw[fill=black] (-0.707,0.707) circle (0.05cm);
    \draw (-0.707,-0.707) circle (0.05cm);

    \draw [-latex] (-1.414,1.414) -- (-0.75,0.75);
    \draw [-latex] (0.75,0.75) -- (1.414,1.414);
    \draw [-latex] (0.75,-0.75) -- (1.414,-1.414);
    \draw [-latex] (-1.414,-1.414) -- (-0.74,-0.74);

    \draw [-latex] (-0.65,0.707) -- (0.65,0.707);
    \draw [-latex] (-0.65,-0.707) -- (0.65,-0.707);
    \draw [-latex] (-0.707,-0.65) -- (-0.707,0.65);
    \draw [-latex] (0.707,0.65) -- (0.707,-0.65);

    \draw (-1,0) node {$a$};
    \draw (0,-1) node {$b$};
    \draw (1,0) node {$c$};
    \draw (0,1) node {$d$};

    % Arrow

    \draw [-latex] (2.5,0) -- (3.5,0);

    % After
    \draw [dashed]    (6+0,0)           circle (2cm);
    \draw[fill=black] (6+0.707,0.707)   circle (0.05cm);
    \draw[fill=white] (6+0.707,-0.707)  circle (0.05cm);
    \draw[fill=white] (6+-0.707,0.707)  circle (0.05cm);
    \draw[fill=black] (6+-0.707,-0.707) circle (0.05cm);

    \draw [-latex] (6+-1.414,1.414)  -- (6+-0.75,0.75);
    \draw [-latex] (6+0.75,0.75)     -- (6+1.414,1.414);
    \draw [-latex] (6+0.75,-0.75)    -- (6+1.414,-1.414);
    \draw [-latex] (6+-1.414,-1.414) -- (6+-0.74,-0.74);

    \draw [-latex] (6+-0.65,0.707)  -- (6+0.65,0.707);
    \draw [-latex] (6+-0.65,-0.707) -- (6+0.65,-0.707);
    \draw [-latex] (6+-0.707,0.65) -- (6+-0.707,-0.65);
    \draw [-latex] (6+0.707,-0.65)   -- (6+0.707,0.65);

    \draw (5-0.25,0) node {$\widetilde{a}$};
    \draw (6+0,-1) node {$\widetilde{b}$};
    \draw (7+0.25,0)  node {$\widetilde{c}$};
    \draw (6+0.08,1.05)  node {$\widetilde{d}$};

\end {tikzpicture}
\end {center}

The new edge weights are given by

\begin {alignat*}{2}
    \widetilde{b} &= b+adc &\hspace{1cm}
    \widetilde{a} &= dc \widetilde{b}^{-1}  \\
    \widetilde{c} &= \widetilde{b}^{-1} ad &\hspace{1cm}
    \widetilde{d} &= \widetilde{a} bc^{-1} 
\end {alignat*}

The double brackets of these new weights are

\[
    \begin {array}{ccc|ccc}
        \db{\widetilde{b}, \widetilde{a}} &=& \frac{1}{2} \, \left( \widetilde{a} \widetilde{b} \otimes 1 \right) &
        \db{\widetilde{a}, \widetilde{d}} &=& \frac{1}{2} \, \left( \widetilde{a} \otimes \widetilde{d} \right) \\[1.2ex]
        \db{\widetilde{b}, \widetilde{c}} &=& \frac{1}{2} \, \left( 1 \otimes \widetilde{b} \widetilde{c} \right) &
        \db{\widetilde{c}, \widetilde{d}} &=& \frac{1}{2} \, \left( \widetilde{d} \otimes \widetilde{c} \right)
    \end {array}
\]

\bigskip

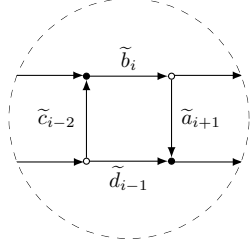
\begin {figure}[h!]
\centering
\caption {Edge weights after white-swap and black-swap moves}
\label {fig:after_postnikov}
\begin {tikzpicture}[scale=0.8, every node/.style={scale=0.8}]
    \clip (0,0) circle (2cm);
    \draw [dashed]     (0,0) circle (2cm);
    \draw (0.707,0.707) circle (0.05cm);
    \draw[fill=black] (0.707,-0.707) circle (0.05cm);
    \draw[fill=black] (-0.707,0.707) circle (0.05cm);
    \draw (-0.707,-0.707) circle (0.05cm);

    \draw [-latex] (-2,0.72) -- (-0.75,0.72);
    \draw [-latex] (0.75,0.72) -- (1.9,0.72);
    \draw [-latex] (0.75,-0.72) -- (1.9,-0.72);
    \draw [-latex] (-1.9,-0.72) -- (-0.74,-0.72);

    \draw [-latex] (-0.65,0.707) -- (0.65,0.707);
    \draw [-latex] (-0.65,-0.707) -- (0.65,-0.707);
    \draw [-latex] (-0.707,-0.65) -- (-0.707,0.65);
    \draw [-latex] (0.707,0.65) -- (0.707,-0.65);

    \draw (-1.2,0) node {$\widetilde{c}_{i-2}$};
    \draw (0,-1) node {$\widetilde{d}_{i-1}$};
    \draw (1.2,0) node {$\widetilde{a}_{i+1}$};
    \draw (0,1) node {$\widetilde{b}_i$};
\end {tikzpicture}
\end {figure}

Next, applying the white-swap and black-swap moves interchanges the square and octagonal faces.
The resulting weights around the square faces is pictured in \textbf{Figure \ref{fig:after_postnikov}}.

Finally, we perform gauge transformations, as we did to get the original $X,Y$ weights, so that all
weights become 1 except those on the bottom and left edges of the square faces. We call these
resulting weights $\widetilde{X}$ and $\widetilde{Y}$. If we define the ``staircase'' monomials
$\xi_k = \widetilde{a}_1 \widetilde{b}_1 \cdots \widetilde{a}_k \widetilde{b}_k$, then we can express
the new weights as

\begin {alignat*}{2}
    \widetilde{X}_i &= \xi_i \, \widetilde{c}_i \widetilde{a}_{i+2}^{-1} \widetilde{b}_{i+1}^{-1} \widetilde{a}_{i+1}^{-1} \, \xi_i^{-1} &\hspace{1cm}
    \widetilde{Y}_i &= \xi_i \, \widetilde{d}_{i+1} \widetilde{a}_{i+3}^{-1}  \widetilde{b}_{i+2}^{-1} \widetilde{a}_{i+2}^{-1} \widetilde{b}_{i+1}^{-1} \widetilde{a}_{i+1}^{-1} \, \xi_i^{-1} 
\end {alignat*}

These weights are the images of $X_i$ and $Y_i$ under the pentagram map. That is, $\widetilde{X}_i = T(X_i)$ and $\widetilde{Y}_i = T(Y_i)$.
As noted in the previous sections, these are given in terms of the original $X,Y$ variables as $\widetilde{X}_i = \sigma_i^{-1} X_i \sigma_{i+2}$
and $\widetilde{Y}_i = \sigma_{i+1}^{-1} Y_{i+1} \sigma_{i+3}$, where $\sigma_i = X_i+Y_i$. We are now ready to prove the main theorem of this section,
but first we introduce some notation.

\medskip

Define the map $S \colon \mathscr{L}_\bullet \to \mathscr{L}_\bullet$ which shifts all the indices. That is, $S(X_i) = X_{i-1}$ and $S(Y_i) = Y_{i-1}$,
where the indices are read cyclically.

\bigskip

\begin {thm} \label{thm:invariance} \ \relax \\
    \begin {tabular}{cl}
    $(a)$ & In $\mathscr{L}_\bullet^{(1)}$, the induced bracket $\left<-,-\right>$ is invariant under the pentagram map. \\[1.2ex]
    $(b)$ & In $\mathscr{L}_\bullet$, the induced bracket is invariant under $S^2 \circ T$.
    \end {tabular}
\end {thm}
\begin {proof}
    We want to show that the induced brackets of $\widetilde{X}_i$ and $\widetilde{Y}_i$ have exactly the same form as the induced brackets
    of the $X_i$ and $Y_i$. Note that for $a,b,c \in \mathscr{L}_\bullet$, we have $\left<a,c\right> = \left<bab^{-1},c\right>$,
    since $a^\natural = (bab^{-1})^\natural$. So we will instead compute the brackets of the conjugate, but equivalent, elements:

    \begin {alignat*}{2}
        \widetilde{x}_i &= \xi_i^{-1} \widetilde{X}_i \xi_i = \widetilde{c}_i \widetilde{a}_{i+2}^{-1} \widetilde{b}_{i+1}^{-1} \widetilde{a}_{i+1}^{-1} &\hspace{1cm}
        \widetilde{y}_i &= \xi_i^{-1} \widetilde{Y}_i \xi_i = \widetilde{d}_{i+1} \widetilde{a}_{i+3}^{-1}  \widetilde{b}_{i+2}^{-1} \widetilde{a}_{i+2}^{-1} \widetilde{b}_{i+1}^{-1} \widetilde{a}_{i+1}^{-1}
    \end {alignat*}

    We will compute the three possible combinations $\db{x_i,x_j}$, $\db{y_i,y_j}$, and $\db{x_i,y_j}$.
    For $\db{x_i,x_j}$, a simple calculation using the Leibniz rules for double brackets gives:
    \begin {align*}
        \db{\widetilde{x}_i, \widetilde{x}_j} &= \frac{1}{2} \delta_{i,j+1} (\widetilde{c}_j \widetilde{a}_{j+2}^{-1} \widetilde{x}_i \otimes \widetilde{a}_{j+2} \widetilde{c}_j^{-1} \widetilde{x}_j
                                                                     + \widetilde{x}_j \widetilde{a}_i \widetilde{b}_i \otimes \widetilde{x}_i \widetilde{b}_i^{-1} \widetilde{a}_i^{-1}) \\
        &\phantom{=} - \frac{1}{2} \delta_{i,j-1} (\widetilde{x}_j \widetilde{b}_j^{-1} \widetilde{a}_j^{-1} \otimes \widetilde{x}_i \widetilde{a}_j \widetilde{b}_j
                                                   + \widetilde{b}_j^{-1} \widetilde{a}_j^{-1} \otimes \widetilde{x}_i \widetilde{a}_j \widetilde{b}_j \tilde{x}_j)
    \end {align*}
    Therefore their induced brackets are
    \[
        \left<\widetilde{x}_i, \widetilde{x}_j\right> = \delta_{i,j+1} (\widetilde{c}_j \widetilde{a}_{j+2}^{-1}) \, \widetilde{x}_i \, (\widetilde{a}_{j+2} \widetilde{c}_j^{-1}) \, \widetilde{x}_j
        - \delta_{i,j-1} (\widetilde{b}_j^{-1} \widetilde{a}_j^{-1}) \, \widetilde{x}_i \, (\widetilde{a}_j \widetilde{b}_j) \, \tilde{x}_j
    \]
    Conjugating the right hand side by $\xi_j$ gives the equivalent formula in terms of the $\widetilde{X}_i$:
    \[ \left< \widetilde{X}_i, \widetilde{X}_j \right> = (\delta_{i,j+1} - \delta_{i,j-1}) \widetilde{X}_i \widetilde{X}_j \]
    This is exactly the same as the bracket formula for the original $X_i$. There are, however, some exceptions. When $i=n$,
    we instead have
    \[ \left< \widetilde{X}_1, \widetilde{X}_n \right> = \widetilde{X}_1 Z^{-1} \widetilde{X}_n Z \]
    
    \medskip

    The calculations for $\left< \widetilde{Y}_i, \widetilde{Y}_j \right>$
    and $\left< \widetilde{X}_i, \widetilde{Y}_j \right>$ are similar, with exceptions when $i=n$ and $i=n-1$. So the exact form of the bracket
    relations are not invariant, since the exceptional cases (which are conjugated by $Z$) do not occur at the same indices.
    But if we set $Z=1$, then the bracket is invariant in $\mathscr{L}_\bullet^{(1)}$. Also, if we shift the indices by 2
    and define $\widehat{X}_i := S^2 (\widetilde{X}_i) = S^2 \circ T (X_i)$, then the bracket is invariant under $X_i \mapsto \widehat{X}_i$.
\end {proof}

\bigskip

\begin {rmk}
    \textbf{Theorem \ref{thm:invariance}} suggests that a more appropriate definition of the pentagram map
    should include a ``rotation'', which shifts the labels/indices of the vertices by 2. This shift of indices
    is necessary to make the Poisson structure invariant. From now on, we will use the notation 
    $\widehat{T} \colon \mathcal{GP}_{n,N} \to \mathcal{GP}_{n,N}$ for this modified version of the map.
    It is worth noting that earlier definitions of the pentagram map in \cite{ost10} \cite{schwartz92} \cite{gsvt_16}
    do not include this shift of indices, and this phenomenon only appears by considering this noncommutative generalization.
\end {rmk}

\bigskip

\subsection {The Invariants}

Let $Q = Q_n$ be the quiver/network for the pentagram map, with the $X_i$ and $Y_i$ weights as before, 
and $\mathscr{B} = \mathscr{B}_Q(\lambda)$ its boundary measurement matrix, whose entries are 
elements of the Laurent polynomial ring $\mathscr{L}_\bullet[\lambda^\pm]$.
For each $i \geq 1$, denote the coefficients of $\lambda^k$ in $\mathrm{tr}(\mathscr{B}^i)$ by $t_{ik}$:
\[ \mathrm{tr}(\mathscr{B}^i) = \sum_k t_{ik} \lambda^k \]
We will spend the remainder of this section proving the following theorem.

\begin {thm} \label{thm:invariants}
    For each $i$ and $j$, the classes $t_{ij}^\natural$ are invariant under both $T$ and $\widehat{T}$.
\end {thm}

\bigskip

To begin proving this, we first make the following simple observation:

\bigskip

\begin {lem}
    Let $R$ be an associative ring, and $p, q \in R[\lambda^\pm]$. Then every coefficient of $[p,q]$ is in $[R,R]$.
\end {lem}
\begin {proof}
    Let $p = \sum_i p_i \lambda^i$ and $q = \sum_j q_j \lambda^j$. Then $[p,q] = \sum_{i,j} [p_i,q_j] \lambda^{i+j}$.
\end {proof}

\medskip

This implies the following

\medskip

\begin {cor} \label{cor:coefs_mod_commutators}
    Let $f = \sum_i f_i \lambda^i$ and $g = \sum_j g_j \lambda^j$ be Laurent polynomials in $A := R[\lambda^\pm]$.
    If $f \equiv g ~ \mathrm{mod} ~ [A,A]$, then $f_i \equiv g_i ~ \mathrm{mod} ~ [R,R]$ for each $i$.
\end {cor}

\bigskip

We now consider traces of powers of matrices over general rings:

\bigskip

\begin {lem} \label{lem:traces_of_powers}
    Let $R$ be an associative ring, and $A,B \in \mathrm{Mat}_n(R)$.
    Then $\mathrm{tr}((AB)^k) \equiv \mathrm{tr}((BA)^k) ~ \mathrm{mod} ~ [R,R]$ for all $k$.
\end {lem}
\begin {proof}
    If the matrices are given by $A = (a_{ij})$ and $B = (b_{ij})$ then
    \begin {align*} 
        \mathrm{tr}(AB)^k &= \sum_{i,j_1,\dots,j_{k-1}} (AB)_{ij_1} (AB)_{j_1j_2} \cdots (AB)_{j_{k-1}i} \\
        &= \sum_{i,j_1,\dots,j_{k-1}} \left( \sum_{\ell_1} a_{i\ell_1}b_{\ell_1 j_1} \right) \cdots \left( \sum_{\ell_k} a_{j_{k-1} \ell_k} b_{\ell_k i} \right) \\
        &= \sum_{i,j_1,\dots,j_{k-1}} \sum_{\ell_1,\dots,\ell_k} a_{i\ell_1}b_{\ell_1 j_1} \cdots a_{j_{k-1} \ell_k} b_{\ell_k i}
    \end {align*}
    Similarly, using $BA$ instead, we get
    \[ \mathrm{tr}(BA)^k = \sum_{i,j_1,\dots,j_{k-1}} \sum_{\ell_1,\dots,\ell_k} b_{i\ell_1}a_{\ell_1 j_1} \cdots b_{j_{k-1} \ell_k} a_{\ell_k i} \]
    After re-indexing by $\ell_t \mapsto j_t$ (for $1 \leq t \leq k-1$), $\ell_k \mapsto i$, $i \mapsto \ell_1$, $j_t \mapsto \ell_{t+1}$, we get
    \[ \mathrm{tr}(BA)^k = \sum_{i,j_1,\dots,j_{k-1}} \sum_{\ell_1,\dots,\ell_k} b_{\ell_1 j_1}a_{j_1 \ell_2} \cdots b_{\ell_k i} a_{i \ell_1} \]
    Clearly, this is the same as the expression for $\mathrm{tr}((AB)^k)$ mod $[R,R]$, since each term with corresponding indices
    is the same after cyclicly shifting the last $a_{i \ell_1}$ to the beginning.
\end {proof}

\bigskip

Now we may prove \textbf{Theorem \ref{thm:invariants}}:
\begin {proof}
    As discussed earlier, the map $T$ is given, in the $X,Y$ variables, by a sequence of Postnikov moves
    and gauge transformations.
    The square move and gauge transformations do not change the boundary measurement matrix at all, and the 
    ``white-swap'' and ``black-swap'' moves only change the boundary measurement matrix up to conjugation,
    since it amounts to a refactorization $\mathscr{B} = \mathscr{B}_1 \mathscr{B_2} \mapsto \widetilde{\mathscr{B}} = \mathscr{B}_2 \mathscr{B}_1$.
    If we let $R = \mathscr{L}_\bullet$ and $A = \mathscr{L}_\bullet[\lambda^\pm]$, 
    then by \textbf{Corollary \ref{cor:coefs_mod_commutators}} and \textbf{Lemma \ref{lem:traces_of_powers}}, 
    $\mathrm{tr}(\widetilde{\mathscr{B}}^k)$ and $\mathrm{tr}(\mathscr{B}^k)$ differ by an
    element of $[A,A]$, and so $t_{ij}^\natural$ is invariant under $T$. The map $\widehat{T}$ is just
    $T$ followed by a shift of indices, which is equivalent to cutting part of the network off one end,
    and glueing it to the opposite end. Again, this amounts to a re-factorization of the boundary measurement
    measurement matrix, which changes it only up to conjugation.
\end {proof}

\bigskip

\begin {rmk}
    In \cite{izosimov_18}, Izosimov interprets the pentagram map, as well as some generalizations,
    in terms of re-factorizations in Poisson-Lie groups. At the end of the paper, he poses the question
    of whether his techniques, when applied to matrix-valued coefficients, give rise to the Grassman
    pentagram map of Mari Beffa and Felipe. It seems that the proof of \textbf{Theorem \ref{thm:invariants}},
    which realizes the pentagram map as a re-factorization of the boundary measurement matrix,
    suggests a positive answer to Izosimov's question.
\end {rmk}

\bigskip

\subsection {Involutivity of the Invariants}

\bigskip

In this section, we will prove that the invariants from the previous section are an involutive family with
respect to the induced bracket. More specifically:

\begin {thm} \label{thm:involutivity}
    Let $Q_n$ be the network for the pentagram map on $\mathcal{GP}_{n,N}$, and $\mathscr{B} = \mathscr{B}_Q(\lambda)$ its boundary measurement matrix,
    with $t_{ik}$ the homogeneous components of $\mathrm{tr}(\mathscr{B}^i)$ as defined before. 
    Then for all $i,j,k,\ell$:
    \[ \left< t_{ik}, t_{j\ell} \right> = 0 \]
\end {thm}

\medskip

Recall that we use as the cut the identified top/bottom edge of the rectangle on which we draw the quiver.
So the element $t_{ik} \in \mathscr{L}$ is the sum over all loops in $Q$ which are homologous to $(i,k)$ cycles on the torus. That is, if we lift
these paths to the universal cover, then they cross fundamental domains $i$ times horizontally and $k$ times vertically (with sign).
Let $A_{ik}$ be the set of all such paths, so that
\[ t_{ik} = \sum_{p \in A_{ik}} p \]
Then by bilinearity and the formula from \textbf{Theorem \ref{thm:bracket_formula}},
\[ \left< t_{ik}, t_{j\ell} \right> = \sum_{f \in A_{ik}} \sum_{g \in A_{j\ell}} \left< f, g \right> 
    = \sum_{f \in A_{ik}} \sum_{g \in A_{j\ell}} \sum_{\bullet \in f \cap g} \varepsilon_\bullet(f,g) f_\bullet g_\bullet
\]

We want to prove that this expression is zero. To do so, we will define a sign-reversing, fixed-point-free permutation on
the set of terms appearing in the sum. That is, each term which appears can be paired with another
term with opposite coefficient and the same monomial in $A^\natural$.

\bigskip

We start by choosing an arbitrary non-zero term in the sum, of the form $\varepsilon_\bullet(f,g) \, f_\bullet g_\bullet$.
To define a permutation as suggested above, we need to know what are the other terms in the sum
that are equivalent to $f_\bullet g_\bullet$ in $A^\natural$.
The answer is given by the following lemmas:

\bigskip

\begin {lem}
    Suppose $f_\bullet g_\bullet \equiv f'_\ast g'_\ast$ mod $[\mathscr{L}_\bullet,\mathscr{L}_\bullet]$,
    and that $\varepsilon_\ast(f',g') \neq 0$.
    Then $\ast$ is in either $f \cap g$, $f \cap f$, or $g \cap g$.
\end {lem}
\begin {proof}
    Certainly $f'$ and $g'$ together have the same combined set of edges as $f$ and $g$ since $f_\bullet g_\bullet = f'_\ast g'_\ast$
    cyclically. Then since two paths come together at $\ast$ (since $\varepsilon_\ast(f',g') \neq 0$), it must be that
    either $f$ and $g$ meet at $\ast$, or $f$ or $g$ meets itself.
\end {proof}

\bigskip

\begin {lem}
    Let $\bullet$ and $\ast$ be as in the previous lemma, with $f_\bullet g_\bullet = f'_\ast g'_\ast$.
    Additionally assume that $(f,g)$ and $(f',g')$ are both in $A_{ik} \times A_{j\ell}$. Then \\
    \begin {tabular}{cl}
        $(a)$ & If $\ast \in f \cap g$, then $(f',g')$ are obtained from $(f,g)$ by swapping the segments
                of $f$ and $g$ between $\bullet$ and $\ast$. \\[1.2ex]
        $(b)$ & If $\ast \in f \cap f$, then $(f',g')$ are obtained from $(f,g)$ by swapping the entire loop $g$
                with the subloop of $f$ based at $\ast$.
    \end {tabular}
\end {lem}
\begin {proof}
    $(a)$ Denote by $a$ and $b$ the segments of $f$
    between $\bullet$ and $\ast$, and similarly let $x$ and $y$ be the segments of $g$ between $\bullet$ and $\ast$, 
    so we may write $f = \bullet a \ast b$ and $g = \bullet x \ast y$. Then 
    \[ f_\bullet g_\bullet = \bullet a \ast b \bullet x \ast y \]
    Since we assumed that $f'_\ast g'_\ast = f_\bullet g_\bullet$ up to cyclic permutation, 
    then we can write it starting at $\ast$ as
    \[ f'_\ast g'_\ast = \ast b \bullet x \ast y \bullet a \]
    This means $f' = \ast b \bullet x$ and $g' = \ast y \bullet a$, or the other way around: $f' = \ast y \bullet a$
    and $g' = \ast b \bullet x$. The two possibilities differ by switching the roles of $f'$ and $g'$.
    But since we assume $(f',g') \in A_{ik} \times A_{j\ell}$, and since $A_{ik} \neq A_{j\ell}$,
    only one of the possibilities will be correct.

    \medskip

    After cyclically permuting, we have
    $f' = \bullet x \ast b$ and $g' = \bullet a \ast y$ (or the other way around).
    Thus we see that $f'$ and $g'$ are obtained from $f$ and $g$ by swapping either $x$ with $a$
    or $y$ with $b$,
    which are the portions of the paths in between their common subpaths $\bullet$ and $\ast$.
    This is illustrated in \textbf{Figure \ref{fig:type_I_swap}}.
    Note that since $(f',g') \in A_{ik} \times A_{j\ell}$, the swapped segments must have the
    same intersection index with the cut.

    
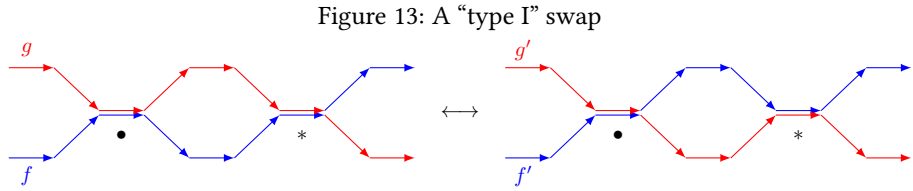
\begin {figure}[h]
    \centering
    \caption {A ``type I'' swap}
    \label {fig:type_I_swap}
    \begin {tikzpicture}[scale=0.6, every node/.style={scale=0.8}]
        % f
        \draw [-latex, blue] (0,-1)    -- (1,-1);
        \draw [-latex, blue] (1,-1)    -- (2,-0.05);
        \draw [-latex, blue] (2,-0.05) -- (3,-0.05);
        \draw [-latex, blue] (3,-0.05) -- (4,-1);
        \draw [-latex, blue] (4,-1)    -- (5,-1);
        \draw [-latex, blue] (5,-1)    -- (6,-0.05);
        \draw [-latex, blue] (6,-0.05) -- (7,-0.05);
        \draw [-latex, blue] (7,0.05) -- (8,1);
        \draw [-latex, blue] (8,1)    -- (9,1);
        \draw [blue] (0.4,-1.4) node {$f$};

        % g
        \draw [-latex, red] (0,1)    -- (1,1);
        \draw [-latex, red] (1,1)    -- (2,0.05);
        \draw [-latex, red] (2,0.05) -- (3,0.05);
        \draw [-latex, red] (3,0.05) -- (4,1);
        \draw [-latex, red] (4,1)    -- (5,1);
        \draw [-latex, red] (5,1)    -- (6,0.05);
        \draw [-latex, red] (6,0.05) -- (7,0.05);
        \draw [-latex, red] (7,-0.05) -- (8,-1);
        \draw [-latex, red] (8,-1)    -- (9,-1);
        \draw [red] (0.4,1.4) node {$g$};

        \draw [black] (2.5,-0.5) node {$\bullet$};
        \draw [black] (6.5,-0.5) node {$\ast$};

        % ---------------------
        \draw (10,0) node {$\longleftrightarrow$};
        % ---------------------

        % f'
        \draw [-latex, blue] (11 + 0,-1)    -- (11 + 1,-1);
        \draw [-latex, blue] (11 + 1,-1)    -- (11 + 2,-0.05);
        \draw [-latex, blue] (11 + 2,-0.05) -- (11 + 3,-0.05);

        \draw [-latex, blue] (11 + 3,0.05) -- (11 + 4,1);
        \draw [-latex, blue] (11 + 4,1)    -- (11 + 5,1);
        \draw [-latex, blue] (11 + 5,1)    -- (11 + 6,0.05);
        \draw [-latex, blue] (11 + 6,0.05) -- (11 + 7,0.05);

        \draw [-latex, blue] (11 + 7,0.05)  -- (11 + 8,1);
        \draw [-latex, blue] (11 + 8,1)     -- (11 + 9,1);
        \draw [blue]         (11 + 0.4,-1.4) node {$f'$};

        % g'
        \draw [-latex, red] (11 + 0,1)     -- (11 + 1,1);
        \draw [-latex, red] (11 + 1,1)     -- (11 + 2,0.05);
        \draw [-latex, red] (11 + 2,0.05)  -- (11 + 3,0.05);

        \draw [-latex, red] (11 + 3,-0.05)  -- (11 + 4,-1);
        \draw [-latex, red] (11 + 4,-1)     -- (11 + 5,-1);
        \draw [-latex, red] (11 + 5,-1)     -- (11 + 6,-0.05);
        \draw [-latex, red] (11 + 6,-0.05)  -- (11 + 7,-0.05);

        \draw [-latex, red] (11 + 7,-0.05) -- (11 + 8,-1);
        \draw [-latex, red] (11 + 8,-1)    -- (11 + 9,-1);
        \draw [red]         (11 + 0.4,1.4) node {$g'$};

        \draw [black] (11 + 2.5,-0.5) node {$\bullet$};
        \draw [black] (11 + 6.5,-0.5) node {$\ast$};
    \end {tikzpicture}
    \end {figure}

    \bigskip

    $(b)$ Now suppose $\ast \in f \cap f$. Let $a,b,c$ be the segments of $f$ between $\bullet$ and $\ast$, so that
    $f = \bullet a \ast b \ast c$. Also let $x$ be the rest of $g$ after $\bullet$, so that $g = \bullet x$. Then
    \[ f_\bullet g_\bullet = \bullet a \ast b \ast c \bullet x \]
    Cyclically permuting, and using the assumption that $f_\bullet g_\bullet = f'_\ast g'_\ast$, we get
    \[ f'_\ast g'_\ast = \ast c \bullet x \bullet a \ast b \]
    Again there are two possibilities. Either $f' = \ast c \bullet x \bullet a$ and $g' = \ast b$,
    or $f'$ and $g'$ are flipped. But the condition that $f' \in A_{ik}$ and $g' \in A_{j\ell}$ ensures
    that only one of the two choices is correct. Suppose it is the first case, and $f' = \ast c \bullet x \bullet a$
    and $g' = \ast b$. Then $f'$ is obtained from $f$ by removing the loop based at $\ast$ and adding the loop $g$,
    and $g'$ is simply the subloop of $f$ based at $\ast$. This is illustrated in \textbf{Figure \ref{fig:type_II_swap}}.
    Note that since $(f',g') \in A_{ik} \times A_{j\ell}$,
    the loops $\bullet x$ and $\ast b$ must cross the cut and rim the same number of times.

    
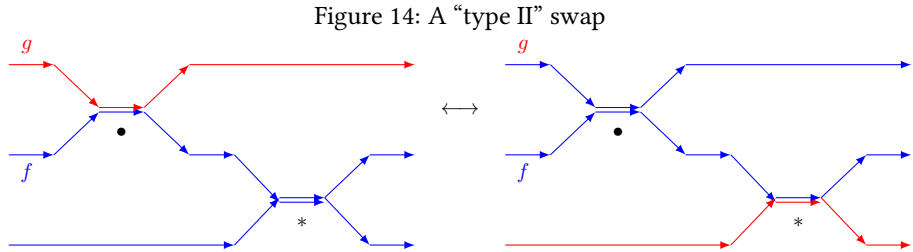
\begin {figure}[h]
    \centering
    \caption {A ``type II'' swap}
    \label {fig:type_II_swap}
    \begin {tikzpicture}[scale=0.6, every node/.style={scale=0.8}]
        % f
        \draw [-latex, blue] (0,-1)    -- (1,-1);
        \draw [-latex, blue] (1,-1)    -- (2,-0.05);
        \draw [-latex, blue] (2,-0.05) -- (3,-0.05);
        \draw [-latex, blue] (3,-0.05) -- (4,-1);
        \draw [-latex, blue] (4,-1)    -- (5,-1);
        \draw [-latex, blue] (5,-1)    -- (6,-2.05);
        \draw [-latex, blue] (6,-2.05) -- (7,-2.05);
        \draw [-latex, blue] (7,-2 + 0.05) -- (8,-3);
        \draw [-latex, blue] (8,-3)    -- (9,-3);

        \draw [-latex, blue] (0,-3) -- (5,-3);
        \draw [-latex, blue] (5,-3) -- (6,-2 + 0.05);
        \draw [-latex, blue] (6,-2 + 0.05) -- (7,-2 + 0.05);
        \draw [-latex, blue] (7,-2+0.05) -- (8,-1);
        \draw [-latex, blue] (8,-1) -- (9,-1);

        \draw [blue] (0.4,-1.4) node {$f$};

        % g
        \draw [-latex, red] (0,1)    -- (1,1);
        \draw [-latex, red] (1,1)    -- (2,0.05);
        \draw [-latex, red] (2,0.05) -- (3,0.05);
        \draw [-latex, red] (3,0.05) -- (4,1);
        \draw [-latex, red] (4,1)    -- (9,1);
        \draw [red] (0.4,1.4) node {$g$};

        \draw [black] (2.5,-0.5) node {$\bullet$};
        \draw [black] (6.5,-2.5) node {$\ast$};

        % ---------------------
        \draw (10,0) node {$\longleftrightarrow$};
        % ---------------------

        % f'
        \draw [-latex, blue] (11 + 0,-1)    -- (11 + 1,-1);
        \draw [-latex, blue] (11 + 1,-1)    -- (11 + 2,-0.05);
        \draw [-latex, blue] (11 + 2,-0.05) -- (11 + 3,-0.05);
        \draw [-latex, blue] (11 + 3,-0.05) -- (11 + 4,-1);
        \draw [-latex, blue] (11 + 4,-1)    -- (11 + 5,-1);
        \draw [-latex, blue] (11 + 5,-1)    -- (11 + 6,-2.05);
        \draw [-latex, red] (11 + 6,-2.05) -- (11 + 7,-2.05);
        \draw [-latex, red] (11 + 7,-1.95) -- (11 + 8,-3);
        \draw [-latex, red] (11 + 8,-3)    -- (11 + 9,-3);

        \draw [-latex, red] (11 + 0,-3)    -- (11 + 5,-3);
        \draw [-latex, red] (11 + 5,-3)    -- (11 + 6,-2 + 0.05);
        \draw [-latex, blue] (11 + 6,-1.95) -- (11 + 7,-1.95);
        \draw [-latex, blue] (11 + 7,-1.95) -- (11 + 8,-1);
        \draw [-latex, blue] (11 + 8,-1)    -- (11 + 9,-1);

        \draw [blue] (11 + 0.4,-1.4) node {$f$};

        % g'
        \draw [-latex, blue] (11 + 0,1)    -- (11 + 1,1);
        \draw [-latex, blue] (11 + 1,1)    -- (11 + 2,0.05);
        \draw [-latex, blue] (11 + 2,0.05) -- (11 + 3,0.05);
        \draw [-latex, blue] (11 + 3,0.05) -- (11 + 4,1);
        \draw [-latex, blue] (11 + 4,1)    -- (11 + 9,1);
        \draw [red] (11 + 0.4,1.4) node {$g$};

        \draw [black] (11 + 2.5,-0.5) node {$\bullet$};
        \draw [black] (11 + 6.5,-2.5) node {$\ast$};
    \end {tikzpicture}
    \end {figure}
\end {proof}

\bigskip

The next lemma is a converse to the previous one, so in fact the conditions given above completely characterize
the set of terms in the expansion of $\left<t_{ik}, t_{j\ell}\right>$ with a given term $f_\bullet g_\bullet$.

\bigskip

\begin {lem}
    Let $(f,g) \in A_{ik} \times A_{j\ell}$ and $\bullet \in f \cap g$ such that $\varepsilon_\bullet(f,g) \neq 0$. \\[1.2ex]
    \begin {tabular}{cp{5in}}
        $(a)$ & If there exists an $\ast \in f \cap g$ such that the segments of $f$ and $g$ between $\bullet$ and $\ast$
                (or between $\ast$ and $\bullet$) cross the cut and rim the same number of times, then swapping those
                segments gives $(f',g') \in A_{ik} \times A_{j\ell}$ such that $f'_\ast g'_\ast = f_\bullet g_\bullet$. \\[1.2ex]
        $(b)$ & If there exists an $\ast \in f \cap f$ such that $g$ and the subloop of $f$ based at $\ast$ cross
                the cut and rim the same number of times, then swapping $g$ and this subloop gives $(f',g') \in A_{ik} \times A_{j\ell}$
                such that $f'_\ast g'_\ast = f_\bullet g_\bullet$.
    \end {tabular}
\end {lem}

\bigskip

We will call the swaps from part $(a)$ of the previous two lemmas ``type I'' swaps, and the swaps from part $(b)$
will be called ``type II'' swaps. We are now ready to define our sign-reversing map, which we will denote $\sigma$,
which acts on the set of terms appearing in the sum. 

\bigskip

\begin {lem} \label{lem:sign_reversing}
    There exists a fixed-point-free permutation $\sigma$ of the non-zero terms appearing in the expansion of $\left<t_{ik},t_{j\ell}\right>$
    such that $\sigma(x) = -x$ for each term $x$.
\end {lem}
\begin {proof}
    Recall that $f$ is a term in $t_{ik}$ and $g$ a term
    in $t_{j\ell}$. Assume that $i \geq j$, so that $f$ goes around the torus more times ``horizontally''.
    Starting at $\bullet$, walk along $f$ until we reach the first admissible swap. If the first admissible swap
    is of type I, and occurs at $\ast$, then we will define the image of the term $\varepsilon_\bullet(f,g) f_\bullet g_\bullet$
    under $\sigma$ to be the term $\varepsilon_\ast(f',g') f'_\ast g'_\ast$, where $f'$ and $g'$ are obtained by
    performing the type I swap between $\bullet$ and $\ast$. In order for $\varepsilon_\ast(f',g') \neq 0$, it must
    be that $f$ and $g$ cross at $\ast$ (rather than touch). This guarantees that $\varepsilon_\bullet(f,g) f_\bullet g_\bullet$
    is not a fixed point of $\sigma$. To see that $\varepsilon_\ast(f',g') = -\varepsilon_\bullet(f,g)$, note
    that by the preceding lemmas, the segments of $f$ and $g$ from $\bullet$ to $\ast$ must cross the cut and rim
    the same number of times. This means that on the universal cover, the segments of $f$ and $g$ between $\bullet$
    and $\ast$ bound a contractible disc. Thus after performing the swap, $f'$ and $g'$ touch with opposite orientation --
    i.e. if $f$ was to the left, and $g$ to the right, then $f'$ is on the right, and $g'$ on the left. This can
    be seen in \textbf{Figure \ref{fig:type_I_swap}}.

    \medskip

    If, on the other hand, the first admissible swap is of type II, we consider two separate cases. Recall
    that a type II swap means $f$ intersects itself at $\ast$, so we may write $f_\ast = f'_\ast f''_\ast$,
    and that $f''_\ast \in A_{j\ell}$ (the same as $g$). First we consider the case that $f''_\ast$ does not
    intersect $g$. In this case, we define the action of $\sigma$ on $\varepsilon_\bullet(f,g) f_\bullet g_\bullet$ to be
    the term corresponding to performing the type II swap. Now, consider the second case, in which
    $f''_\ast$ does intersect $g$. In this case, let $\star$ be the first intersection point after $\ast$.
    Then we define the action of $\sigma$ on $\varepsilon_\bullet(f,g) f_\bullet g_\bullet$ to be
    the term corresponding to the type I swap from $\bullet$ to $\star$.
\end {proof}

\bigskip

\textbf{Theorem \ref{thm:involutivity}} now follows immediately as a corollary, since the lemma implies that
there are an even number of terms equivalent to any given $f_\bullet g_\bullet$, half with coefficient $+1$
and half with coefficient $-1$.

\bigskip

\section {Recovering the Lax Representation}

\bigskip

In \cite{bf15}, the authors gave a Lax representation of the Grassmann pentagram map, which gives a family
of invariants. However, they did not formulate Liouville integrability, since there was no associated Poisson
structure. 
In this section, we explain how to apply the results of the present paper to recover the invariants
from \cite{bf15}, and additionally give a Poisson structure in which these invariants commute.

\bigskip

As described earlier, we may lift the vertices of a twisted Grassmann $n$-gon to $3N$-by-$N$ matrices $V_i$,
and we get a relation of the form
\[ V_{i+3} = V_i A_i + V_{i+1} B_i + V_{i+2} C_i \]
As before, we define the $3N$-by-$3N$ matrices $\mathbf{V}_i = (V_i V_{i+1} V_{i+2})$. 
Then the formula above can be re-phrased by
saying $\mathbf{V}_{i+1} = \mathbf{V}_i L_i$, where $L_i$ is the block matrix

\[ L_i = 
    \left( \begin {array} {ccc}
        0 & 0 & A_i \\
        \mathrm{Id}_N & 0 & B_i \\
        0 & \mathrm{Id}_N & C_i
    \end {array} \right)
\]

Assuming the lift was twisted, then
the monodromy matrix of the polygon is related to the $L_i$ by 
\[ M \, \mathbf{V}_i = \mathbf{V}_i \, L_i L_{i+1} \cdots L_{i+n-1} \]
In that paper, the authors proved the following

\medskip

\begin {thm} \cite{bf15}
    The Grassmann pentagram map is invariant under the scaling $B_i \mapsto \lambda^{-1} B_i$, $C_i \mapsto \lambda C_i$.
\end {thm}

\medskip

In light of this theorem, we may define the modified matrix with parameter $\lambda$:
\[ L_i(\lambda) = 
    \left( \begin {array}{ccc}
        0 & 0 & A_i \\
        \mathrm{Id}_N & 0 & \lambda^{-1} B_i \\
        0 & \mathrm{Id}_N & \lambda C_i
    \end {array} \right)
\]

In \cite{bf15}, the authors proved that the conjugacy class of $L(\lambda) := L_1(\lambda) \cdots L_n(\lambda)$ is preserved under the
pentagram map, so the spectral invariants of $L(\lambda)$ are invariants of the map. It is not hard to see
that the pentagram map is also invariant under the scaling $A_i \mapsto \lambda A_i$, $B_i \mapsto \lambda B_i$.
We will call the corresponding matrix $\widetilde{L}_i(\lambda)$:
\[ \widetilde{L}_i(\lambda) = 
    \left( \begin {array}{ccc}
        0 & 0 & \lambda A_i \\
        \mathrm{Id}_N & 0 & \lambda B_i \\
        0 & \mathrm{Id}_N &  C_i
    \end {array} \right)
\]

\bigskip

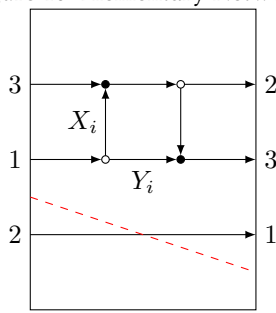
\begin {figure}[h!]
    \centering
    \caption{Elementary Networks}
    \label{fig:elementary}
    \begin {tikzpicture}
        \draw (0,0) -- (3,0) -- (3,4) -- (0,4) -- cycle;

        \draw [fill=black] (1,3) circle (0.05cm);
        \draw [fill=white] (2,3) circle (0.05cm);
        \draw [fill=white] (1,2) circle (0.05cm);
        \draw [fill=black] (2,2) circle (0.05cm);

        \draw [-latex] (0,3) -- (0.95,3);
        \draw [-latex] (1.05,3) -- (1.95,3);
        \draw [-latex] (2.05,3) -- (3,3);
        \draw [-latex] (0,2) -- (0.95,2);
        \draw [-latex] (1.05,2) -- (1.95,2);
        \draw [-latex] (2.05,2) -- (3,2);
        \draw [-latex] (0,1) -- (3,1);
        \draw [-latex] (1,2.05) -- (1,2.95);
        \draw [-latex] (2,2.95) -- (2,2.05);

        \draw [dashed, red] (0,1.5) -- (3,0.5);

        \draw (0.7,2.5) node {$X_i$};
        \draw (1.5,1.7) node {$Y_i$};

        \draw (3.2,1) node {$1$};
        \draw (3.2,2) node {$3$};
        \draw (3.2,3) node {$2$};

        \draw (-0.2,1) node {$2$};
        \draw (-0.2,2) node {$1$};
        \draw (-0.2,3) node {$3$};
    \end {tikzpicture}
\end {figure}

In \cite{gsvt_16}, the authors present a very similar Lax representation for the pentagram map on $\Bbb{P}^2$,
and explicitly describe the connection with the boundary measurement matrix used in the combinatorial proof
of integrability. We now mimic this approach to connect the earlier results from this paper to the Lax
representation described above from \cite{bf15}. In order to more easily generalize this approach,
we change some of our notations and conventions from earlier to more closely resemble those from \cite{gsvt_16}.
In particular, we re-index the $Y_i$'s by 1, so that our relation reads
\[ V_{i+3} = V_i Y_{i-1} + V_{i+1} X_i + V_{i+2} \]
Also, recall we use an oriented curve called the ``cut'' on the cylinder to determine the powers of $\lambda$
in the boundary measurements. Earlier in the paper, we chose the convention that when we draw the cylinder
as a rectangle (identifying the top/bottom edges), we take the cut to be the the top/bottom edge. We change this
convention now so that the cut goes diagonally down and to the right, crossing each edge that connects two white vertices.
With this convention, the quiver can be realized as the concatenation of the elementary networks pictured in
\textbf{Figure \ref{fig:elementary}}.

The boundary measurement matrix of the $i^\mathrm{th}$ elementary network is given by
\[ \mathscr{B}_i(\lambda) = \left( \begin {array}{ccc}
        0 & X_i & X_i + Y_i \\
        \lambda & 0 & 0 \\
        0 & 1 & 1
    \end {array} \right)
\]
The boundary measurement matrix of the entire network $Q_n$ is then the product $\mathscr{B}(\lambda) := \mathscr{B}_1(\lambda) \cdots \mathscr{B}_n(\lambda)$. 
Although our matrices have entries in the non-commutative ring $\mathscr{L}_\bullet[\lambda^\pm]$,
the following result from \cite{gsvt_16} is true in this more general context:

\medskip

\begin {prop} \cite{gsvt_16}
    $\mathscr{B}_i(\lambda) = \mathscr{A}_i(\lambda) \widetilde{L}_{i+1}(\lambda) \mathscr{A}_{i+1}(\lambda)^{-1}$, where the matrices $\mathscr{A}_i(\lambda)$ are given by
    \[ \mathscr{A}_i(\lambda) := \left( \begin{array}{ccc} \lambda^{-1} & 0 & X_i \\ 0 & 1 & 0 \\ 0 & 0 & 1 \end{array} \right) \]
    In particular, $\mathscr{B}(\lambda)$ is conjugate to $\widetilde{L}(\lambda) := \widetilde{L}_1(\lambda) \cdots \widetilde{L}_n(\lambda)$.
\end {prop}

\bigskip

Therefore, the spectral invariants of $\mathscr{B}(\lambda)$ are the same as those of $\widetilde{L}(\lambda)$.
The coefficients of $\lambda$ in the expansions of $\mathrm{tr}(\mathscr{B}(\lambda)^k)$ can be written as non-commutative polynomials
in the $X_i$ and $Y_i$. Thus the (traces of the) spectral invariants can be interpreted as functions
on $\mathcal{GP}_{n,N}$, after identifying it with $\mathrm{GL}_N^{2n+1}/\mathrm{Ad} \, \mathrm{GL}_N$, as in \textbf{Theorem \ref{thm:moduli_space_identification}}. The Poisson bracket is
given by taking the trace of the induced bracket on $\mathscr{L}_\bullet^\natural$:
\[ \{\mathrm{tr}(X_i), \, \mathrm{tr}(Y_j) \} = \mathrm{tr} \left<X_i,Y_j\right> \]
As before, we let $t_{ij}$ denote the components of the spectral invariants, so that $\mathrm{tr}(\mathscr{B}(\lambda)^j) = \sum_i t_{ij} \lambda^i$.
The discussion above tells us that $\mathrm{tr}(t_{ij})$ may be interpreted as functions on $\mathcal{GP}_{n,N}$.
Although in the present paper we chose a different normalization convention (namely the $X_i$, $Y_i$, $Z$ matrices)
than Mari Beffa and Felipe, these invariant functions are still essentially the same as those mentioned in \cite{bf15}.
\textbf{Theorem \ref{thm:involutivity}} then implies that these functions form an involutive family with respect to
the induced Poisson structure.

\bigskip

\section {Conclusions and Suggestions for Future Work}

We have shown in \textbf{Theorem \ref{thm:invariants}} and \textbf{Theorem \ref{thm:involutivity}} that the coefficients
$t_{ij}$ of the traces of powers of the boundary measurement matrix are non-commutative invariants of the pentagram map,
and that they form an involutive family under the induced Lie bracket $\left<-,-\right>$ on $\mathscr{L}^\natural$.
In this sense, we have established a form of non-commutative integrability for the Grassmann pentagram map, in terms
of the $X_i,Y_i$ variables/matrices. Furthermore, this induces a Poisson structure on the moduli space $\mathcal{GP}_{n,N}$
in which Mari-Beffa and Felipe's invariants Poisson-commute.

\medskip

However, in the usual commutative setting, the definition of Liouville integrability requires not only a Poisson-commuting family
of invariants, but also that these invariants are independent, and that they are a maximal family of such independent
commuting invariants. These last two aspects --- independence and maximality --- were not addressed in the present paper.
It would be nice to have a proof that among this infinite family of commuting invariants, one can find a maximal independent subset,
so that we have a Liouville integrable system.

\medskip

The proof of \textbf{Theorem \ref{thm:involutivity}} is very combinatorial. It would be interesting to give an alternate
proof using $R$-matrices, which would more closely resemble the method of proof in the classical case \cite{gsvt_16}.
This would involve formulating an $R$-matrix double bracket on the noncommutative space of matrices, and formulating
some result saying that the spectral invariants are in involution with respect to the induced $H_0$-Poisson structure.

\section {Acknowledgments}

\bigskip

The author would like to give thanks to Michael Shapiro for suggesting the problem, and for helpful feedback
and conversations; and to Semeon Artamonov and Leonid Chekhov for helpful discussions.
The author was partially supported by the NSF grant DMS-1702115.

\vfill

\bibliographystyle {alpha}
\bibliography {paper}

\begin{thebibliography}{GGRW02}

\bibitem[Art18]{artamonov_18}
Semen Artamonov.
\newblock {\em Generalized quasi Poisson structures and noncommutative
  integrable systems}.
\newblock PhD thesis, Rutgers University-School of Graduate Studies, 2018.

\bibitem[Ber08]{bergh_08}
Michel Van~Den Bergh.
\newblock Double poisson algebras.
\newblock {\em Transactions of the American Mathematical Society},
  360(11):5711--5769, 2008.

\bibitem[CB11]{crawley-boevey_11}
William Crawley-Boevey.
\newblock Poisson structures on moduli spaces of representations.
\newblock {\em Journal of Algebra}, 325(1):205--215, 2011.

\bibitem[Coh77]{cohn77}
Paul~Moritz Cohn.
\newblock {\em Skew field constructions}, volume~27.
\newblock CUP Archive, 1977.

\bibitem[FG06]{fock_goncharov}
Vladimir Fock and Alexander Goncharov.
\newblock Moduli spaces of local systems and higher teichm{\"u}ller theory.
\newblock {\em Publications Math{\'e}matiques de l'Institut des Hautes
  {\'E}tudes Scientifiques}, 103(1):1--211, 2006.

\bibitem[FMB15]{bf15}
Ra{\'u}l Felipe and Gloria Mari~Beffa.
\newblock The pentagram map on grassmannians.
\newblock {\em arXiv preprint arXiv:1507.04765}, 2015.

\bibitem[FZ02]{fz1}
Sergey Fomin and Andrei Zelevinsky.
\newblock Cluster algebras i: foundations.
\newblock {\em Journal of the American Mathematical Society}, 15(2):497--529,
  2002.

\bibitem[FZ07]{fz4}
Sergey Fomin and Andrei Zelevinsky.
\newblock Cluster algebras iv: coefficients.
\newblock {\em Compositio Mathematica}, 143(1):112--164, 2007.

\bibitem[GGRW02]{grw02}
Israel Gelfand, Sergei Gelfand, Vladimir Retakh, and Robert Wilson.
\newblock Quasideterminants.
\newblock {\em arXiv preprint math/0208146}, 2002.

\bibitem[Gli11]{glick11}
Max Glick.
\newblock The pentagram map and y-patterns.
\newblock {\em Advances in Mathematics}, 227(2):1019--1045, 2011.

\bibitem[Gol84]{goldman_84}
William~M Goldman.
\newblock The symplectic nature of fundamental groups of surfaces.
\newblock {\em Advances in Mathematics}, 54(2):200--225, 1984.

\bibitem[Gol86]{goldman_86}
William~M Goldman.
\newblock Invariant functions on lie groups and hamiltonian flows of surface
  group representations.
\newblock {\em Inventiones mathematicae}, 85(2):263--302, 1986.

\bibitem[GSTV16]{gsvt_16}
Michael Gekhtman, Michael Shapiro, Serge Tabachnikov, and Alek Vainshtein.
\newblock Integrable cluster dynamics of directed networks and pentagram maps.
\newblock {\em Advances in Mathematics}, 300:390--450, 2016.

\bibitem[GSV09]{gsv_09}
Michael Gekhtman, Michael Shapiro, and Alek Vainshtein.
\newblock Poisson geometry of directed networks in a disk.
\newblock {\em Selecta Mathematica, New Series}, 15(1):61--103, 2009.

\bibitem[GSV10a]{gsv_book}
Michael Gekhtman, Michael Shapiro, and Alek Vainshtein.
\newblock {\em Cluster algebras and Poisson geometry}.
\newblock Number 167. American Mathematical Soc., 2010.

\bibitem[GSV10b]{gsv_10}
Michael Gekhtman, Michael Shapiro, and Alek Vainshtein.
\newblock Poisson geometry of directed networks in an annulus.
\newblock {\em arXiv preprint arXiv:0901.0020}, 2010.

\bibitem[Izo18]{izosimov_18}
Anton Izosimov.
\newblock Pentagram maps and refactorization in poisson-lie groups.
\newblock {\em arXiv preprint arXiv:1803.00726}, 2018.

\bibitem[Lew74]{lewin74}
Jacques Lewin.
\newblock Fields of fractions for group algebras of free groups.
\newblock {\em Transactions of the American Mathematical Society},
  192:339--346, 1974.

\bibitem[MKS76]{magnus_76}
Wilhelm Magnus, Abraham Karrass, and Donald Solitar.
\newblock {\em Combinatorial group theory}.
\newblock Dover Publications, 1976.

\bibitem[MT12]{turaev_12}
Gwenael Massuyeau and Vladimir Turaev.
\newblock Quasi-poisson structures on representation spaces of surfaces.
\newblock {\em International Mathematics Research Notices}, 2014(1):1--64,
  2012.

\bibitem[OST10]{ost10}
Valentin Ovsienko, Richard Schwartz, and Serge Tabachnikov.
\newblock The pentagram map: a discrete integrable system.
\newblock {\em Communications in Mathematical Physics}, 299(2):409--446, 2010.

\bibitem[Pos06]{postnikov_06}
Alexander Postnikov.
\newblock Total positivity, grassmannians, and networks.
\newblock {\em arXiv preprint math/0609764}, 2006.

\bibitem[Pro76]{procesi_76}
Claudio Procesi.
\newblock The invariant theory of n $\times$ n matrices.
\newblock {\em Advances in Mathematics}, 19(3):306 -- 381, 1976.

\bibitem[Sch92]{schwartz92}
Richard Schwartz.
\newblock The pentagram map.
\newblock {\em Experimental Mathematics}, 1(1):71--81, 1992.

\bibitem[Sib68]{sibirskii_68}
KS~Sibirskii.
\newblock Algebraic invariants for a set of matrices.
\newblock {\em Siberian Mathematical Journal}, 9(1):115--124, 1968.

\end{thebibliography}

\end {document}